\def\VERSION{26.7.2024}
\def\WHO{nbd} 
\def\users{us}    
\def\users{world} 
\documentclass[11pt]{article} 

\usepackage{color}
\usepackage{amsmath}
\usepackage{amsthm}
\usepackage{amssymb}
\usepackage{epsfig}
\usepackage{psfrag}
\usepackage{graphicx}
\usepackage{textcomp}
\usepackage{bm}
\usepackage{pgf}
\numberwithin{equation}{section}
\usepackage{upgreek} 
\newtheorem{theorem}{Theorem}[section]

\newtheorem{definition}[theorem]{Definition}
\newtheorem{example}[theorem]{Example}
\newtheorem{proposition}[theorem]{Proposition}

\newtheorem{remark}[theorem]{Remark}

\usepackage{mathrsfs,cite}
\usepackage{ifthen}
\usepackage{ulem}
\usepackage[colorlinks=true, pdfstartview=FitV, linkcolor=blue, 
            citecolor=blue, urlcolor=blue]{hyperref}
\usepackage{cancel}\ifthenelse{\equal{\users}{world}}
{
\newcommand{\REM}[1]{}

	\newcommand{\DELETE}[1]{}

        \newcommand{\COMMENT}[1]{}
        \newcommand{\TCOMMENT}[1]{}
    \newcommand{\MARGINOTE}[1]{}
}	
{
\usepackage[notcite,notref,color]{showkeys}

\definecolor{brown}{rgb}{0.6,0.2,0.2}
\newcommand{\REM}[1]{\marginpar{\bfseries\tiny{\color{blue}#1}}}

 \newcommand{\COMMENT}[1]{{\color{blue}\uuline{#1}\color{black}}} 
 \newcommand{\DELETE}[1]{{\color{brown}\cancel{#1}\color{black}}}

 \newcommand{\TCOMMENT}[1]{{\color{blue}{ #1}}}
\newcommand{\MARGINOTE}[1]{\marginpar{\color{red}\tiny\texttt{#1}}}
\usepackage{fancyhdr}
\pagestyle{fancy}
\headheight=28pt\headwidth=17cm
\definecolor{gray}{gray}{0.5}
\newcount\hour \newcount\minute
\hour=\time
\divide \hour by 60
\minute=\time
\loop \ifnum \minute > 59 \advance \minute by -60 \repeat

\rhead{\color{gray}time discretization Eulerian\\
by T.Roub\'\i\v cek}
\chead{}
\lhead{Version \VERSION, file: \jobname.tex, \underline{\Large\color{red}\WHO's working on} \\
compiled:
\number\day.\number\month.\number\year\ at
\the\hour:\ifnum\minute<10 0\fi\the\minute\ h }
}

\newcommand{\R}{\mathbb{R}}
\newcommand{\N}{\mathbb{N}}
\newcommand{\bbI}{\mathbb{I}}

\newcommand{\bbD}{\mathbb{D}}

\newcommand\DT[1]{\mathchoice
                 {{\buildrel{\hspace*{.1em}\text{\LARGE.}}\over{#1}}}
                 {{\buildrel{\hspace*{.1em}\text{\LARGE.}}\over{#1}}}
                 {{\buildrel{\hspace*{.1em}\text{\Large.}}\over{#1}}}
                 {{\buildrel{\hspace*{.1em}\text{\large.}}\over{#1}}}}

\newcommand{\lineunder}[2]{\LU{\begin{array}[t]{c}\underbrace{#1}\vspace*{.5em}\end{array}}{\mbox{\footnotesize\rm #2}}}
\newcommand{\linesunder}[3]{\LSU{\begin{array}[t]{c}\underbrace{#1}\vspace*{.5em}\end{array}}{\mbox{\footnotesize\rm #2}}{\mbox{\footnotesize\rm#3}}}
\newcommand{\LU}[2]{\begin{array}[t]{c}#1\vspace*{-1em}\\_{#2}\end{array}}
\newcommand{\LSU}[3]{\begin{array}[t]{c}#1\vspace*{-1em}\\_{#2}\vspace*{-.5em}\\_{#3}\end{array}}
\renewcommand{\d}{{\rm d}}

\newcommand{\NablaS}{\Nabla_{\scriptscriptstyle\textrm{\hspace*{-.3em}S}}^{}}
\newcommand{\divS}{\mathrm{div}_{\scriptscriptstyle\textrm{\hspace*{-.1em}S}}^{}}
\newcommand{\eq}[1]{(\ref{#1})}
\newcommand{\Cdot}{\hspace{-.1em}\cdot\hspace{-.1em}}
\newcommand{\Colon}{\hspace{-.15em}:\hspace{-.15em}}
\newcommand{\UUU}[3]{\begin{array}[b]{c}\vspace*{-1.2mm}_{\text{\scriptsize{#2}}}\vspace*{-.6mm}\\[-.0em]_{\text{\scriptsize{#3}}}\vspace*{.6mm}\\[-.5em]#1\end{array}}

\def\vv{{\bm v}}
\def\pp{{\bm p}}
\def\uu{{\bm u}}

\def\AA{{\bm A}}

\def\xx{{\bm x}}
\def\yy{{\bm y}}
\def\nn{{\bm n}}

\newcommand{\DD}{\bm D}
\newcommand{\TT}{\bm T}
\newcommand{\Se}{\bm S}
\newcommand{\strain}{{\boldsymbol\varepsilon}}
\newcommand{\GRAVITY}{\bm g}

\def\vvk{\vv_\etau^k}
\def\vvkk{\vv_\etau^{k-1}}
\def\overlineEetau{\hspace*{.2em}\overline{\hspace*{-.2em}\bm E}_\etau^{}}
\def\overlineFetau{\hspace*{.2em}\overline{\hspace*{-.2em}\bm F}_{\!\etau}^{}}
\def\overlineJetau{\hspace*{.2em}\overline{\hspace*{-.2em}J}_\etau^{}}

\def\overlineLptau{\hspace*{.3em}\overline{\hspace*{-.3em}\Lp\hspace*{-.1em}}_{\,\etau^{}}}

\def\overlinevvtau{\hspace*{.15em}\overline{\hspace*{-.15em}\vv}_{\etau}^{}}

\def\FF{{\bm F}}
\def\HH{{\bm H}}
\def\MM{{\bm M}}
\def\WW{{\bm W}}
\def\Lp{{\bm \varPi}}

\def\HYPER{\mu}
\def\CUT{{\mathcal{K}}}

\def\rhoMAX{{\rho_{\rm max}^{}}}

\newcommand\ZJ[1]{\mathchoice
                 {{\buildrel{\hspace*{.1em}{_{\,\boldsymbol\circ}}}\over{#1}}}
                 {{\buildrel{\hspace*{.1em}{_{\,\boldsymbol\circ}}}\over{#1}}}
                 {{\buildrel{\hspace*{.1em}{\boldsymbol\circ}}\over{#1}}}
                 {{\buildrel{\hspace*{.1em}{\boldsymbol\circ}}\over{#1}}}}
\newcommand{\barOmega}{\,\overline{\!\varOmega}}
\newcommand{\nablaS}{\nabla_{\scriptscriptstyle\textrm{\hspace*{-.3em}S}}^{}}

\newcommand{\Nabla}{\nabla}
\newcommand{\Rsym}{\mathbb R^{3\times3}_{\rm sym}}

\def\Vdots{\!\mbox{\setlength{\unitlength}{1em}
\begin{picture}(0,0)
\put(-.07,0){.}
\put(-.07,.3){.}
\put(-.07,.6){.}
\end{picture}
}
}

\newcommand{\wt}[1]{\mathchoice{\hspace*{-.09em}\text{\large$\hspace*{.09em}\tilde{\text{\normalsize$#1$}}\hspace*{.05em}$}\hspace*{-.05em}}
{\hspace*{-.09em}\text{\large$\hspace*{.09em}\tilde{\text{\normalsize$#1$}}\hspace*{.05em}$}\hspace*{-.05em}}
{\text{\normalsize$\hspace*{.08em}\tilde{\text{\scriptsize$#1$}}\hspace*{.06em}$}}
{\text{\small$\tilde{\text{\tiny$#1$}}$}}}

\newcounter{myfigure}
\newenvironment{my-picture}[3]{\refstepcounter{myfigure}\label{#3}\setlength{\unitlength}{1em}\begin{picture}(#1,#2)}{\end{picture}}

\newcommand\DELETEDELETE[1]{}

\newcommand\pdt[1]{\frac{\partial{#1}}{\partial t}}
\newcommand\Ee{{\bm E}}              
\newcommand\Ep{{\bm P}}              

\begin{document}

\def\EPS{\varepsilon}
\def\DELTA{\delta}
\def\etau{{\EPS\DELTA\tau}}
\def\EEps{{\EPS\DELTA}}

\def\TTtauk{\TT^k_{\!\etau}}
 \def\Eek{\Ee_\etau^k}
 \def\Eetau{\Ee_\etau^{}}
\def\overlineDetau{\hspace*{.2em}\overline{\hspace*{-.2em}\bm D}_\etau^{}}
\def\overlineTetau{\hspace*{.2em}\overline{\hspace*{-.1em}\bm T}_\etau^{}}
\def\overlineLpetau{\hspace*{.2em}\overline{\hspace*{-.3em}\bm\varPi\hspace*{-.1em}}_\etau^{}}
\def\overlineMetau{\hspace*{.4em}\overline{\hspace*{-.4em}\bm M\hspace*{-.1em}}_\etau^{}}
\def\overlineFetau{\hspace*{.2em}\overline{\hspace*{-.2em}\bm F}_{\!\etau}^{}}
\def\overlineHetau{\hspace*{.2em}\overline{\hspace*{-.2em}\bm H}_\etau^{}}
\def\overlineEetau{\hspace*{.2em}\overline{\hspace*{-.2em}\bm E}_\etau^{}}
\def\overlinevvtau{\hspace*{.1em}\overline{\hspace*{-.1em}\vv}_{\etau}^{}}
\def\overlineppetau{\hspace*{.15em}\overline{\hspace*{-.15em}\pp}_\etau^{}}
\def\overlinerhoetau{\hspace*{.15em}\overline{\hspace*{-.15em}\varrho}_\etau^{}}

\def\Ke{K_\text{\sc e}^{}}
\def\Ge{G_\text{\sc e}^{}}
\def\Kv{K_\text{\sc v}^{}}
\def\Gv{G_\text{\sc v}^{}}
\def\Gm{G_\text{\sc m}^{}}

\allowdisplaybreaks

\noindent{\LARGE\bf Time discretization in visco-elastodynamics
\\[.2em]at large displacements and strains
\\[.2em]in the Eulerian frame.\footnote{{\tt ORCID}: 0000-0002-0651-5959.}
}

\bigskip\bigskip

\noindent{\large\sc Tom\'{a}\v{s} Roub\'\i\v{c}ek}\\
{\it Mathematical Institute, Charles University, \\Sokolovsk\'a 83,
CZ--186~75~Praha~8,  Czech Republic
}\\and\\
{\it Institute of Thermomechanics, Czech Academy of Sciences,\\Dolej\v skova~5,
CZ--182~08 Praha 8, Czech Republic
}

\bigskip\bigskip

\begin{center}\begin{minipage}[t]{14.5cm}

{\small

\noindent{\bfseries Abstract.}
The fully-implicit time discretization (i.e.\ the 
backward Euler formula) is applied to compressible nonlinear
dynamical models of viscoelastic solids in the Eulerian
description, i.e.\ in the actual deforming configuration.
The Kelvin-Voigt rheology or also, in the deviatoric part, the
Jeffreys rheology are considered. Both a linearized convective
model at large displacements with a convex stored energy
and the fully nonlinear large strain variant with a (possibly
generalized) polyconvex stored energy are considered.
The time-discrete suitably regularized schemes are devised for
both cases. The numerical stability and, considering the multipolar
2nd-grade viscosity, also convergence towards weak solutions are proved,
exploiting the convexity of the kinetic energy when written in terms
of linear momentum instead of velocity. In the fully nonlinear case,
the examples of neo-Hookean and Mooney–Rivlin materials are presented.
A comparison with models of viscoelastic barotropic fluids is also made.

\medskip

\noindent{\it Keywords}: finitely-strained solids,
visco-elastodynamics, Kelvin-Voigt rheology, anti-Zener rheology, Euler description,
backward Euler time discretization, Rothe method, weak solutions.

\medskip

\noindent{\small{\it AMS Subject Classification}:
35Q74, 
65M99, 
74A30, 
74B20, 
74H20. 
74S99. 
}

} 
\end{minipage}
\end{center}

\bigskip

\section{Introduction}

The models of {\it visco-elastodynamics} in continuum mechanics {\it at finite}
(also called {\it large}) {\it strains} lead to heavily nonlinear systems of
evolution partial differential
equations. In solid mechanics, the Lagrangian approach (using a referential
configuration) is mostly used, but sometimes the {\it Eulerian approach}
(using the actual deforming configuration) has specific advantages. In particular,
sometimes no referential configuration has a justified sense (e.g.\ in
geophysical models on long time scales of millions years) and the
interaction with external spatial fields (such as gravitational or electromagnetical)
or fluid-solid interaction is more straightforward in the Eulerian frame.
Also, the Stokes viscosity is more straightforward in the Eulerian frame than
in the Lagrangian frame where frame-indifferency needs special
nonlinearities and special analytical effort. Time discretization of
visco-elastodynamics in the Lagrangian frame seems problematic at large strains
due to the inevitable nonconvexity which cannot be compensated by viscosity,
cf.\ \cite[Exercises~9.2.9 and 9.2.15]{KruRou19MMCM}, while in the
Eulerian frame the time discretization is duable, which is the main
goal of this article, shown in Section~\ref{sec-finite-strain} below.
An added value of the description in the actual Eulerian frame is the revelation
of the actual physics, in particular the actual (i.e.\ Cauchy) stress instead of the
rather fictitious Piola-Kirchhoff stress. Both approaches can benefit from involving
some higher-order gradients either in the conservative or in the dissipative
part of the models, respectively, which facilitates the rigorous analysis and
allows for various dispersion of velocity of propagation of elastic waves, cf.\
\cite{Roub24SGTL}.
These higher (here 2nd order) gradients lead to the concept of (here 2nd-grade)
{\it nonsimple media}, as often occurs in literature since the works by
R.A.\,Toupin \cite{Toup62EMCS} and R.D.\,Mindlin \cite{Mind64MSLE}. In the dissipative
part as used in this paper, it was also devised by E.\,Fried and M.\,Gurtin
\cite{FriGur06TBBC} and earlier, even more generally and nonlinearly, as
{\it multipolar fluids} by J.\,Ne\v cas at al.\
\cite{BeBlNe92PBMV,NeNoSi89GSIC,NecRuz92GSIV}. 

Due to the mentioned nonlinear character of convective Eulerian models,
the Faedo-Galerkin space discretization is most commonly used, although
various nonlinear tests needed in energy methods make is not straightforward
and quite cumbersome, cf.\ \cite{Roub22VELS}.
For this reason, the time discretization (Rothe's) method seems to be a
good alternative, although the mentioned nonlinearities make
it not easy, either. Such {\it implicit time discretization} has sometimes
been used for the compressible fluid dynamics, mostly merely computationally
while for analytically rigorous treatment we refer to 
\cite{GaMaNo19EEIM,FHMN17EENM,FeKaPo16MTCV,FLMS21NACF,Kar13CFEM,Zato12ASCN}.
The novelty of this paper is the usage of the implicit time discretization
for compressible solid mechanics. Also, as a consequence of the
multipolar 2nd-grade viscosity, the convergence is proved
in 3D under suitable qualification of data without additional assumptions 
on a solution itself, in contrast to e.g.\ \cite{FLMS21NACF} where
strict positivity of mass density has to be supposed or to
\cite{Zato12ASCN} restricted to 2D situations, cf.\ also
Remark~\ref{rem-barotropic-fluids} below.

For completeness, let us mention the implicit (or semi-implicit) time
discretization in the (much simpler) incompressible cases e.g.\ in
\cite{ISTT11IFEM,NoSaTo14MNSE}.

In Section~\ref{sec-linearized} we investigate the linearized
convective model involving small strain tensor at large deformations.
This is a very popular approach especially in geophysical modelling
or also mechanical-engineering modelling, where the objective strain
rate due to the Zaremba-Jaumann time derivative is most often used.
The essential point for the time discretization is the convexity of
the kinetic energy when expressed in terms of linear momentum instead
of velocity, which allows us to present the main analytical aspects first
on a relatively simpler but still well aplicable model. Then, 
Section~\ref{sec-finite-strain} treats the fully
nonlinear models which also ultimately suffer from the non-convexity of the
frame-invariant stored energy. The time discretization then requires
some additional evolution equations for quantities entering into (possibly
generalized) polyconvexity and some additional analytical aspects.

For readers' convenience, let us summarize the basic notation used in what
follows:
\begin{center}
\fbox{
\begin{minipage}[t]{17em}\small\smallskip
$\yy$ deformation,\\
  $\vv$ velocity,\\
$\varrho$ mass density,\\
$\pp=\varrho\vv$ the linear momentum,\\
$\FF$ elastic distortion,\\
$\Ee$ small strain,\\
$\GRAVITY$ gravity acceleration,\\
$\Ke,\Ge$ elastic bulk and shear moduli,\\
$\Kv,\Gv$ viscosity bulk and shear moduli,\\
$\Gm$ Maxwellian viscosity modulus,\\
$\HYPER$ the hyper-viscosity coefficient,\\
$I=[0,T]$ a time interval, $T>0$,\\
$\DELTA$, $\EPS$ regularizing perameters,\\
$\R_{\rm sym}^{3\times3}$ set of symmetric matrices,\\
$\varphi:\R^{3\times 3}\to\R$ stored energy (actual),\\[.1em]
$(\cdot)'$ (partial) derivative of a function,\\[.1em]
$(\cdot)\!\DT{^{\,}}$ convective time derivative,
\end{minipage}
\begin{minipage}[t]{21em}\small\smallskip
$\TT$ Cauchy stress,\\
$\DD$ dissipative stress,\\
$\MM$ Mandel's stress,\\
$\Lp$ inelastic distortion rate,\\
$\FF_{\rm tot}=\nabla\yy$ deformation gradient,\\
${\rm Cof}\FF=J\FF^{-\top}\!=J(\FF^{\top})^{-1}\!$ = cofactor of $\FF$,\\
$J=\det\FF$ Jacobian\,=\,determinant of $\FF$,\\
tr$(\cdot)$ trace of a matrix,\\
dev$(\cdot)$ deviatoric part of a matrix,\\
$R=1/\Gm$ reciprocal Maxwellian creep modulus,\\
$\nu$ coefficient for inelastic distortion rate gradient,\\
$\bbI$ the identity matrix,\\
$\R_{\rm dev}^{3\times3}=\{A\in\R_{\rm sym}^{d\times d};\ {\rm tr}A=0\}$,\\
${\rm GL}_3^+=\{A\in\R^{3\times3};\ \det A>0\}$,\\
$\upvarphi:\R^{3\times3}\to\R$ stored energy (referential),\\
$\tau>0$ a time step for discretization,\\
$(\cdot)\!\ZJ{^{\,}}$ Zaremba-Jaumann corotational time derivative.
\smallskip \end{minipage}
}\end{center}

\vspace{-.9em}

\begin{center}
{\small\sl Table\,1.\ }
{\small
Summary of the basic notation used. 
}
\end{center}

We will use the standard notation concerning the Lebesgue and the Sobolev
spaces of functions on the Lipschitz bounded domain $\varOmega\subset\R^3$,
namely $L^p(\varOmega;\R^n)$ for Lebesgue measurable $\R^n$-valued functions
$\varOmega\to\R^n$ whose Euclidean norm is integrable with $p$-power, and
$W^{k,p}(\varOmega;\R^n)$ for functions from $L^p(\varOmega;\R^n)$ whose
all derivative up to the order $k$ have their Euclidean norm integrable with
$p$-power. We also write briefly $H^k=W^{k,2}$.
We have the embedding $H^1(\varOmega)\subset L^6(\varOmega)$. Moreover, for a Banach
space $X$ and for $I=[0,T]$, we will use the notation $L^p(I;X)$ for the Bochner
space of Bochner measurable functions $I\to X$ whose norm is in $L^p(I)$, 
and $H^1(I;X)$ for functions $I\to X$ whose distributional derivative is in $L^2(I;X)$.
Occasionally, we will use $L_{\rm w*}^p(I;X)$ the space of weakly* measurable
mappings $I\to X$ if $X$ has a predual, i.e.\ there is $X'$ such that $X=(X')^*$
where $(\cdot)^*$ denotes the dual space.
Note also the embedding $L^\infty(I;L^2(\varOmega))\,\cap\,L^p(I;W^{1,p}(\varOmega))
\subset L^{5p/3}(I{\times}\varOmega)$.
The space of continuous functions on the closure $\barOmega$ of
$\varOmega$ will be denoted by $C(\barOmega)$.

\section{Linearized large-deformation convective model}\label{sec-linearized}

Most materials typically cannot withstand too much large elastic strains without
initiating inelastic processes (such as damage or creep or plastification).
Thus, the elastic strains will always remain rather
small (and the elastic distortion tensor is close to the identity $\bbI$),
and a usual small-strain linearization is well acceptable and
widely used in many applications. Anyhow, small strains do not preclude
large displacements, typically in fluids but also in solids, in the latter case
particularly when the Kelvin-Voigt model is combined with Maxwellian-type rheology in
the deviatoric part. This suggests the use of the Eulerian small-strain models combined
with a properly designed transport of the strain tensor, in addition to
the usual transport of the mass density. Of course, thermodynamic consistency
in terms of mass, momentum, and energy are concerned is the ultimate attribute that
any rational model must respect.

This compromise combination of small elastic strains with large deformations and
displacements requires appropriate formulations in a convected
coordinate system, in particular a proper choice and treatment
of objective rates. Objectivity here means that the time
derivatives do not depend on the evolving reference frame.
There are many possibilities used in the literature for different models.
It is reasonable to require that the tensor time derivatives, i.e.\ the tensor rates
(in particular the stress rate), are so-called identical and corotational,
meaning that the stress rate vanishes for all rigid body motions and commutes
with index raising and index lowering, respectively. The simplest corotational
variant is the {\it Zaremba-Jaumann time derivative} \cite{Jaum11GSPC,Zare03FPTR},
denoted it by a circle $(\cdot)\!\ZJ{^{\,}}$, which is well justified when
applied to the Cauchy stress tensor by Biot \cite[p.494]{Biot65MID}, cf.\ also
\cite{Bruh09EEBI,Fial11GSSM,Fial20OTDR,MorGio22OREM}. This choice is 
most commonly used, especially in geophysics
although, for some other applications where cycling regimes are expected, it may
exhibit undesired ``ratchetting'' effects \cite{JiaFis17ADRD,MeXiBrMe03ESRC}.
For {\it isotropic materials}, this derivative also affects  the symmetric
small-strain tensor $\Ee$, cf.\ \cite{Roub23SPTC}. For a survey of objective
corrotational strain rates see \cite{XiBrMe98SRMS,MeShBr00SCOR}.
Specifically, having the Eulerian velocity $\vv$, we define 
\begin{align}\label{ZJ}
\ZJ\Ee=
\pdt\Ee+(\vv\Cdot\nabla)\Ee-\WW\Ee+\Ee\WW\ \ \ \text{ with }\ \WW={\rm skw}(\nabla\vv)=
\frac{\nabla\vv{-}(\nabla\vv)^\top}2\,.
\end{align}
The Eulerian velocity $\vv$ is used also in the convective time derivative
\begin{align}
(\bm\cdot)\!\DT{^{}}=\pdt{}(\bm\cdot)+(\vv{\cdot}\nabla)(\bm\cdot)
\end{align}
to be used for scalars and, component-wise, for vectors or tensors.
Then \eq{ZJ} can be written shortly as
$\ZJ\Ee=\DT\Ee-{\rm skw}(\nabla\vv)\Ee+\Ee{\rm skw}(\nabla\vv)$.

\subsection{The system and its energetics}\label{sec-linearized-system}

We demonstrate the time-discretization method on the visco-elastodynamics
in {\it Kelvin-Voigt rheology} in the volumetric part and a {\it Jeffreys}
(also called {\it anti-Zener}) {\it rheology} in the deviatoric (isochoric)
part, which is a fairly general model that allows for isochoric creep. For the
displacement $\uu$, we implement the {\it Green-Naghdi additive decomposition}
\cite{GreNag65GTEP} of the total small strain $\strain(\uu)$ into the elastic
and the inelastic strain $\strain(\uu)=\Ee{+}\Ep$, expressed in objective rates as
$\strain(\vv)=\ZJ\Ee{+}\ZJ\Ep$ with denoting the inelastic strain (creep) rate
$\ZJ\Ep=:\Lp$.

The ingredients are the stored energy $\varphi:\Rsym\to\R$ acting on $\Ee$ and the
potential of dissipative forces $\zeta:\Rsym\times\Rsym\to\R$ acting on the rates
$\strain(\vv)$ and $\Lp$. We consider the isotropic material in which these
functionals depend only on the spherical (also called volumetric or hydrostatic)
and the deviatoric (also called shear or isochoric) parts separately. For the
(generally nonlinear) St.Venant–Kirchhoff-type models, we consider
$\varphi(\Ee)=\wt\varphi_1({\rm sph}\,\Ee)+\wt\varphi_2(|{\rm dev}\,\Ee|)$.
For the linear material, this means
\begin{subequations}\label{St.Venenat-Kirchhoff-linear}\begin{align}
&\varphi(\Ee)
=\frac32\Ke|{\rm sph}\,\Ee|^2\!+\Ge|{\rm dev}\,\Ee|^2\ \ \text{ with }\ \ 
{\rm sph}\,\Ee=\frac13({\rm tr}\Ee)\bbI\ \text{ and }\ 
{\rm dev}\,\Ee=\Ee-\frac13({\rm tr}\Ee)\bbI\,
\intertext{with the bulk elastic modulus $\Ke$ and the shear elastic modulus
$\Ge$ in Pa=J/m$^3$, and}
&\zeta(\strain(\vv),\Lp)=
\frac32\Kv|{\rm sph}\,\strain(\vv)|^2\!+\Gv|{\rm dev}\,\strain(\vv)|^2
+\frac12\Gm|\Lp|^2
\label{dissipation-pot}\end{align}\end{subequations}
with the bulk modulus $\Kv$ and the shear modulus $\Gv$ responsible for the
Stokes-type viscosity and $\Gm$ the modulus responsible for the Maxwellian
viscosity in Pa\,s=J\,s/m$^3$. The model \eq{St.Venenat-Kirchhoff-linear}
is schematically illustrated in Figure~\ref{fig-rheology-mixed}.

\begin{center}
\begin{my-picture}{25}{10.5}{fig-rheology-mixed}
\psfrag{VOLUMETRIC PART}{\begin{minipage}[t]{14em}\footnotesize
{\sf VOLUMETRIC PART}\\[-.2em]\footnotesize\hspace*{-.2em}({\sf Kelvin-Voigt rheology})\end{minipage}}
\psfrag{DEVIATORIC PART}{\begin{minipage}[t]{14em}\footnotesize
{\sf DEVIATORIC  PART}\\[-.2em]\footnotesize\hspace*{.5em}({\sf Jeffreys rheology})\end{minipage}}
\psfrag{e}{\small$\strain(\uu)$}
\psfrag{e1}{\small${\rm dev}\,\strain(\uu)$}
\psfrag{e2}{\begin{minipage}[t]{8em}\small${\rm sph}\,\strain(\uu)$\\[-.2em]\small$\hspace*{.2em}=({\rm div}\,\uu)\bbI/3$\end{minipage}}
\psfrag{E}{\small${\rm dev}\,\Ee$}
\psfrag{P}{\small$\Ep$}
\psfrag{c}{\small$\Ke$}
\psfrag{c1}{\small$\Ge$}
\psfrag{d2}{\small$\Kv$}
\psfrag{d11}{\small$\Gv$}
\psfrag{d21}{\small$\Gm$}
\hspace{-7em}\includegraphics[width=31em]{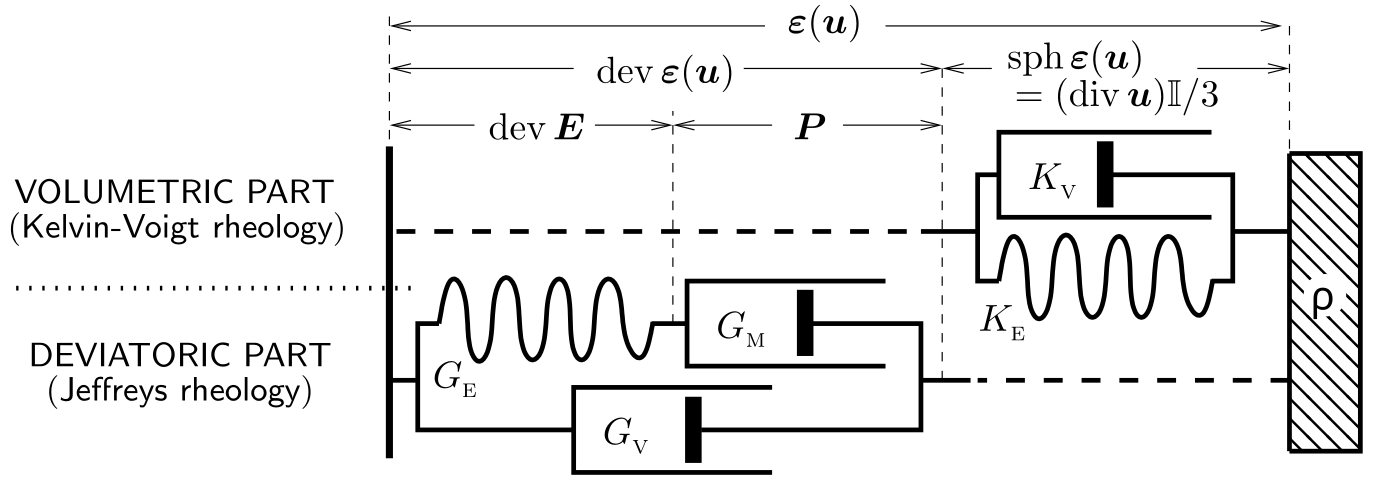}
\end{my-picture}
\nopagebreak
\\
{\small\sl\hspace*{-3em}Fig.\,\ref{fig-rheology-mixed}:~\begin{minipage}[t]{25em}
A schematic 1-dimensional illustration of combining the Kelvin-Voigt solid
rheologies in the volumetric part with the Jeffreys fluidic rheology in the
deviatoric (shear) part.
\end{minipage}
}
\end{center}

Rather for the notational simplicity, we confine ourselves to the
linear dissipation with the quadratic potential \eq{dissipation-pot},
and invent the 4th-order viscous-moduli tensor $\bbD$ defined as
$\bbD_{ijkl}=\Kv\delta_{ij}\delta_{kl}+\Gv(\delta_{ik}\delta_{jl}
{+}\delta_{il}\delta_{jk}{-}\frac23\delta_{ij}\delta_{kl})$ with ``$\delta$'' being
here the Kronecker symbol. In general, the creep rate is governed by the flow rule
of the type $\zeta_{\Lp}'(\strain(\vv),\Lp)=-{\rm dev}\,\TT$ where $\TT$ denotes
the Cauchy stress, which is given by $\TT=\varphi'(\Ee)+\varphi(\Ee)\bbI$. 
Considering \eq{dissipation-pot}, this flow rule reads as
\begin{align}\label{Euler-small-flow-rule}
\Gm\Lp=-{\rm dev}\,\TT\,,\ \text{ which couples with }\ \ZJ\Ee=\strain(\vv)-\Lp\,.
\end{align}
This merges into the single equation $\ZJ\Ee=\strain(\vv)-R\,{\rm dev}\,\TT$ with
$R:=1/\Gm$.

In addition, we include a nonlinear higher-order dissipation acting,
for simplicity, equally on the volumetric and on the shear parts. 

We express the system in terms of rates, i.e.\ we consider the mass density
$\varrho$, velocity $\vv$ (which is the Eulerian time derivative of the
displacement $\yy$), and the elastic strain $\Ee$, while do not explicitly
involve the displacement $\yy$ and, when considering the creep, neither the
inelastic strain. The system then consists of the continuity equation
for the mass density, the momentum equation, and the kinematic/rheological
equation for the triple $(\varrho,\vv,\Ee)$:
\begin{subequations}\label{Euler-small-viscoelastodyn+}
\begin{align}\label{Euler-small-viscoelastodyn+0}
&\!\!\DT\varrho=-\varrho{\rm div}\,\vv\,,
\\&\nonumber
\!\!\varrho\DT\vv={\rm div}\big(\TT{+}\DD
\big)+\varrho\GRAVITY\ \text{ with }\ \TT=\varphi'(\Ee)+\varphi(\Ee)\bbI\ \ 
\\&
\hspace*{4em}\text{ and }\ \ \DD=\bbD\strain(\vv)
-{\rm div}\mathfrak{H}
\ \text{ with }\ \mathfrak{H}=\HYPER|\nabla^2\vv|^{p-2}\nabla^2\vv\,,\!\!\!
\label{Euler-small-viscoelastodyn+1}
\\[-.3em]&
\!\!\ZJ\Ee=\strain(\vv)-R\,{\rm dev}\,\TT\,.
\label{Euler-small-viscoelastodyn+2}
\end{align}\end{subequations}
When $R=0$, the rheology degenerates to the Kelvin-Voigt model also in the deviatoric
part. The motivation for this ``hyper-viscosity'', which leads to the so-called
hyperstress $\mathfrak{H}$, is to facilitate mathematical analysis towards
conventional weak solutions and, only as a side effect, it would allow for certain
tuning of the normal dispersion of wave speed propagation, cf.\ \cite{Roub24SGTL},
and, if the $\nabla^2\vv$-dependent hyperviscosity coefficient $\HYPER|\nabla^2\vv|^{p-2}$
would be extended as $(\HYPER_1{+}\HYPER_2|\nabla^2\vv|^{p-2})$, also uniqueness
of the weak solution to the system \eq{Euler-small-viscoelastodyn+}, which is not
investigated in this paper, however.

The momentum equation \eq{Euler-small-viscoelastodyn+1} is to be supplemented
by the boundary conditions on the boundary $\varGamma$ of the domain $\varOmega$.
We consider the fixed shape of $\varOmega$ by zero normal velocity while otherwise
the boundary is free, i.e.
\begin{align}
&\vv\Cdot\nn=0,\ \ \ [(\TT{+}\DD)\nn+\divS(\mathfrak{H}\nn)]_\text{\sc t}^{}
\bm0\,,\ \ \text{ and }\ \ \mathfrak{H}\Colon(\nn{\otimes}\nn)={\bm0}\ \
\text{ on }\ \varGamma\,.
\label{Euler-small-BC-hyper}\end{align}
The energy balance behind this system can then be revealed by
testing \eq{Euler-small-viscoelastodyn+1} by $\vv$ and using 
the Green formula twice over $\varOmega$ and once over $\varGamma$, together with
\eq{Euler-small-viscoelastodyn+0} multiplied by $|\vv|^2/2$ and 
\eq{Euler-small-viscoelastodyn+2} multiplied by $\TT$.

For the inertial form in \eq{Euler-small-viscoelastodyn+1} tested by $\vv$, we
employ \eq{Euler-small-viscoelastodyn+0} in the form $\pdt{}\varrho=
-{\rm div}(\varrho\vv)$ tested by $\frac12|\vv|^2$:
\begin{align}
  \pdt{}\bigg(\frac\varrho2|\vv|^2\bigg)=\varrho\vv\Cdot\pdt\vv
  +\pdt\varrho\frac{|\vv|^2}2\!\stackrel{{\scriptsize\rm\eq{Euler-small-viscoelastodyn+0}}}{=}\!\!\varrho\vv\Cdot\pdt\vv
-{\rm div}(\varrho\vv)\frac{|\vv|^2}2\,.
\label{rate-of-kinetic}\end{align}
Then, using \eq{rate-of-kinetic} and Green's formula with $\vv\Cdot\nn=0$ on
$\varGamma$, it holds
\begin{align}\nonumber
\!\!\!\!\!\int_\varOmega\varrho(\vv\Cdot\nabla)\vv\Cdot\vv\,\d\xx
&=\int_\varGamma\varrho|\vv|^2\!\!\!\!\lineunder{\!\!\!\!\vv\Cdot\nn\!\!\!\!}{$=0$}\!\!\!\!\d S
-\!\!\int_\varOmega\!\vv\Cdot{\rm div}(\varrho\vv{\otimes}\vv)\,\d\xx
=-\!\!\int_\varOmega\varrho(\vv\Cdot\nabla)\vv\Cdot\vv+{\rm div}(\varrho\vv)|\vv|^2\,\d\xx
\\&
=-\int_\varOmega\!
{\rm div}(\varrho\vv)\frac{|\vv|^2}2\,\d\xx
\!\stackrel{{\scriptsize\rm\eq{rate-of-kinetic}}}{=}\!\frac{\d}{\d t}
\int_\varOmega\frac\varrho2|\vv|^2\,\d\xx-\!\int_\varOmega\!\varrho\pdt\vv\Cdot\vv\,\d\xx\,,\!\!\!\!
\label{calculus-convective}\end{align}
i.e.\ $\int_\varOmega\varrho\DT\vv\Cdot\vv\,\d\xx=\frac{\d}{\d t}\int_\varOmega\frac12\varrho|\vv|^2\,\d\xx$.

For the Cauchy stress $\TT$, we have use the calculus
\begin{align}\nonumber
\int_\varOmega&\!({\rm div}\TT)\Cdot\vv\,\d\xx=\!\int_\varGamma\!\!\vv\Cdot\TT\nn\,\d S
-\!\!\int_\varOmega\!\!\TT\Colon\strain(\vv)\,\d\xx
\!\stackrel{\eq{Euler-small-viscoelastodyn+1}}{=}\!\!
\int_\varGamma\!\!\vv\Cdot\TT\nn\,\d S
-\!\!\int_\varOmega\!\!\big(\varphi'(\Ee){+}\varphi(\Ee)\bbI\big)\Colon\strain(\vv)\,\d\xx
\\[-.3em]&\nonumber\!\stackrel{\eq{Euler-small-viscoelastodyn+2}}{=}\!\!\!\int_\varGamma\!\!\vv\Cdot\TT\nn\,\d S
-\int_\varOmega\varphi'(\Ee)\Colon\Big(\ZJ\Ee{-}R\,{\rm dev}\,\TT\Big)+\varphi(\Ee){\rm div}\,\vv\,\d\xx
\\&\nonumber\!\stackrel{\eq{ZJ}}{=}\!\int_\varGamma\!\!\vv\Cdot\TT\nn\,\d S
-\int_\varOmega\!\varphi'(\Ee)\Colon\Big(\pdt{\Ee\!}+(\vv\Cdot\Nabla)\Ee-\WW\Ee+\Ee\WW\Big)
-R\,\big|{\rm dev}\,\varphi'(\Ee)\big|^2+\varphi(\Ee){\rm div}\,\vv\,\d\xx
\\[-1.1em]&\!\!\!\!\UUU{=}{\eq{Euler-small-calc}}{\eq{E:WS=E:SW}}\!\!\!\!
\int_\varGamma\!\!\vv\Cdot\TT\nn\,\d S
-\frac{\d}{\d t}\int_\varOmega\varphi(\Ee)\,\d\xx
-\int_\varOmega R\,\big|{\rm dev}\,\varphi'(\Ee)\big|^2\,\d\xx\,.
\label{Euler-small-divT.v++}\end{align}
Here, we used the calculus
\begin{align}\label{Euler-small-calc}
\int_\varOmega\!\varphi'(\Ee)\Colon(\vv\Cdot\Nabla)\Ee+\varphi(\Ee){\rm div}\,\vv\,\d\xx=\!\!\int_\varOmega\!\nabla\varphi(\Ee)\Cdot\vv+\varphi(\Ee){\rm div}\,\vv\,\d\xx=\!\!\int_\varGamma\!\varphi(\Ee)\!\!\lineunder{\!\!(\vv\Cdot\nn)\!\!}{$=0$}\!\!\d S
\end{align}
and also, taking into account the form of the spin
$\WW=\frac12\Nabla\vv-\frac12(\Nabla\vv)^\top$, we have used 
$\varphi'(\Ee)\Colon(\WW\Ee-\Ee\WW)=\bm0$. More in detail, by abbreviating
$\Se=\varphi'(\Ee)$ and using the algebra
$A\Colon (BC)=(B^\top\!A)\Colon C=(AC^\top)\Colon B$, it holds 
\begin{align}\nonumber
\Se\Colon(\WW\Ee-\Ee\WW)&=
\frac12\Se\Colon\big((\Nabla\vv)\Ee-(\Nabla\vv)^\top\Ee
-\Ee(\Nabla\vv)+\Ee(\Nabla\vv)^\top\big)
\\&\nonumber=\frac12\big(\Se\Ee^\top\!-\Ee^\top\Se\big)\Colon\Nabla\vv
-\frac12\big((\Nabla\vv)\Se-\Se(\Nabla\vv)\big)\Colon\Ee
\\&=\frac12\big(\Se\Ee^\top\!\!-\Ee\Se^\top\!\!-\Ee^\top\Se+\Se^\top\Ee\big)\Colon\Nabla\vv
=
\big(\!\!\!\lineunder{\Se\Ee\!-\Ee\Se\!\!}{$={\bm0}$}\!\!\!\big)\Colon\Nabla\vv
={\bm 0}\,;
\label{E:WS=E:SW}
\end{align}
here we have to assume that the initial condition for $\Ee$ is symmetric and
exploit that the Zaremba-Jaumann corotational derivative in
\eq{Euler-small-viscoelastodyn+2} keeps symmetry of $\Ee$ during
the whole evolution 
and that the material is considered isotropic so that $\Se$ commutes with $\Ee$.
Indeed, always we have
$\Se\Ee=({\rm sph}\,\Se+{\rm dev}\,\Se)({\rm sph}\,\Ee+{\rm dev}\,\Ee)
=\Ee\Se$ whenever $({\rm dev}\,\Se)({\rm dev}\,\Ee)=({\rm dev}\,\Ee)({\rm dev}\,\Se)$,
which is indeed true in particular for the mentioned isotropic 
St.Venant–Kirchhoff-type models.

The further contribution from the dissipative part of the Cauchy stress
uses Green's formula over $\varOmega$ twice and the surface Green formula
over $\varGamma$. We abbreviate the ``simple stress'' by
$\Se=\TT+\bbD(\theta)\strain(\vv)$. 
Then
\begin{align}\nonumber
&\int_\varOmega{\rm div}
\big(\Se-{\rm div}\mathfrak{H}\big){\cdot}\vv\,\d\xx
=\int_\varGamma\vv{\cdot}
\big(\Se{-}{\rm div}\mathfrak{H}\big)\nn\,\d S
-\!\int_\varOmega\!\big(\Se-{\rm div}\mathfrak{H}\big){:}\Nabla\vv\,\d\xx
\\&\nonumber=\int_\varGamma\vv{\cdot}
\big(\Se{-}{\rm div}\mathfrak{H}\big)\nn
-\nn{\cdot}\mathfrak{H}{:}\Nabla\vv\,\d S
-\!\int_\varOmega\!
\Se{:}\strain(\vv)+\mathfrak{H}\Vdots\Nabla^2\vv\,\d\xx
\\&\nonumber
=\int_\varGamma\mathfrak{H}{:}(\nn{\otimes}\nn)
+\nn{\cdot}\mathfrak{H}{:}\NablaS\vv+\vv{\cdot}
\big(\Se{-}{\rm div}\mathfrak{H}\big)\nn\,\d S
-\!\int_\varOmega\Se{:}\strain(\vv)
+\mathfrak{H}\Vdots\Nabla^2\vv\,\d\xx
\\&=\int_\varGamma\mathfrak{H}{:}(\nn{\otimes}\nn)-
\big(\divS(\nn{\cdot}\mathfrak{H})+
\big(\Se{-}{\rm div}\mathfrak{H}\big)\nn\big)
{\cdot}\vv\,\d S
-\!\int_\varOmega\Se{:}\strain(\vv)
+\HYPER|\nabla^2\vv|^p\,\d\xx\,,
\label{Euler-test-momentum++}\end{align}
where the hyperstress $\mathfrak{H}=\HYPER|\nabla^2\vv|^{p-2}\nabla^2\vv$ is from
\eq{Euler-small-viscoelastodyn+1} and where
we used the decomposition of $\Nabla\vv$ into its normal and tangential
parts, i.e.\ written componentwise
$\nabla\vv_i=(\nn{\cdot}\nabla\vv_i)\nn+\nablaS\vv_i$. This is to be merged with
\eq{Euler-small-divT.v++}.

Relying on the boundary conditions \eq{Euler-small-BC-hyper}
we finally obtain the {\it energy-dissipation balance}
\begin{align}
\!\!\!\!\!\frac{\d}{\d t}\int_\varOmega\!\!\!\!
\linesunder{\frac\varrho2|\vv|^2}{kinetic}{energy}
\!\!\!\!\!+\!\!\!\!\!\linesunder{\varphi(\Ee)}{stored}{energy}\!\!\!\!
\,\d\xx
+\!\int_\varOmega\hspace{-.7em}
\linesunder{\bbD\strain(\vv)\Colon\strain(\vv){+}\HYPER|\nabla^2\vv|^p}{disipation rate due to}{the Stokes-type viscosity}\hspace{-.7em}
+\hspace{-1em}\linesunder{R\,\big|{\rm dev}\,\varphi'(\Ee)\big|^2}{disipation rate due to}{Maxwellian viscosity}\hspace{-1em}\d\xx
=\int_\varOmega\hspace{-1.4em}\linesunder{\varrho\GRAVITY\Cdot\vv}{\ power of}{gravity field\ }\hspace{-1.4em} \,\d\xx
  \,.\ 
\nonumber\\[-2.5em]\label{Euler-small-energy-balance-stress}\end{align}

\medskip

\subsection{Weak solutions to an initial-value problem for \eq{Euler-small-viscoelastodyn+}}

We consider an initial-value problem for the evolutionar boundary-value problem 
\eq{Euler-small-viscoelastodyn+}--\eq{Euler-small-BC-hyper} by prescribing
the initial conditions
\begin{align}
\varrho|_{t=0}^{}=\varrho_0\,,\ \ \ \ \vv|_{t=0}^{}=\vv_0\,,
\ \ \text{ and }\ \ \Ee|_{t=0}^{}=\Ee_0\,.
\label{Euler-small-IC-hyper}\end{align}

The weak formulation of the momentum equation uses the by-part integration in
time and the Green formula over $\varOmega$ twice and once over $\varGamma$
exploiting also the boundary conditions \eq{Euler-small-BC-hyper}.

\begin{definition}[Weak formulation of 
\eq{Euler-small-viscoelastodyn+}--\eq{Euler-small-BC-hyper} with \eq{Euler-small-IC-hyper}]\label{def-ED-Ch4}
We call $(\varrho,\vv,\Ee)$ with
$\varrho\in L^\infty(I{\times}\varOmega)\cap W^{1,1}(I{\times}\varOmega)$,
$\vv\in L^p(I;W^{2,p}(\varOmega;\R^3))$,
$\Ee\in 
W^{1,1}(I{\times}\varOmega;\Rsym)$ a weak solution to the system
\eq{Euler-small-viscoelastodyn+} with the boundary conditions
\eq{Euler-small-BC-hyper} and the initial conditions \eq{Euler-small-IC-hyper} if 
\begin{align}
&\nonumber
\int_0^T\!\!\!\int_\varOmega\bigg(\Big(
\varphi'(\Ee)+\bbD\strain(\vv)-\varrho\vv{\otimes}\vv\Big){:}\strain(\widetilde\vv)
+\varphi(\Ee){\rm div}\,\widetilde\vv
-\varrho\vv{\cdot}\pdt{\widetilde\vv}
\\[-.6em]&\hspace{7em}
+\HYPER|\nabla^2\vv|^{p-2}\nabla^2\vv\Vdots\nabla^2\widetilde\vv\bigg)\,\d\xx\d t
=\!\int_0^T\!\!\!\int_\varOmega\varrho\GRAVITY{\cdot}\widetilde\vv\,\d\xx\d t
+\!\int_\varOmega\!\varrho_0\vv_0{\cdot}\widetilde\vv(0)\,\d\xx
\label{def-ED-Ch4-momentum}\end{align}
holds for any $\widetilde\vv$ smooth with $\widetilde\vv{\cdot}\nn={\bm0}$ on
$I{\times}\varGamma$ and $\widetilde\vv(T)=0$, and if
\eq{Euler-small-viscoelastodyn+0} and \eq{Euler-small-viscoelastodyn+2} hold
a.e.\ on $I{\times}\varOmega$ together with the respective initial conditions for
$\varrho$ and $\Ee$ in \eq{Euler-small-IC-hyper}.
\end{definition}

Let us summarize the assumptions used in what follows:
\begin{subequations}\label{Euler-small-ass}
\begin{align}\label{Euler-small-ass-phi}
&\!\!\varphi\in C^1(\Rsym):\ \ \inf_{E\ne0}\frac{\varphi(E)}{|E|^2}>0\,,\ \ 
\sup_{E}\frac{\varphi'(E)}
{1{+}|E|\!}<+\infty\,,\ \ 
\exists \alpha>0{:}\ \varphi(\cdot)-\alpha|\cdot|^2\,\text{ is convex},\!\!
\\[-.3em]&\!\!\label{Euler-small-ass-D}
\bbD:\Rsym\to\Rsym\ \text{ linear symmetric}\,,\ \
\min_{|E|=1}\bbD E\Colon E=\nu_1>0\,,\ \ \ \HYPER>0\,,
\\[-.4em]&\!\!\label{Euler-small-ass-IC}
\varrho_0\in W^{1,r}(\varOmega)\,,\ \ \ {\rm min}_{\barOmega}^{}\varrho_0>0\,,\ \ \
\vv_0\in L^2(\varOmega;\R^3)\,,\ \ \ \Ee_0\in H^1(\varOmega;\Rsym)\,,
\\&\!\!\label{Euler-small-ass-g}
\GRAVITY\in L^1(I;L^\infty(\varOmega;\R^3))
\,.\end{align}\end{subequations}
The last condition in \eq{Euler-small-ass-phi} means the {\it strong convexity}
of $\varphi$. For the generalization of \eq{Euler-small-ass-phi} to a more
general $\varphi$ and \eq{Euler-small-ass-g} towards
$\GRAVITY\in L^1(I;L^\infty(\varOmega;\R^3))$
see Remarks~\ref{rem-general-phi} and \ref{rem-f-ne-0} below.

The following statement will be proved in Sect.\,\ref{sec-linearized-proof} below
by a rather constructive time discretization:

\begin{proposition}[Existence of weak solutions to \eq{Euler-small-viscoelastodyn+}]
\label{prop-ED-Ch4-existence}
Let \eq{Euler-small-ass} hold with $r>3$ and let $p>3$. Then the initial-boundary-value
problem for the system \eq{Euler-small-viscoelastodyn+} has at least one weak
solution $(\varrho,\vv,\Ee)$ in the sense of Definition~\ref{def-ED-Ch4} such that also
$\varrho\in L^\infty(I;W^{1,r}(\varOmega))\,\cap\,C(I{\times}\barOmega)$ with
$\min_{I{\times}\barOmega}\varrho>0$, $\vv\in L^\infty (I;L^2(\varOmega;\R^3))$,
$\Ee\in L^\infty(I;H^1(\varOmega;\Rsym))$ and such that the energy-dissipation balance
\eq{Euler-small-energy-balance-stress} integrated over the time interval
$[0,t]$ holds for any $t\in I$.
\end{proposition}

\subsection{Reformulation of the system \eq{Euler-small-viscoelastodyn+} and time discretization}

We express the system in terms of the linear momentum $\pp=\varrho\vv$ and,
later, exploit the fact that the kinetic energy $(\pp,\varrho)\mapsto
\frac12|\pp|^2/\varrho$ is convex, in contrast to the equivalent form
$(\vv,\varrho)\mapsto\frac12\varrho|\vv|^2$ which is nonconvex.
The system \eq{Euler-small-viscoelastodyn+} then reads as
\begin{subequations}\label{Euler-small-viscoelastodyn+trans}
\begin{align}\label{Euler-small-viscoelastodyn+0trans}
&
\!\!\pdt\varrho=-{\rm div}\,\pp\,,
\\&\nonumber
\!\!\pdt\pp={\rm div}\big(\TT{+}\DD
-\pp{\otimes}\vv\big)+\varrho\GRAVITY
\ \text{ with }\ \TT=\varphi'(\Ee)+\varphi(\Ee)\bbI
\\[-.4em]&
\hspace*{4em}\text{ and }\ \DD=\bbD\strain(\vv)
-{\rm div}\mathfrak{H}\,,
\ \text{ where }\ 
\mathfrak{H}=\HYPER|\nabla^2\vv|^{p-2}\nabla^2\vv\ \text{ and }\ \ \pp=\varrho\vv\,,
\label{Euler-small-viscoelastodyn+1trans}
\\[-.3em]&
\!\!\pdt\Ee=\strain(\vv)-R\,{\rm dev}\,\TT-\bm B_\text{\sc zj}^{}(\vv,\bm E)\,,
\label{Euler-small-viscoelastodyn+2trans}
\end{align}\end{subequations}
where we have used a shorthand notation for the bi-linear operator involved in
the Zaremba-Jaumann derivative:
 \begin{align}
 \bm B_\text{\sc zj}^{}(\vv,\bm E)=(\vv{\cdot}\nabla)\bm E-
      {\rm skew}(\nabla\vv)\Ee+\Ee\,{\rm skew}(\nabla\vv)\,.
\label{ZJ-def}\end{align}

The energetics can be revealed by testing \eq{Euler-small-viscoelastodyn+0trans}
by $|\pp|^2/(2\varrho^2)$ and \eq{Euler-small-viscoelastodyn+1trans}
by $\pp/\varrho$, which gives
\begin{align}
  \pdt{}\bigg(\frac{|\pp|^2}{2\varrho}\bigg)
  &=\frac{\pp}{\varrho}\Cdot\pdt\pp
  -\frac{|\pp|^2}{2\varrho^2}\pdt\varrho
=\frac{\pp}{\varrho}\Cdot
\Big({\rm div}\big(\TT{+}\DD-\vv{\otimes}\pp\big)+\varrho\GRAVITY\Big)
+\frac{|\pp|^2}{2\varrho^2}{\rm div}\,\pp
\,.\label{rate-of-kinetic-modif}\end{align}
Then, using the calculus
\begin{align}\label{calculus-for-kinetic}
\!\!\int_\varOmega\frac{\pp}{\varrho}\Cdot{\rm div}\big(\vv{\otimes}\pp\big)
-\frac{|\pp|^2\!}{2\varrho^2\!}\,{\rm div}\,\pp\,\d\xx
=\int_\varOmega\!\vv\Cdot{\rm div}\big(\varrho\vv{\otimes}\vv\big)
-\frac{|\vv|^2\!}{2}\,{\rm div}(\varrho\vv)\,\d\xx
\!\!\stackrel{\scriptsize\eq{calculus-convective}}{=}\!0\,,
\end{align}
we again arrive to the energy-dissipation balance
\eq{Euler-small-energy-balance-stress} with the kinetic energy in the form
$\frac12|\pp|^2/\varrho$.

Now we devise the time discrete scheme for the system
\eq{Euler-small-viscoelastodyn+trans} regularized by adding gradient terms
into the transport equations for $\varrho$ and $\Ee$ and an artifical nonlinear
0-order dissipative term into the momentum equation together
with a compensating gradient term. We use the so-called {\it Rothe method},
i.e.\ the fully implicit time discretisation with an equidistant partition of
the time interval $I$ with the time step $\tau>0$ such that $T/\tau\in\N$. 
We denote by $\varrho_\etau^k$, $\vv_\etau^k$, $\Ee_\etau^k$, $\TTtauk$ ...
the approximate values of $\varrho$, $\vv$, $\Ee$, $\TT$ ...
at time instants $t=k\tau$ with $k=1,2,...,T/\tau$. For $\DELTA>0$ and $\EPS>0$,
we will then use the following recursive regularized time-discrete scheme 
\begin{subequations}\label{Euler-small-viscoelastodyn+disc}
\begin{align}\label{Euler-small-viscoelastodyn+0disc}
&\!\!\frac{\varrho_\etau^k{-}\varrho_\etau^{k-1}\!\!}\tau\,
={\rm div}\Big(\DELTA|\nabla\varrho_\etau^k|^{r-2}
\nabla\varrho_\etau^k-\,\CUT(\varrho_\etau^k)\pp_\etau^k\Big)\,,
\\[-.2em]&\nonumber
\!\!\frac{\pp_\etau^k{-}\pp_\etau^{k-1}\!\!}\tau\,=
{\rm div}\Big(\TTtauk{+}\DD_\etau^k{-}\CUT(\varrho_\etau^k)\pp_\etau^k{\otimes}\vvk\Big)
+\varrho_\etau^k\GRAVITY_{\DELTA\tau}^k
\\[-.3em]&\nonumber\hspace*{7em}
-\EPS|\vvk|^{p-2}\vvk
-\DELTA|\nabla\varrho_\etau^k|^{r-2}(\nabla\vvk)\nabla\varrho_\etau^k
\\[-.1em]&\nonumber
\hspace*{5em}\text{ with }\ \TTtauk=\varphi'(\Ee_\etau^k)+\varphi(\Ee_\etau^k)\bbI
\ \ \text{ and }\ \DD_\etau^k=\bbD\strain(\vvk)-{\rm div}\mathfrak{H}_\etau^k\,,
\\[-.0em]&
\hspace*{5em}
\text{ where }\ \mathfrak{H}_\etau^k=\HYPER\big|\nabla^2\vvk\big|^{p-2}
\nabla^2\vvk
\ \ \ \text{ and }\ \ \ \pp_\etau^k=\varrho_\etau^k\vvk\,,
\label{Euler-small-viscoelastodyn+1disc}
\\[-.5em]&
\!\!\frac{\Ee_\etau^k{-}\Ee_\etau^{k-1}\!\!}\tau\,
=\strain(\vvk)-R\,{\rm dev}\,\TTtauk-\bm B_\text{\sc zj}^{}(\vvk,\Ee_\etau^k)
+\DELTA\Delta\Ee_\etau^k
\label{Euler-small-viscoelastodyn+2disc}
\end{align}\end{subequations}
with $\GRAVITY_{\DELTA\tau}^k=\frac1\tau\int_{(k-1)\tau}^{k\tau}\GRAVITY_\DELTA(t)\,\d t$
with some $\GRAVITY_\DELTA\in L^{p'}(I;L^\infty(\varOmega;\R^3))$, with
$\bm B_\text{\sc zj}^{}$ from \eq{ZJ-def}, and with a
``cut-off'' function 
$\CUT\in L^\infty(\R)\cap C^1(\R\setminus\{0\})$ such that  
\begin{align}\label{cut-off}
\CUT(\varrho):=\begin{cases}\ \ 1&\text{for }\ \ 
0\le\varrho\le\rhoMAX,\\
\ \ 0&\text{for }\ \varrho<0\ \text{\ or }\ \varrho\ge\rhoMAX{+}1,
\\
\in[0,1]&\text{for }\ \rhoMAX<\varrho<\rhoMAX{+}1
\end{cases}
\end{align}
with some $0<\rhoMAX$ to be chosen later. Note that, inspite of the discontinuity
of $\CUT$ at 0, the function $\varrho\mapsto\CUT(\varrho)\varrho$ is a
Lipschitz continuous, non-negative function. Having in mind the qualification
\eq{Euler-small-ass-g}, we assume
\begin{align}
\GRAVITY_\DELTA\to\GRAVITY\ \ \text{ strongly in }\ L^1(I;L^\infty(\varOmega;\R^3))
\ \text{ for $\DELTA\to0$}\,.
\end{align}

The corresponding boundary
conditions now expanded due to the regularizing terms in
\eq{Euler-small-viscoelastodyn+0disc} and in
\eq{Euler-small-viscoelastodyn+2disc}:
\begin{subequations}\label{BC-disc}\begin{align}\label{BC-disc-a}
&\!\pp_\etau^k{\cdot}\nn=0\,,\ \ \,
\Big[\Big(\TTtauk{+}\bbD\strain(\vvk){-}{\rm div}\mathfrak{H}_\etau^k\Big)\nn
-\divS\big(\mathfrak{H}_\etau^k\Cdot\nn\big)\Big]_\text{\sc t}^{}\!\!=\bm0\,,\ \ \,
\nabla^2\vvk{:}(\nn{\otimes}\nn)=0\,,
\\&\!\nn\Cdot\nabla\varrho_\etau^k=0\,,\ \ \text{ and }\ \ \
(\nn\Cdot\nabla)\Ee_\etau^k=\bm0\,.
\end{align}\end{subequations}

Such a time discretization was used for compressible Navier-Stokes equations
in \cite{Kar13CFEM,Zato12ASCN} or \cite[Ch.7]{FeKaPo16MTCV}. The explicit usage
of the convexity of the kinetic energy expressed in terms of the momentum $\pp$
is in \cite{FLMS21NACF}. The regularization in
\eq{Euler-small-viscoelastodyn+0disc} together with the ``compensating'' term
in \eq{Euler-small-viscoelastodyn+1disc} and the cut-off is similar as in
\cite{Zato12ASCN} although the cut-off here is different.

The philosophy of our cut-off by $\CUT$ is that, if
$0\le\varrho_\etau^{k-1}\le\rhoMAX{+}1$ on $\varOmega$, then the discrete
regularized continuity equation \eq{Euler-small-viscoelastodyn+0disc} cannot 
have a solution $\varrho_\etau^k$ with values below
$0$ (or even equal to 0), nor above $\rhoMAX{+}1$. This can be seen
by a contradiction argument: assuming $\varrho_\etau^{k-1}>0$ on $\barOmega$ and
the minimum of a (momentarily smooth) solution $\varrho_\etau^k$ is attained at
some $\xx\in\varOmega$ and $\varrho_\etau^k(\xx)\le0$,
so that $\varrho_\etau^k(\xx)-\varrho_\etau^{k-1}(\xx)<0$ and
$\nabla\varrho_\etau^k(\xx)=\bm0$ and also
$\nabla(\CUT(\varrho_\etau^k(\xx))\varrho_\etau^k(\xx))=\bm0$
and also ${\rm div}(\DELTA|\nabla\varrho_\etau^k|^{r-2}\nabla\varrho_\etau^k)\ge0$,
from \eq{Euler-small-viscoelastodyn+0disc} written at $\xx$ in the
form $\varrho_\etau^k(\xx)-\varrho_\etau^{k-1}(\xx)
=-\tau\vv_\etau^k(\xx)\Cdot
\nabla(\CUT(\varrho_\etau^k(\xx))\varrho_\etau^k(\xx))
-\tau({\rm div}\vv_\etau^k(\xx))\CUT(\varrho_\etau^k(\xx))\varrho_\etau^k(\xx)
+\tau{\rm div}(\DELTA|\nabla\varrho_\etau^k|^{r-2}\nabla\varrho_\etau^k)\ge0$,
we obtain the contradiction showing that $\varrho_\etau^k(\xx)>0$.
Moreover, if the maximum of $\varrho_\etau^k$ would be attained
at some $\xx\in\varOmega$ and $\varrho_\etau^k(\xx)>\rhoMAX{+}1$
while $\varrho_\etau^{k-1}(\xx)\le\rhoMAX{+}1$,
we would have $\varrho_\etau^k(\xx)-\varrho_\etau^{k-1}(\xx)>0$ and
$\nabla\varrho_\etau^k(\xx)=\bm0$ and also
${\rm div}(\DELTA|\nabla\varrho_\etau^k(\xx)|^{r-2}\nabla\varrho_\etau^k(\xx))\le0$,
so that, from \eq{Euler-small-viscoelastodyn+0disc} written at $\xx$,
we would now obtain $\varrho_\etau^k(\xx)-\varrho_\etau^{k-1}(\xx)\le0$,
a contradiction. Thus, altogether, $0<\varrho_\etau^k(\xx)\le\rhoMAX{+}1$
for all $\xx\in\varOmega$.

This system of boundary-value problems for the triple
$(\varrho_\etau^k,\pp_\etau^k,\Ee_\etau^k)$ is to be solved recursively for
$k=1,2,...,T/\tau$, starting for $k=1$ with the initial conditions: 
\begin{align}\label{IC-disc}
\varrho_\etau^0=\varrho_0\,,\ \ \ \ \ \ \ \ \pp_\etau^0=\varrho_0\vv_{0}\,,\  
    \ \ \text{ and }\ \ \ {\Ee}_\etau^0={\Ee}_0\,.
\end{align}
Thus, from \eq{Euler-small-viscoelastodyn+1disc}, we obtain also
$\vvk=\pp_\etau^k/\varrho_\etau^k$ providing $\varrho_\etau^k>0$, as indeed
proved later. 

\subsection{Stability and convergence, the proof of Proposition~\ref{prop-ED-Ch4-existence}}\label{sec-linearized-proof}

It is important that the regularization in \eq{Euler-small-viscoelastodyn+disc}
is devised so that:
\begin{itemize}
\vspace*{-.7em}\item
we avoid using the discrete Gronwall inequality on the
time-discrete scheme, which could not work on an equidistant time
discretization. More in detail,
$\nabla\vv_\etau^k\in L^{\infty}(\varOmega;\R^{3\times3})$, which would act as
coefficients in the discrete Gronwall inequality, is not uniformly bounded
in $k$, while a conservative-type regularization of the momentum equation,
which would guarantee such uniform bound, does not seem feasible for the limit
passage. In particular, we avoid testing \eq{Euler-small-viscoelastodyn+2disc}
by $\Ee_\etau^k$ and read the estimate for $\nabla\Ee_\etau^k$ from the
strong convexity of $\varphi$ instead. 
\vspace*{-.7em}\item
The other technicality is a
careful two-step regularization of the system to avoid testing of
the transport equations by the Laplacian on the time-discrete level
(due to the aforementioned troubles with the discrete Gronwall inequality) and
to have some estimate of the velocity $\vv$ even if the regularized mass-density
continuity equation itself does not guarantee that the mass density
$\varrho$ remains well far 0.
\end{itemize}\vspace*{-.7em}

For clarity, we will divide the following argumentation into eight steps.

\medskip{\it Step 1: The choice of $\rhoMAX$ 
in \eq{cut-off}}.
We use (formally at this point) the energy-dissipation balance
\eq{Euler-small-energy-balance-stress} together with $\varrho\ge0$
and the continuity equation \eq{Euler-small-viscoelastodyn+0trans} and the
impenetrability of the boundary, which guarantees 
\begin{align}\label{Euler-small-est-1}
\int_\varOmega\varrho(t)\,\d\xx=\int_\varOmega\varrho_0\,\d\xx=:M\,.
\end{align}
Treating the gravity loading as in \eq{Euler-small-est-Gronwall} below,
we obtain the a-priori bounds
$\|\strain(\overlinevvtau)\|_{L^2(I{\times}\varOmega;\R^{3\times3})}^{}\le C$ and 
$\|\nabla^2\overlinevvtau\|_{L^p(I{\times}\varOmega;\R^{3\times3\times3})}^{}\le C$.
This a-priori quality of the velocity field together with the
regularity of the initial condition for $\varrho_0$ in 
\eq{Euler-small-ass-IC} allows for using the transport-and-evolution
equation \eq{Euler-small-viscoelastodyn+0} to obtain the a-priori bound
$\|\varrho\|_{L^\infty(W^{1,r}(\varOmega))}\le C_r$, which suggest to
put $\rhoMAX>N_rC_r$ with $N_r$ denoting the norm of the embedding
$W^{1,r}(\varOmega)\subset L^\infty(\varOmega)$ holding for $r>3$. For details
see \cite[Lemma\,5.1]{Roub24TVSE}; actually, rather for analytical
arguments we need $\vv$ also in $L^1(I{\times}\varOmega;\R^3)$ although
not a-priori bounded, which is seen from \eq{Euler-small-energy-balance-stress}
provided $\min_{I{\times}\barOmega}\varrho>0$ as will actually be granted
in what follows.

\medskip{\it Step 2: Basic stability of the scheme \eq{Euler-small-viscoelastodyn+disc}
and first a-priori estimates}.
The discrete analog of the dissipation-energy balance
\eq{Euler-small-energy-balance-stress} can be obtained by testing
\eq{Euler-small-viscoelastodyn+1disc} by $\vv_\etau^k:=\pp_\etau^k/\varrho_\etau^k$
while using \eq{Euler-small-viscoelastodyn+0disc} tested by
$\frac12|\pp_\etau^k|^2/(
\varrho_\etau^k)^2$ and
\eq{Euler-small-viscoelastodyn+2disc} tested by $\varphi'(\Ee_\etau^k)$.

The essential point is the mentioned convexity of the function
$(\pp,\varrho)\mapsto\frac12|\pp|^2/\varrho:\R^3\times\R^+\mapsto\R$ with
$\R^+=(0,+\infty)$. To prove it, we investigate its 2nd-derivative
(Hessian) which is the $(3{+}1)\times(3{+}1)$-matrix
$$
\bigg(\frac{|\pp|^2}{2\varrho}\bigg)''=
\bigg(\begin{array}{c}\pp/\varrho\\\!\!-|\pp|^2/(2\varrho^2)\\\end{array}\bigg)'
=\bigg(\begin{array}{cc}\!(1/\varrho)\bbI\:, &\  -\pp/\varrho^2
\\\!-\pp/\varrho^2\,, &\ |\pp|^2/\varrho^3\end{array}\bigg)\,.
$$
The convexity to be proved is equivalent to the positive semi-definiteness of this
matrix. Since surely the matrix $(1/\varrho)\bbI$ is positive definite,
by Sylvester's criterion, it suffices to check that the determinant of this
$(3{+}1)\times(3{+}1)$-matrix is non-negative. By the
 ``block determinant'' formula, this total determinant is equal to
$\det((1/\varrho)\bbI)\,
\det(|\pp|^2/\varrho^3-\pp^\top\!/\varrho^2[(1/\varrho)\bbI]^{-1}\pp/\varrho^2)
=(1/\varrho^3)(|\pp|^2/\varrho^3-\pp\Cdot\pp/\varrho^3)=0$.

Due to this convexity, the discrete analog of the equality
\eq{rate-of-kinetic-modif} reads as an inequality 
\begin{align}\nonumber
  \!\frac1\tau&\bigg(\frac{|\pp_\etau^k|^2}{2\varrho_\etau^k}
  -\frac{|\pp_\etau^{k-1}|^2}{2\varrho_\etau^{k-1}}\bigg)
\le\frac{\pp_\etau^k}{\varrho_\etau^k}\Cdot\frac{\pp_\etau^k-\pp_\etau^{k-1}\!\!}\tau
-\frac{|\pp_\etau^k|^2}{2(\varrho_\etau^k)^2\!}\,\frac{\varrho_\etau^k-\varrho_\etau^{k-1}\!\!}\tau
\\[-.1em]&\nonumber\hspace{1em}
=\vvk\Cdot\frac{\pp_\etau^k-\pp_\etau^{k-1}\!\!}\tau-\frac{|\vvk|^2\!}2
\,\frac{\varrho_\etau^k-\varrho_\etau^{k-1}\!\!}\tau
\\[-.6em]&\nonumber\hspace{-.1em}
\stackrel{{\scriptsize\rm(\ref{Euler-small-viscoelastodyn+disc}a,b)}}{=}\!\!
\vvk\Cdot
\bigg({\rm div}\Big(\TTtauk{+}\DD_\etau^k{-}\CUT(\varrho_\etau^k)\pp_\etau^k{\otimes}\vv_\etau^k\Big)+
\varrho_\etau^k\GRAVITY_{\DELTA\tau}^k
-\DELTA|\nabla\varrho_\etau^k|^{r-2}(\nabla\vvk)\nabla\varrho_\etau^k
\\[-.2em]&\hspace{8.3em}
-\EPS|\vvk|^{p-2}\vvk
\bigg)+\frac{|\vvk|^2\!\!}{2}\,{\rm div}\Big(\CUT(\varrho_\etau^k)\pp_\etau^k
-\DELTA|\nabla\varrho_\etau^k|^{r-2}\nabla\varrho_\etau^k\Big)\,.
\label{Euler-ED-basic-convexity-calculus}
\end{align}
When integrated over $\varOmega$, we can use the calculus
\eq{calculus-for-kinetic} for
$(\varrho_\etau^k,\vv_\etau^k,\pp_\etau^k)$ instead of $(\varrho,\vv,\pp)$
and the cancellation of the regularization terms:
\begin{align}\nonumber
\int_\varOmega&\vvk\Cdot\big(|\nabla\varrho_\etau^k|^{r-2}
(\nabla\vvk)\nabla\varrho_\etau^k\big)
+\frac{|\vvk|^2\!\!}{2}\,{\rm div}\big(|\nabla\varrho_\etau^k|^{r-2}
\nabla\varrho_\etau^k\big)\,\d\xx
\\[-.4em]&=\int_\varOmega\vvk\Cdot\big(|\nabla\varrho_\etau^k|^{r-2}
(\nabla\vvk)\nabla\varrho_\etau^k\big)
-\nabla\frac{|\vvk|^2\!\!}{2}\Cdot\big(|\nabla\varrho_\etau^k|^{r-2}
\nabla\varrho_\etau^k\big)\,\d\xx=0\,.
\label{Eulerian-small-reg-cancel}\end{align}
Thus the inequality \eq{Euler-ED-basic-convexity-calculus} turns into
\begin{align}
  \!\frac1\tau&\int_\varOmega\frac{|\pp_\etau^k|^2}{2\varrho_\etau^k}
  -\frac{|\pp_\etau^{k-1}|^2}{2\varrho_\etau^{k-1}}\,\d\xx
\le\!\Big\langle{\rm div}\big(\TTtauk{+}\DD_\etau^k\big),\vvk\Big\rangle
+\!\int_\varOmega\varrho_\etau^k\GRAVITY_{\DELTA\tau}^k\Cdot\vvk
-\DELTA|\vvk|^p\,\d\xx\,.
\label{Eulerian-small-reg-cancel+}\end{align}
Furthermore, by the convexity of $\varphi$ used for
$\varphi'(\Ee_\etau^k)\Colon(\Ee_\etau^k{-}\Ee_\etau^{k-1})
\le\varphi(\Ee_\etau^k)-\varphi(\Ee_\etau^{k-1})$,
we can execute \eq{Euler-small-divT.v++} as an inequality where, additionally,
we handle the regularization terms by using the Green formula for
\begin{align}\nonumber
-\int_\varOmega\varphi'(\Ee_\etau^k)\Colon\Delta\Ee_\etau^k\,\d\xx
&=
\int_\varOmega\nabla\varphi'(\Ee_\etau^k)\Vdots\nabla\Ee_\etau^k\,\d\xx
\\&\label{Euler-small-calc-Delta}
=\int_\varOmega\varphi''(\Ee_\etau^k)\nabla\Ee_\etau^k\Vdots\nabla\Ee_\etau^k\,\d\xx
\ge\alpha\int_\varOmega|\nabla\Ee_\etau^k|^2\,\d\xx\,,
\end{align}
where $\alpha>0$ is from \eq{Euler-small-ass-phi}; here we have formally
used the 2nd-order derivative $\varphi''$ but, for $\varphi\in C^1(\Rsym)$,
this inequality holds by a smoothening argument, too.
Thus, we arrive to the discrete (and regularized) analog of the energy balance
\eq{Euler-small-energy-balance-stress}, specifically
\begin{align}\nonumber
\!\!\!&\int_\varOmega\frac{|\pp_\etau^k|^2\!}{2\varrho_\etau^k\!}
+\varphi(\Ee_\etau^k)\,\d\xx+
\tau\sum_{l=1}^k\int_\varOmega\!\bigg(\bbD\strain(\vv_\etau^l)\Colon\strain(\vv_\etau^l)
+R\,\big|{\rm dev}\,\varphi'(\Ee_\etau^l)\big|^2\!+\HYPER|\nabla^2\vv_\etau^l|^p\!
\\[-.7em]&\hspace{1em}
+\DELTA|\vv_\etau^l|^p
+\alpha\DELTA|\nabla\Ee_\etau^l|^2\bigg)\d\xx
\le\int_\varOmega\!\frac{|\pp_0|^2\!}{2\varrho_0\!}+\varphi(\Ee_0)\,\d\xx
+\tau\sum_{l=1}^k
\int_\varOmega\!
\varrho_\etau^l\GRAVITY_{\DELTA\tau}^l\Cdot\vv_\etau^l\,\d\xx\,.\!\!
\label{Euler-ED-basic-energy-balance-disc}\end{align}

To derive the basic a-priori estimates, we first test
\eq{Euler-small-viscoelastodyn+0disc} by 1. Due to the boundary condition
$\pp_\etau^k\Cdot\nn=0$, this gives the total-mass conservation 
\begin{align}\label{Euler-small-est-1+}
\int_\varOmega\varrho_\etau^k\,\d\xx=M\
\end{align}
with $M$ from \eq{Euler-small-est-1}.
We have at our disposal the energetics \eq{Euler-ED-basic-energy-balance-disc}.
When estimating the bulk-load term $\varrho_\etau^l\GRAVITY_{\DELTA\tau}^l$ as
\begin{align}
\nonumber
\!\!\int_\varOmega\varrho_\etau^l\GRAVITY_{\DELTA\tau}^l\Cdot\vv_\etau^l\,\d\xx&\le
\Bigg\|\frac{\pp_\etau^l}{\sqrt{\varrho_\etau^l}}\Bigg\|_{L^2(\varOmega;\R^3)}^{}\!\!
\big\|\!\sqrt{\varrho_\etau^l}\GRAVITY_{\DELTA\tau}^l\big\|_{L^2(\varOmega;\R^3)}^{}
\\[-.1em]&\nonumber
\le\Bigg\|\frac{\pp_\etau^l}{\sqrt{\varrho_\etau^l}}\Bigg\|_{L^2(\varOmega;\R^3)}^{}\!\!\!
\|\varrho_\etau^l\|_{L^1(\varOmega)}^{}\|\GRAVITY_{\DELTA\tau}^l\|_{L^\infty(\varOmega;\R^3)}^{}
\\[-.1em]&
\le M\Bigg(\!1+\bigg\|\frac{\pp_\etau^l}{\sqrt{\varrho_\etau^l}}\bigg\|_{L^2(\varOmega;\R^3)}^2
\Bigg)\!
\|\GRAVITY_{\DELTA\tau}^l\|_{L^\infty(\varOmega;\R^3)}^{}
\label{Euler-small-est-Gronwall}\end{align}
and taking into account the assumption \eq{Euler-small-ass-g},
we can read the estimates 
\begin{subequations}\label{Euler-small-est}\begin{align}
\label{Euler-small-est1}
&\bigg\|\frac{\overlineppetau}{\sqrt{\overlinerhoetau}}
\bigg\|_{L^\infty(I;L^2(\varOmega;\R^3))}^{}\le C
\ \ \text{ and }\ \ 
\|\overlinevvtau\|_{L^p(I\times\varOmega;\R^3)}\le\frac C{\sqrt[p]\EPS}\,,
\\&\label{Euler-small-est2}
\|\strain(\overlinevvtau)\|_{L^2(I{\times}\varOmega;\R^{3\times3})}^{}\le C
\ \ \ \,\text{ and }\ \ \ \,
\|\nabla^2\overlinevvtau\|_{L^p(I{\times}\varOmega;\R^{3\times3\times3})}^{}\le C\,,
\\&\label{Euler-small-est3}
\|\varphi(\overlineEetau)\|_{L^\infty(I;L^1(\varOmega))}^{}\le C
\ \ \ \,\text{ and }\ \ \ \,
\|\nabla\overlineEetau\|_{L^2(I{\times}\varOmega;\R^{3\times3\times3})}^{}
\le\frac{C}{\sqrt\DELTA}\,.
\intertext{Note that in the last estimate, we used
the strong convexity assumption \eq{Euler-small-ass-phi}. By the coercivity and
growth \eq{Euler-small-ass-phi} of $\varphi$, the last estimate implies}
\label{Euler-small-est3+}
&\|\overlineEetau\|_{L^\infty(I;L^2(\varOmega;\R^{3\times3}))}^{}\!\le C,\ \,
\|\overlineTetau\|_{L^\infty(I;L^1(\varOmega;\R^{3\times3}))}^{}\!\le C,\ \,
\|{\rm dev}\overlineTetau\|_{L^\infty(I;L^2(\varOmega;\R^{3\times3}))}^{}\!\le C.
\intertext{Then, from the former estimate in \eq{Euler-small-est1} and
from $\overlinerhoetau\le\rhoMAX{+}1$, we have also}
&\nonumber\|\overlineppetau\|_{L^\infty(I;L^1(\varOmega;\R^3))}^{}\le\big\|\!\sqrt{\overlinerhoetau}\big\|_{L^\infty(I;L^2(\varOmega;\R^3))}^{}\bigg\|\frac{\overlineppetau}{\sqrt{\overlinerhoetau}}\bigg\|_{L^\infty(I;L^2(\varOmega;\R^3))}^{}
\\[-.3em]&\hspace{8.2em}
={\rm meas}(\varOmega)(\rhoMAX{+}1)\bigg\|\frac{\overlineppetau}{\sqrt{\overlinerhoetau}}\bigg\|_{L^\infty(I;L^2(\varOmega;\R^3))}^{}
\le C\,\!\!\!\!\!\!\!
\label{Euler-small-est3+++}\end{align}\end{subequations}
with $C$ here and below denoting a generic constant;
it is important that $C$ in these estimates \eq{Euler-small-est} is independent
of $\DELTA$ and of $\EPS$. Here we have used the notation that exploits the interpolants
defined as follows: using the values $(\vvk)_{k=0}^{T/\tau}$, we define the
piecewise constant and the piecewise affine (forward or backward) interpolants
respectively as
\begin{align}\nonumber
&\vv_\etau(t)\!:=\Big(\frac t\tau{-}k{+}1\Big)\vvk
\!+\Big(k{-}\frac t\tau\Big)\vvkk,\ \ \
\\[-.0em]&\label{def-of-interpolants}
\overlinevvtau(t)\!:=\vvk\,,\ \ \text{ and }\ \ 
\underline\vv_\etau(t)\!:=\vvkk\ \text{ for }\ (k{-}1)\tau<t\le k\tau
\end{align}
for $k=0,1,...,T/\tau$. Analogously, we define also $\overlinerhoetau$,
$\varrho_\etau^{}$, $\overlineppetau$, $\pp_\etau^{}$, $\overlineEetau$, $\Eetau$, etc.

In terms of these interpolants, the recursive system \eq{Euler-small-viscoelastodyn+disc}
can be written ``more compactly'' as
\begin{subequations}\label{Euler-small-viscoelastodyn+discr}
\begin{align}\label{Euler-small-viscoelastodyn+0discr}
&\!\!\pdt{\varrho_\etau}={\rm div}\Big(\DELTA
|\nabla\overlinerhoetau|^{r-2}\nabla\overlinerhoetau-\CUT(\overlinerhoetau)\overlineppetau\Big)\,,
\\&\nonumber
\!\!\pdt{\pp_\etau}=
{\rm div}\Big(\overlineTetau{+}\overlineDetau
{-}\CUT(\overlinerhoetau)\overlineppetau{\otimes}\overlinevvtau\Big)
+\overlinerhoetau\overline\GRAVITY_{\DELTA\tau}
\\[-.3em]&\nonumber
\hspace*{5.5em}
-\EPS|\overlinevvtau|^{p-2}\overlinevvtau
-\DELTA|\nabla\overlinerhoetau|^{r-2}(\nabla\overlinevvtau)\nabla\overlinerhoetau
\\[-.2em]&\nonumber
\hspace*{3em}\text{ with }\
\overlineTetau=\varphi'(\overlineEetau)+\varphi(\overlineEetau)\bbI
\ \ \text{ and }\ \ \overlineDetau=\bbD\strain(\overlinevvtau)
-{\rm div}\overline{\mathfrak{H}}_\etau\,,\ \ 
\\[-.3em]&
\hspace*{3em}
\text{ where }\ \overline{\mathfrak{H}}_\etau=
\HYPER\big|\nabla^2\overlinevvtau\big|^{p-2}\nabla^2\overlinevvtau
\ \ \text{ with $\ {\overlineppetau}=\overlinerhoetau\overlinevvtau$}\,,
\label{Euler-small-viscoelastodyn+1discr}
\\[-.3em]&
\!\!\pdt{\Ee_\etau}
=\strain(\overlinevvtau)-\bm B_\text{\sc zj}^{}(\overlinevvtau,\overlineEetau)
-R\,{\rm dev}\overlineTetau+\DELTA\Delta\overlineEetau
\label{Euler-small-viscoelastodyn+2discr}
\end{align}\end{subequations}
with the corresponding boundary conditions and with the initial conditions
for $(\varrho_\etau,\pp_\etau,\Ee_\etau)$.

Except the latter estimate in \eq{Euler-small-est2}, all the above estimates
show the stability of the discrete regularized scheme even for the
simple-material variant with $\HYPER=0$. However, for the rest, we will need the
non-simple material model with a fixed hyper-viscosity coeffcient $\HYPER>0$
and $p>3$.

\medskip{\it Step 3: Further a-priori estimates relying on $\HYPER>0$}.
Now we can test \eq{Euler-small-viscoelastodyn+0disc} by $\varrho_\etau^k$.
Using the algebra
$(\varrho_\etau^k{-}\varrho_\etau^{k-1})\varrho_\etau^k
=\frac12|\varrho_\etau^k|^2-\frac12|\varrho_\etau^{k-1}|^2
+\frac12|\varrho_\etau^k{-}\varrho_\etau^{k-1}|^2$ and Young's inequality, we obtain
\begin{align}\nonumber
&\frac12\|\varrho_\etau^k\|_{L^2(\varOmega)}^2-\frac12\|\varrho_\etau^{k-1}\|_{L^2(\varOmega)}^2
+\DELTA\!\int_\varOmega|\nabla\varrho_\etau^k|^{r}
\,\d\xx
\le\int_\varOmega\CUT(\varrho_\etau^k)\pp_\etau^k\Cdot\nabla\varrho_\etau^k\,\d\xx
\\[-.3em]&
\qquad\ =\int_\varOmega\hspace{-.7em}
\lineunder{\CUT(\varrho_\etau^k)\varrho_\etau^k}{$\ \ \ \le\rhoMAX{+}1$}\hspace{-.7em}
\vvk\Cdot\nabla\varrho_\etau^k\,\d\xx
\le C_{\DELTA,r,p}^{}\big(1{+}\|\vvk\|_{L^p(\varOmega;\R^3)}^p\big)
+\frac\DELTA2\|\nabla\varrho_\etau^k\|_{L^r(\varOmega;\R^3)}^r
\,.
\label{Euler-small-estimation-cont-eq}\\[-2em]\nonumber
\end{align}
Due to (\ref{Euler-small-est}a,b), we obtain
$\|\overlinerhoetau\|_{L^\infty(I;L^2(\varOmega))}\le C_\EEps$ and
$\|\nabla\overlinerhoetau\|_{L^r(I\times\varOmega;\R^3)}\le C_\EEps$,
although the first estimate is granted by the already proved
bounds $0\le\overlinerhoetau\le\rhoMAX{+}1$. 
These estimates can further be improved by using the
convexity of the potential of the quasilinear operator
$\varrho\mapsto-{\rm div}(|\nabla\varrho|^{r-2}\nabla\varrho)$
when applying $\nabla$-operator to 
\eq{Euler-small-viscoelastodyn+0disc} and testing it by
$|\nabla\overlinerhoetau|^{r-2}\nabla\overlinerhoetau$ still allows,
in view of (\ref{Euler-small-est}a,b), for the estimate
\begin{align}
&\nonumber
\frac1r\|\nabla\varrho_\etau^k\|_{L^r(\varOmega;\R^3)}^r
-\frac1r\|\nabla\varrho_\etau^{k-1}\|_{L^r(\varOmega;\R^3)}^r
+\DELTA\!\int_\varOmega\big|{\rm div}
\big(|\nabla\varrho_\etau^k|^{r-2}\nabla\varrho_\etau^k\big)\big|^2\,\d\xx
\\[-.4em]&\nonumber\qquad\le
\int_\varOmega{\rm div}\big(\CUT(\varrho_\etau^k)\pp_\etau^k\big)\Cdot
{\rm div}\big(|\nabla\varrho_\etau^k|^{r-2}\nabla\varrho_\etau^k\big)\,\d\xx
\\[-.1em]&\nonumber\qquad
\le\frac1{2\DELTA}\big\|{\rm div}\big(\CUT(\varrho_\etau^k)\pp_\etau^k\big)\big\|_{L^2(\varOmega)}^2
+\frac\DELTA2\big\|{\rm div}\big(|\nabla\varrho_\etau^k|^{r-2}\nabla\varrho_\etau^k\big\|_{L^2(\varOmega)}^2
\\[-.1em]&\nonumber\qquad\le
\frac1{\DELTA}\Big(\big\|\CUT(\varrho_\etau^k)\varrho_\etau^k{\rm div}\,\vvk\big\|_{L^2(\varOmega)}^2\!
+\big\|\big(\CUT(\varrho_\etau^k){+}\varrho_\etau^k\CUT'(\varrho_\etau^k)\big)\nabla\varrho_\etau^k\Cdot\vvk\big\|_{L^2(\varOmega)}^2\Big)
\\[-.1em]&\label{Euler-small-est10++calcul}\hspace{22em}
+\frac\DELTA2\big\|{\rm div}\big(|\nabla\varrho_\etau^k|^{r-2}\nabla\varrho_\etau^k\big)\big\|_{L^2(\varOmega)}^2\,.
\end{align}
This yields
\begin{subequations}\label{Euler-small-est+}
\begin{align}
&\qquad\|\overlinerhoetau\|_{L^\infty(I;W^{1,r}(\varOmega))}^{}\le C_\EEps
\ \ \ \ \text{ and }\ \ \ \
\big\|{\rm div}(|\nabla\overlinerhoetau|^{r-2}\nabla\overlinerhoetau)
\big\|_{L^2(I\times\varOmega)}^{}\le C_\EEps\,.
\label{Euler-small-est10++}
\intertext{Realizing that $\nabla\overlineppetau=\overlinerhoetau\nabla\overlinevvtau +\nabla\overlinerhoetau{\otimes}\overlinevvtau$, this allows for the estimate}
&\qquad\label{Euler-small-est10+}
\|\overlineppetau\|_{L^p(I;W^{1,r}(\varOmega;\R^3))}^{}\le C_\EEps\,. 
\end{align}\end{subequations}

The existence of a weak solution $(\varrho_\etau^k,\pp_\etau^k,\Ee_\etau^k)\in
W^{1,r}(\varOmega)\times W^{2,p}(\varOmega;\R^3)\times H^1(\varOmega;\Rsym)$
with $\varrho_\etau^k>0$ of the coupled quasi-linear boundary-value problem
\eq{Euler-small-viscoelastodyn+disc}--\eq{BC-disc} can thus be proved by a combination
of the quasilinear technique for \eqref{Euler-small-viscoelastodyn+0disc} and
\eqref{Euler-small-viscoelastodyn+1disc} involving the quasi-linear strictly
monotone terms $-{\rm div}(\DELTA|\nabla\varrho|^{r-2}\nabla\varrho)$ and
${\rm div}^2(\HYPER|\nabla^2\vv|^{p-2}\nabla^2\vv)+\EPS|\vv|^{p-2}\vv$, using the usual
semi-linear technique for the rest of the system, provided $\varrho_\etau^{k-1}>0$,
$\pp_\etau^{k-1}$, and $\Ee_\etau^{k-1}$ are known from the previous time step. 
The mentioned coercivity is a particular consequence of the a-priori estimates
\eq{Euler-small-est}. Let us note that, due to the convective terms,
this system does not have any potential so the rather nonconstructive Brouwer
fixed-point arguments combined with the Galerkin approximation are to be used.

By the comparison from \eq{Euler-small-viscoelastodyn+0discr}, using
\eq{Euler-small-est10+}, we obtain
\begin{subequations}\label{Euler-small-est++}\begin{align}
&
\Big\|\pdt{\varrho_\etau}\Big\|_{L^2(I\times\varOmega)}\!
\le C_\EEps\,. 
\label{Euler-small-est6}
\intertext{Moreover, by the comparison from
\eq{Euler-small-viscoelastodyn+1discr} when using (\ref{Euler-small-est}a,b,d), we
obtain also} 
\nonumber
&\Big\|\pdt{\pp_\etau}\Big\|_{L^{p'}(I;W^{2,p}(\varOmega;\R^3)^*)}^{}\!
=\!\!\sup_{\|\wt\vv\|_{L^p(I;W^{2,p}(\varOmega;\R^3))}\le1}\int_0^T\!\!\!\int_\varOmega
\Big(\overlineTetau{+}\bbD\strain(\overlinevvtau)
{-}\CUT(\overlinerhoetau)\overlineppetau{\otimes}\overlinevvtau\Big)
\Colon\strain(\wt\vv)\!\!\!
\\[-.2em]&\hspace{2.3em}
+\overline{\mathfrak{H}}_\etau\Vdots\nabla^2\wt\vv
+\Big(\EPS|\overlinevvtau|^{p-2}\overlinevvtau\!
{-}\DELTA|\nabla\overlinerhoetau|^{r-2}(\nabla\overlinevvtau)\nabla\overlinerhoetau\!
{-}\overlinerhoetau\overline\GRAVITY_{\DELTA\tau}\Big)\Cdot\wt\vv
\,\d\xx\d t\le C_\EEps\,;\label{Euler-small-est6+}
\intertext{note that
$\CUT(\overlinerhoetau)\overlineppetau{\otimes}\overlinevvtau
\in L^{p/2}(I;L^\infty(\varOmega;\Rsym))$ is indeed in duality
with $\strain(\wt\vv)\in L^p(I;L^\infty(\varOmega;\Rsym))$ as $p\ge3$ is assumed.
Eventually, by the comparison from \eq{Euler-small-viscoelastodyn+2discr},
we obtain also}\nonumber
&\Big\|\pdt{\Ee_\etau}\Big\|_{L^{6/5}(I;H^1(\varOmega;\R^{3\times3})^*)}^{}\!
=\!\!\sup_{\|\wt\Ee\|_{L^{6}(I;H^1(\varOmega;\R^{3\times3}))}\le1}\int_0^T\!\!\!\int_\varOmega
\DELTA\nabla\overlineEetau\Vdots\nabla\wt\Ee
\\[-.2em]&\hspace{13.5em}
+\big(\bm B_\text{\sc zj}^{}(\overlinevvtau,\overlineEetau)
{+}R\,{\rm dev}\,\overlineTetau{-}\strain(\overlinevvtau)\big)
\Colon\wt\Ee\,\d\xx\d t\le C_\EEps\,;
\label{Euler-small-est7}\end{align}\end{subequations}
here we used the estimates of 
$(\overlinevvtau{\cdot}\nabla)\overlineEetau\in L^{6/5}(I;L^2(\varOmega;\Rsym))$
if $p\ge3$ and that ${\rm skw}(\nabla\overlinevvtau)\overlineEetau$
$-\overlineEetau{\rm skw}(\nabla\overlinevvtau)\in L^p(I;L^2(\varOmega;\Rsym))$,
and also the bound \eq{Euler-small-est3+} for
${\rm dev}\,\overlineTetau={\rm dev}\varphi'(\overlineEetau)$.
The subscripts $\EPS$ and $\DELTA$ in $C_\EEps$ in
\eq{Euler-small-est+}--\eq{Euler-small-est++} want to indicate the dependence of these
estimates on $\EPS$ in the latter estimate in \eq{Euler-small-est1} and $\DELTA$ in the
latter estimate in \eq{Euler-small-est3} or through \eq{Euler-small-estimation-cont-eq}
and \eq{Euler-small-est10++calcul} which also contain $\DELTA$.

\medskip{\it Step 4: Limit passage for $\tau\to0$}.
By the Banach selection principle, we obtain a subsequence
converging weakly* with respect to the topologies indicated in
\eq{Euler-small-est} and \eq{Euler-small-est+} to some limit
$(\varrho_\EEps,\pp_\EEps,\vv_\EEps,\Ee_\EEps)$. Specifically,
\begin{subequations}\label{Euler-small-converge}
\begin{align}
&&&\!\!\overlinerhoetau\to\varrho_\EEps&&\text{weakly* in $\ L^\infty(I;W^{1,r}(\varOmega))$\,,}
\\&&&\!\!\varrho_\etau\to\varrho_\EEps&&\text{weakly* in $\ L^\infty(I;W^{1,r}(\varOmega))\,\cap\, 
H^1(I;L^2(\varOmega))$\,,}
\\&&&\label{Euler-small-converge-bar-p}
\!\!\overlineppetau\to\pp_\EEps&&\text{weakly\ \;in $\
L^p(I;W^{1,r}(\varOmega;\R^3))$\,,}
\\&&&\label{Euler-small-converge-p}
\!\!\pp_\etau\to\pp_\EEps&&\text{weakly\ \;in $\
L^p(I;W^{1,r}(\varOmega;\R^3))\,\cap\,W^{1,p'}(I;W^{2,p}(\varOmega;\R^3)^*)$\,,}
\\&&&\label{Euler-small-converge-bar-v}
\!\!\overlinevvtau\to\vv_\EEps&&\text{weakly\ \;in $\ L^p(I;W^{2,p}(\varOmega;\R^3))$\,,}
\\
&&&\!\!\overlineEetau\to\Ee_\EEps\!\!\!&&\text{weakly* in $\ 
L^2(I;H^1(\varOmega;\Rsym))\,\cap\,L^\infty(I;L^2(\varOmega;\Rsym))$\,,}\!\!
\\
&&&\!\!\Ee_\etau\to\Ee_\EEps\!\!\!&&\text{weakly\ \;in $\ 
L^2(I;H^1(\varOmega;\Rsym))\,\cap\,W^{1,6/5}(I;H^1(\varOmega;\Rsym)^*)$\,.}&&
\end{align}\end{subequations}
Notably, the limit of $\overlinerhoetau$ and $\varrho_\etau$ is indeed the same
due to the control of $\pdt{}\varrho_\etau$ in \eq{Euler-small-est6}; cf.\
\cite[Sect.8.2]{Roub13NPDE}. The same is true also for $\overlineppetau$ and
$\pp_\etau$, 
and for $\overlineEetau$ and $\Ee_\etau$, too.

By the compact embedding $W^{1,r}(\varOmega)\subset C(\barOmega)$ and the
(generalized) Aubin-Lions theorem,  cf.\ \cite[Corollary~7.9]{Roub13NPDE}
relying on the compact embedding $W^{1,r}(\varOmega)\subset
C(\barOmega)$ for $r>3$,
we have also 
\begin{subequations}\label{Euler-small-converge-strong}
\begin{align}\label{Euler-small-converge-strong-rho}
&&&\!\!\overlinerhoetau\to\varrho_\EEps\!\!\!&&\text{strongly in }\
L^a(I;C(\barOmega))\ 
\ \text{ for any }\ 1\le a<\infty\,,&&
\\\label{Euler-small-converge-strong-p-}
&&&\!\!\overlineppetau\to\pp_\EEps\!\!\!\!\!\!&&\text{strongly in }\ L^p(I;C(\barOmega;\R^3))\,,
\\\label{Euler-small-converge-strong-p}
&&&\!\!\CUT(\overlinerhoetau)\overlineppetau\to
\CUT(\varrho_\EEps)\pp_\EEps\!\!\!\!\!\!&&\text{strongly in }\ L^p(I;C(\barOmega;\R^3))\,,
\ \text{ and}
\\\label{Euler-small-converge-strong-E}
&&&\!\!\overlineEetau\to\Ee_\EEps\!\!\!\!\!&&\text{strongly in }\
L^a(I;L^2(\varOmega;\Rsym))\ \text{ for any }\ 1\le a<\infty\,.
\intertext{Moreover, using the Arzel\`a-Ascoli-type theorem, cf.\
\cite[Lemma~7.10]{Roub13NPDE}, we have also}
\label{Euler-small-converge-strong-rho+}
&&&\!\!\varrho_\etau\to\varrho_\EEps\!\!\!&&\text{strongly in }\ C(I{\times}\barOmega)\,.
\end{align}\end{subequations}

By \eq{Euler-small-converge-strong-E}, we can pass
to the limit in some nonlinear terms. Specifically, denoting
$\TT_{\!\EEps}:=\varphi'(\Ee_\EEps)+\varphi(\Ee_\EEps)\bbI$, we have
\begin{subequations}\label{Euler-small-converge-strong+}
\begin{align}\label{Euler-small-converge-strong-T}
&\overlineTetau=\varphi'(\overlineEetau){+}\varphi(\overlineEetau)\bbI\to
\TT_{\!\EEps}
&&\hspace{-10em}\text{strongly\ in }\, 
L^a(I;L^1(\varOmega;\Rsym))\ \text{ for any }1\le a<\infty,
\\&{\rm dev}\,\overlineTetau\!={\rm dev}\varphi'(\overlineEetau)
\to{\rm dev}\,\TT_{\!\EEps}
&&\hspace{-10em}\text{strongly\ in }\,
L^a(I;L^2(\varOmega;\Rsym))\ \text{ for any }1\le a<\infty,\!
\\&\bm B_\text{\sc zj}^{}(\overlinevvtau,\overlineEetau)
\to\bm B_\text{\sc zj}^{}(\vv_\EEps,\bm E_\EEps)\!\!\!
&&\hspace{-10em}\text{weakly in }\ L^{6/5}
(L^2(\varOmega;\R^{3\times3}))\,;
\intertext{notably, the last convergence is only weak due to the term
$\nabla\overlineEetau$ contained in
$\bm B_\text{\sc zj}^{}(\overlinevvtau,\overlineEetau)$.
More in detail, exploiting \eq{Euler-small-converge-strong-E}, for any
$\wt\Ee\in L^6(I;L^2(\varOmega;\Rsym))$, it holds}
&\nonumber
\int_\varOmega\big((\overlinevvtau\Cdot\nabla)\overlineEetau\big)\Colon\wt\Ee\,\d\xx=
-\int_\varOmega({\rm div}\overlinevvtau)\overlineEetau\Colon\wt\Ee
+\overlineEetau\Colon\big((\overlinevvtau\Cdot\nabla)\wt\Ee\big)\,\d\xx\hspace{-25em}
\\[-.4em]&\nonumber\hspace{6em}
\to-\int_\varOmega({\rm div}\vv_\EEps)\Ee_\EEps\Colon\wt\Ee
+\Ee_\EEps\Colon\big((\vv_\EEps\Cdot\nabla)\wt\Ee\big)\,\d\xx=\int_\varOmega\!
\big((\vv_\EEps\Cdot\nabla)\Ee_\EEps\big)\Colon\wt\Ee\,\d\xx\,.\hspace{-25em}
\intertext{Then, we will prove that}
&\overlinevvtau=\frac{\overlineppetau}{\overlinerhoetau}\to\frac{\pp_\EEps}{\varrho_\EEps}=\vv_\EEps\!\!\!
&&
\label{Euler-small-converge-strong-vv}
\hspace{-10em}\text{strongly in $\ L^q
(I{\times}\varOmega;\R^3)\,$ for any $1\le q<p$\,,}
\intertext{and thus,
by \eq{Euler-small-converge-strong-p}, also}
&\CUT(\overlinerhoetau)\overlineppetau{\otimes}\overlinevvtau\to
\CUT(\varrho_\EEps)\pp_\EEps{\otimes}\vv_\EEps
=\CUT(\varrho_\EEps)\varrho_\EEps\vv_\EEps{\otimes}\vv_\EEps\!
&&\text{strongly in }\ L^{q/2}(I{\times}\varOmega;\Rsym)\,.
\label{Euler-small-converge-strong-pxv}\end{align}\end{subequations}
For this, we have to use a tiny argument: from
(\ref{Euler-small-converge-strong}a,b), for a subsequence of $\tau\to0$, we
know that $\overlineppetau/\overlinerhoetau$ converges to some $\wt\vv_\EEps$
a.e.\ on $I{\times}\varOmega$. Simultaneously, we know
$\sup_{\tau>0}\int_0^T\int_\varOmega|\wt\vv_\EEps-\overlineppetau/\overlinerhoetau|^p\,\d\xx\d t<\infty$
with $p>1$ so that, by the de\,la\,Vall\'ee\,Poussin theorem,
$\{|\wt\vv_\EEps{-}\overlineppetau/\overlinerhoetau|^{q}\}_{\tau>0}$ is relatively
weakly compact in $L^1(I{\times}\varOmega)$ for any $1\le q<p$. Then, by the
Dunford-Pettis theorem, it is uniformly integrable and, by the Vitali theorem,
it converges a.e.\ to its limit which equals 0, i.e.\ it converges to 0 in
$L^q(I{\times}\varOmega)$. By \eq{Euler-small-converge-bar-v}, we can identify
$\wt\vv_\EEps=\vv_\EEps$ so that all the (already chosen subsequence for
\eq{Euler-small-converge}) converges to $\vv_\EEps$, which proves
\eq{Euler-small-converge-strong-vv}. Actually, by the interpolation with
$L^p(I;W^{2,p}(\varOmega;\R^3))$, we have the strong convergence
\eq{Euler-small-converge-strong-vv} even in a better space
$L^q(I;C(\barOmega;\R^3))$.

Now we can make a limit passage in \eq{Euler-small-viscoelastodyn+discr}.
First, we can pass to the limit in the quasilinear equation
\eq{Euler-small-viscoelastodyn+0discr} just by the weak convergence using
monotonicity in the main part and so-called Minty's trick and the strong
convergence \eq{Euler-small-converge-strong-p}, obtaining 
the equation $\pdt{}\varrho_\EEps=
{\rm div}(\DELTA|\nabla\varrho_\EEps|^{r-2}\nabla\varrho_\EEps-\CUT(\varrho_\EEps)\pp_\EEps)$
holding in the weak sense involving the respective initial and the boundary
conditions. We now prove the strong convergence of $\nabla\overlinerhoetau$ by
subtracting this equations from \eq{Euler-small-viscoelastodyn+0discr} and test
it by $\varrho_\EEps{-}\overlinerhoetau$. Like \eq{Euler-small-estimation-cont-eq},
after integration over $[0,t]$, we can estimate 
\begin{align}\nonumber
&\frac12\|\varrho_\EEps(t){-}\varrho_\etau(t)\|_{L^2(\varOmega)}^2
+\DELTA\int_0^t\!\!\int_\varOmega c_r|\nabla\varrho_\EEps{-}\nabla\overlinerhoetau|^r
\,\d\xx\d t
\\[-.3em]&\nonumber\qquad
\le
\int_0^t\!\!\int_\varOmega\big(\CUT(\varrho_\EEps)\pp_\EEps{-}\CUT(\overlinerhoetau)\overlineppetau
\big)\Cdot\nabla(\varrho_\EEps{-}\overlinerhoetau)
\,\d\xx\d t
\\[-.3em]&\nonumber\qquad
\le
C_{\EEps,r}
\big\|\CUT(\varrho_\EEps)\pp_\EEps{-}\CUT(\overlinerhoetau)\overlineppetau\big\|_{L^{r'}(I\times\varOmega;\R^3)}^{r'}\!
+\frac{\DELTA}2\!\int_0^t\!\!\int_\varOmega\! c_r|\nabla\varrho_\EEps{-}\nabla\overlinerhoetau|^r\d\xx\,,
\end{align}
cf.\ also the estimation in \cite[Remark 8.11]{Roub13NPDE}; the constant
$c_r>0$ is from the inequality
$c_r|g-\wt g|^r\le(|g|^{r-2}g-|\wt g|^{r-2}\wt g)\Cdot(g-\wt g)$ which
holds for $r\ge2$, cf.\ \cite[Chap.3, Example~1.7]{HuPap97HMAT}.
Thus, reminding $r>3$ and $p>3$ (so that surely $p>r'$) and using
(\ref{Euler-small-converge-strong}a,b), we obtain
\begin{align}\label{Euler-small-nabla-rho-strong}
\overlinerhoetau\to\varrho_\EEps\quad
\text{strongly in $\ L^\infty(I;L^2(\varOmega))\,\cap\,L^r(I;W^{1,r}(\varOmega))$\,.}
\end{align}

The quasi-linear momentum equation \eq{Euler-small-viscoelastodyn+1discr} is
more difficult. We also use the weak convergence and the Minty trick for
the limit passage in the quasi-linear monotone terms
${\rm div}^2(
\HYPER|\nabla^2\vv|^{p-2}\nabla^2\vv)+\EPS|\vv|^{p-2}\vv$. For any
$\wt\vv\in L^p(I;W^{2,p}(\varOmega;\R^3))$, this monotonicity together with
\eq{Euler-small-viscoelastodyn+1discr} tested by $\overlinevvtau{-}\wt\vv$
yields 
\begin{align}\nonumber
&\nonumber\hspace*{0em}
0\le\int_0^T\!\!\!\int_\varOmega\!\bigg(\bbD\strain(\overlinevvtau{-}\wt\vv)
\Colon\strain(\overlinevvtau{-}\wt\vv)
+\EPS\Big(|\overlinevvtau|^{p-2}\overlinevvtau
-|\wt\vv|^{p-2}\wt\vv\Big)\Cdot(\overlinevvtau{-}\wt\vv)
  \\[-.6em]&\hspace*{10em}\nonumber
 +\HYPER\Big(|\nabla^2\overlinevvtau|^{p-2}\nabla^2\overlinevvtau
-|\nabla^2\wt\vv|^{p-2}\nabla^2\wt\vv\Big)\Vdots\nabla^2(\overlinevvtau{-}\wt\vv)
\bigg)\,\d\xx\d t
 \\[-.5em]&=\nonumber
 \int_0^T\!\!\!\int_\varOmega\bigg(
 \Big(\overlinerhoetau\overline\GRAVITY_{\DELTA\tau}
-\DELTA|\nabla\overlinerhoetau|^{r-2}(\nabla\overlinevvtau)\nabla\overlinerhoetau
-\pdt{\pp_\etau}\Big)
 \Cdot(\overlinevvtau{-}\wt\vv)
\\[-.2em]&\nonumber\hspace{5em}
-\Big(\overlineTetau{-}\CUT(\overlinerhoetau)\overlineppetau{\otimes}\overlinevvtau\Big)\Colon\strain(\overlinevvtau{-}\wt\vv)
  -\bbD\strain(\wt\vv)\Colon\strain(\overlinevvtau{-}\wt\vv)
\\[-.2em]&\nonumber\hspace{10em}
+\EPS|\wt\vv|^{p-2}\wt\vv\Cdot(\overlinevvtau{-}\wt\vv)
-\HYPER|\nabla^2\wt\vv|^{p-2}\nabla^2\wt\vv\Vdots\nabla^2(\overlinevvtau{-}\wt\vv)
\bigg)\,\d\xx\d t
\\[-.5em]&\le\nonumber
 \int_0^T\!\!\!\int_\varOmega\bigg(
 \Big(\overlinerhoetau\overline\GRAVITY_{\DELTA\tau}
-\DELTA|\nabla\overlinerhoetau|^{r-2}(\nabla\overlinevvtau)\nabla\overlinerhoetau\Big)
 \Cdot(\overlinevvtau{-}\wt\vv)
 +\pdt{\pp_\etau}\Cdot\wt\vv
\\[-.2em]&\nonumber\hspace{3em}
-\CUT(\overlinerhoetau)(\overlineppetau{\otimes}\overlinevvtau)\Colon\strain(\wt\vv)
-\big(\bbD\strain(\wt\vv){+}\overlineTetau\big)\Colon\strain(\overlinevvtau{-}\wt\vv)
-\EPS|\wt\vv|^{p-2}\wt\vv\Cdot(\overlinevvtau{-}\wt\vv)
\\[-.1em]&\hspace{10em}
-\HYPER|\nabla^2\wt\vv|^{p-2}\nabla^2\wt\vv\Vdots\nabla^2(\overlinevvtau{-}\wt\vv)
\bigg)\,\d\xx\d t+
\int_\varOmega\frac{|\pp_0|^2}{2\varrho_0}-\frac{|\pp_\etau(T)|^2}{2\varrho_\etau(T)}\,\d\xx\,,
\label{Euler-small-strong-hyper+}\end{align}
where the last inequality has again exploited the convexity of the
kinetic energy in the calculus:
\begin{align}\nonumber
\int_\varOmega&\frac{|\pp_\etau(T)|^2\!}{2\varrho_\etau(T)}
-\frac{|\pp_0|^2}{2\varrho_0}\,\d\xx{\buildrel{\eq{Euler-ED-basic-convexity-calculus}}\over{\le}}
\int_0^T\!\!\!\int_\varOmega\!\pdt{\pp_\etau\!}\Cdot\overlinevvtau
-\frac{|\overlinevvtau|^2}2\pdt{\varrho_\EEps}\,\d\xx\d t
\\&\nonumber
{\buildrel{\eq{Euler-small-viscoelastodyn+0discr}}\over{=}}
\int_0^T\!\!\!\int_\varOmega\!\pdt{\pp_\etau\!}\Cdot\overlinevvtau
+\frac{|\overlinevvtau|^2}2{\rm div}\big(\CUT(\overlinerhoetau)\overlineppetau
-\DELTA|\nabla\overlinerhoetau|^{r-2}\nabla\overlinerhoetau\big)\,\d\xx\d t
\\&\nonumber{\buildrel{\eq{Eulerian-small-reg-cancel}}\over{=}}
\int_0^T\!\!\!\int_\varOmega\!\pdt{\pp_\etau\!}\Cdot\overlinevvtau
+\frac{|\overlinevvtau|^2}2{\rm div}\big(\CUT(\overlinerhoetau\overlineppetau\big)
+\DELTA\overlinevvtau\Cdot\Big(
|\nabla\overlinerhoetau|^{r-2}(\nabla\overlinevvtau)\nabla\overlinerhoetau\Big)
\,\d\xx\d t
\\&\nonumber{\buildrel{\eq{calculus-for-kinetic}}\over{=}}
\int_0^T\!\!\!\int_\varOmega\!\bigg(\pdt{\pp_\etau\!}\Cdot\overlinevvtau
+\overlinevvtau\Cdot{\rm div}
\big(\CUT(\overlinerhoetau)\overlineppetau{\otimes}\overlinevvtau\big)
\\[-.7em]&\hspace{18em}
+\DELTA\overlinevvtau\Cdot\Big(
|\nabla\overlinerhoetau|^{r-2}(\nabla\overlinevvtau)\nabla\overlinerhoetau\Big)
\bigg)\,\d\xx\d t\,.
\label{Euler-small-strong-hyper-calc}\end{align}
Now we want to pass to the limit in \eq{Euler-small-strong-hyper+} or, more
precisely, to estimate the limit superior from above. For this, we again use the
kinetic-energy convexity, which causes the weak lower semicontinuity of
$(\varrho,\pp)\mapsto\int_\varOmega|\pp|^2/\varrho\,\d\xx$ as a convex functional
$\{\rho{\in}L^1(\varOmega);\,\rho\ge0\}\times L^1(\varOmega;\R^3)\to[0,+\infty]$.
Here we rely also on that $|\pp_\etau(T)|^2/\varrho_\etau(T)$ is bounded in
$L^1(\varOmega)$ due to the former estimate in \eq{Euler-small-est1} and on that
$\varrho_\etau(T)\to\varrho_\EEps(T)$ even strongly in $C(\barOmega)$ due to
\eq{Euler-small-converge-strong-rho+}, and on that $\pp_\etau(T)$ converges
weakly* in $C(\barOmega)^*$, i.e.\ as measures on $\barOmega$ due to
\eq{Euler-small-est3+++} to its limit which is $\pp_\EEps(T)$ because
simultaneously $\pp_\etau(T)\to\pp_\EEps(T)$ weakly in $W^{2,p}(\varOmega;\R^3)^*$
due to \eq{Euler-small-converge-p}.
Beside, we use the convergence \eq{Euler-small-converge-strong-T}, 
 \eq{Euler-small-converge-strong-pxv}, and, exploiting
 \eq{Euler-small-nabla-rho-strong}, also the convergence
\begin{align}
|\nabla\overlinerhoetau|^{r-2}(\nabla\overlinevvtau)\nabla\overlinerhoetau\to
|\nabla\varrho_\EEps|^{r-2}(\nabla\vv_\EEps)\nabla\varrho_\EEps\ \ \text{ weakly in }\
L^{p}(I;L^{r'}(\varOmega;\R^3))\,.
\end{align}
All this allows us to estimate of the limit superior of
\eq{Euler-small-strong-hyper+} from above:
\begin{align}
\nonumber
0&\le\int_0^T\!\!\Big\langle\pdt{\pp_\EEps},\wt\vv\Big\rangle
+\!\int_\varOmega\bigg(\Big(\varrho_\EEps\GRAVITY_\EEps
-\DELTA|\nabla\varrho_\EEps|^{r-2}(\nabla\vv_\EEps)\nabla\varrho_\EEps\Big)
 \Cdot(\vv_\EEps{-}\wt\vv)
\\[-.1em]&\nonumber\hspace{3em}
-\CUT(\varrho_\EEps)(\pp_\EEps{\otimes}\vv_\EEps)\Colon\strain(\wt\vv)
  -\big(\TT_\EEps{+}\bbD\strain(\wt\vv)\big)
  \Colon\strain(\vv_\EEps{-}\wt\vv)
  -\EPS|\wt\vv|^{p-2}\wt\vv\Cdot(\vv_\EEps{-}\wt\vv)
\\[-.1em]&\nonumber\hspace{8em}
-\HYPER|\nabla^2\wt\vv|^{p-2}\nabla^2\wt\vv\Vdots\nabla^2(\vv_\EEps{-}\wt\vv)
\bigg)\,\d\xx\d t+
\int_\varOmega\frac{|\pp_0|^2}{2\varrho_0}-\frac{|\pp_\EEps(T)|^2}{2\varrho_\EEps(T)}\,\d\xx
\\[-.5em]&=\nonumber
 \int_0^T\!\!\Big\langle\pdt{\pp_\EEps},\wt\vv-\vv_\EEps\Big\rangle
 +\int_\varOmega\bigg(
 \Big(\varrho_\EEps\GRAVITY_\EEps
-\DELTA|\nabla\varrho_\EEps|^{r-2}(\nabla\vv_\EEps)\nabla\varrho_\EEps\Big)
 \Cdot(\vv_\EEps{-}\wt\vv)
\\[-.3em]&\nonumber\hspace{5em}
-\Big(\TT_\EEps{-}\CUT(\varrho_\EEps)\pp_\EEps{\otimes}\vv_\EEps\Big)
\Colon\strain(\vv_\EEps{-}\wt\vv)-\bbD\strain(\wt\vv)\Colon\strain(\vv_\EEps{-}\wt\vv)
\\[-.4em]&\hspace{10em}
+\EPS|\wt\vv|^{p-2}\wt\vv\Cdot(\vv_\EEps{-}\wt\vv)
-\HYPER|\nabla^2\wt\vv|^{p-2}\nabla^2\wt\vv\Vdots\nabla^2(\vv_\EEps{-}\wt\vv)
\bigg)\,\d\xx\d t\,,
\label{Euler-small-strong-hyper+++}\end{align}
where $\langle\cdot,\cdot\rangle$ denotes the duality between
$W^{2,p}(\varOmega;\R^3)^*$ and $W^{2,p}(\varOmega;\R^3)$ and where, for
the last equality, we used the calculus like \eq{Euler-small-strong-hyper-calc}
but for the continuous-in-time limit which holds as an equality. To complete the
Minty-trick arguments, we replace $\wt\vv$ by $\vv_\EEps\pm a\wt\vv$ with some
$a>0$ and then pass $a\to0+$. Thus we prove that $\vv_\EEps$ satisfies the
regularized momentum equation \eq{Euler-small-viscoelastodyn+1discr} in the
weak sense.

In addition, we can easily pass to the limit in the semi-linear equation
\eq{Euler-small-viscoelastodyn+2discr}.
Altogether, we have proved that $(\varrho_\EEps,\vv_\EEps,\Ee_\EEps)$
is a weak solution to the $(\EPS,\DELTA)$-regularized system
\begin{subequations}\label{Euler-small-viscoelastodyn+eps}
\begin{align}\label{Euler-small-viscoelastodyn+0eps}
&\!\pdt{\varrho_\EEps}={\rm div}\big(
\DELTA|\nabla\varrho_\EEps|^{r-2}\nabla\varrho_\EEps
-\CUT(\varrho_\EEps)\varrho_\EEps\vv_\EEps\big)\,,
\\&\nonumber
\!\pdt{}\big(\varrho_\EEps\vv_\EEps\big)={\rm div}\Big(\TT_\EEps{+}\DD_\EEps
{-}\CUT(\varrho_\EEps)\varrho_\EEps\vv_\EEps{\otimes}\vv_\EEps\Big)
+\varrho_\EEps\GRAVITY_\EEps\!
-\EPS|\vv_\EEps|^{p-2}\vv_\EEps\!-\DELTA
|\nabla\varrho_\EEps|^{r-2}(\nabla\vv_\EEps)\nabla\varrho_\EEps
\\[-.2em]&\nonumber
\hspace*{9em}\text{ with }\ \TT_\EEps=\varphi'(\Ee_\EEps)+\varphi(\Ee_\EEps)\bbI
\\[-.3em]&
\hspace*{9em}
\text{ and }\ \DD_\EEps=\bbD\strain(\vv_\EEps)
-{\rm div}\mathfrak{H}_\EEps\,,\ \ 
\text{ where }\ \mathfrak{H}_\EEps=\HYPER\big|\nabla^2\vv_\EEps\big|^{p-2}
\nabla^2\vv_\EEps
\,,
\label{Euler-small-viscoelastodyn+1eps}
\\[-.3em]&
\!\pdt{\Ee_\EEps}
=\strain(\vv_\EEps)+\bm B_\text{\sc zj}^{}(\vv_\EEps,\Ee_\EEps)
-R\,{\rm dev}\,\TT_\EEps+\DELTA\Delta\Ee_\EEps
\label{Euler-small-viscoelastodyn+2eps}
\end{align}\end{subequations}
with the corresponding boundary conditions as the continuous analog of
\eq{BC-disc}.

\medskip{\it Step 5: Elimination of the cut-off $\CUT$ and further a-priori
estimates}. Having now \eq{Euler-small-viscoelastodyn+0eps} time-continuous, 
we can use \cite[Lemma\,5.1]{Roub24TVSE} like already mentioned in Step~1.
Instead of the transport-and-evolution equation \eq{Euler-small-viscoelastodyn+0},
let us consider the parabolic equation
\begin{align}\label{Euler-small-viscoelastodyn+0-reg}
\pdt{\varrho_\EEps}={\rm div}(\DELTA|\nabla\varrho_\EEps|^{r-2}\nabla\varrho_\EEps
-\varrho_\EEps\vv_\EEps\big)\,.
\end{align}
Applying the arguments based on the tests of \eq{Euler-small-viscoelastodyn+0-reg}
by $\varrho_\EEps$, by $\pdt{}\varrho_\EEps$, and by
${\rm div}(|\nabla\varrho_\EEps|^{r-2}\nabla\varrho_\EEps)$ as in the mentioned
\cite[Lemma\,5.1]{Roub24TVSE}, we again obtain the bound
$\|\varrho_\EEps\|_{L^\infty(W^{1,r}(\varOmega))}\le C_r$, which implies
$\varrho_\EEps<\rhoMAX$ with $\rhoMAX>N_rC_r$ and $N_r$ from in Step~1. It is important
that, analyzing the arguments in \cite{Roub24TVSE} where this bound is eventually
derived for $\DELTA=0$, we can see that this bound applies also for $\DELTA>0$ 
even independently of $\DELTA$. It is also important that the solution
to \eq{Euler-small-viscoelastodyn+0-reg} is unique not only for $\DELTA=0$ as
proved in \cite{Roub24TVSE} but also for each particular $\DELTA>0$.

Therefore, the solution $(\varrho_\EEps,\vv_\EEps,\Ee_\EEps)$
obtained in Step~4 which satisfies \eq{Euler-small-viscoelastodyn+0eps},
complying with $\varrho_\EEps|_{t=0}=\varrho_0<\rhoMAX$, satisfies
$0\le\varrho_\EEps(t,\cdot)<\rhoMAX$ at least for some sufficient small $t$'s
due to that $\varrho_\EEps\in C(I{\times}\barOmega)$, so that
$\CUT(\varrho_\EEps(t,\cdot))=1$. Thus it satisfies also
\eq{Euler-small-viscoelastodyn+0-reg}. By uniqueness for
\eq{Euler-small-viscoelastodyn+0-reg}, also the solution $\varrho_\EEps$ to
\eq{Euler-small-viscoelastodyn+0eps} is, for a given $\vv_\EEps$, determined
uniquely and, by a continuation argument and the choice of $\rhoMAX$,
it remains below $\rhoMAX$ on the entire time interval $I$. Therefore, we can
forgot $\CUT(\varrho_\EEps)$ both in \eq{Euler-small-viscoelastodyn+0eps} and in
\eq{Euler-small-viscoelastodyn+1eps}.

In addition to the estimates \eq{Euler-small-est}
with $C$ independent of $\EPS>0$ and of $\DELTA>0$ and thus, except the latter
estimates in \eq{Euler-small-est1} and \eq{Euler-small-est3},
inherited for $(\varrho_\EEps,\pp_\EEps,\Ee_\EEps)$,
we can now prove further estimates uniform in $\EPS>0$ and $\DELTA>0$.
The other estimates \eq{Euler-small-est+} and \eq{Euler-small-est++}
should now be replicated but benefiting from having the time-continuous problem
\eq{Euler-small-viscoelastodyn+0eps} for which the Gronwall inequality can
be used in a standard way.

Coming back to the test of \eq{Euler-small-viscoelastodyn+0-reg} by
${\rm div}(|\nabla\varrho_\EEps|^{r-2}\nabla\varrho_\EEps)$, in contrast to
\eq{Euler-small-est10++calcul}, we can handle the term arising from
${\rm div}\,\pp_\EEps$ by the Green formula and the Gronwall inequality, cf.\ the
calculations in \cite[Lemma\,5.1]{Roub24TVSE} for details. This leads to an
improvement of the estimate \eq{Euler-small-est10++} by the specific asymptotic
with respect to $\DELTA>0$ as 
\begin{subequations}\label{Euler-small-est-rho-E-}
\begin{align}\label{Euler-small-est-rho}
\|\varrho_\EEps\|_{L^\infty(I;W^{1,r}(\varOmega))}^{}\le C\ \ \ \text{ and }\ \ \
\|{\rm div}(|\nabla\varrho_\EEps|^{r-2}\nabla\varrho_\EEps)\|_{L^2(I\times\varOmega))}^{}
\le\frac C{\sqrt{\DELTA}}\,.
\intertext{Similarly, we can improve also the estimates (\ref{Euler-small-est}c,d) for
$\Ee$ when testing \eq{Euler-small-viscoelastodyn+2eps} by $\Delta\Ee_\EEps$.
Recalling the
qualification of the initial condition $\Ee_0\in H^1(\varOmega;\Rsym)$, we thus obtain}
\label{Euler-small-est-E}
\|\Ee_\EEps\|_{L^\infty(I;H^1(\varOmega;\R^{3\times3}))}^{}\le C\ \ \ \text{ and }\ \ \
\|\Delta\Ee_\EEps\|_{L^2(I\times\varOmega;\R^{3\times3}))}^{}\le\frac C{\sqrt{\DELTA}}\,.
\end{align}\end{subequations}
For this test, in addition, we have needed to estimate the term arising from
$R\,{\rm dev}\,\TT_\EEps$ as
\begin{align}\nonumber
&\int_\varOmega R\nabla({\rm dev}\TT_\EEps)\Vdots\nabla\Ee_\EEps\,\d\xx
=
\int_\varOmega  R\nabla\Big(\varphi'(\Ee_\EEps)
{-}\frac13({\rm tr}\varphi'(\Ee_\EEps)\bbI)\Big)\Vdots\nabla\Ee_\EEps\,\d\xx
\\[-.4em]
&\qquad\quad=\int_\varOmega R\Big(\varphi''(\Ee_\EEps)\nabla\Ee_\EEps
{-}\frac13({\rm tr}\varphi''(\Ee_\EEps)\nabla\Ee_\EEps)\bbI)\Big)\Vdots\nabla\Ee_\EEps\,\d\xx\le C\|\nabla\Ee_\EEps\|_{L^2(\varOmega;\R^{3\times3\times3})}^2
\label{Euler-small-est-E-}\end{align}
which is then to be treated by the Gronwall inequality. Here again we have formally
used the 2nd-order derivative $\varphi''$, but the estimate \eq{Euler-small-est-E-}
holds by smoothening argument for the at most linear growth of $\varphi'$ as assumed in
\eq{Euler-small-ass-phi}.

From \eq{Euler-small-est-rho}, we can further read a bound for the ``compensating
force'' $|\nabla\varrho_\EEps|^{r-2}(\nabla\vv_\EEps)\nabla\varrho_\EEps$ in
$L^{p}(I;L^{r'}(\varOmega;\R^3))$, so that 
\begin{subequations}\label{Euler-small-est-rho-E}
\begin{align}
&\big\|\DELTA|\nabla\varrho_\EEps|^{r-2}(\nabla\vv_\EEps)\nabla\varrho_\EEps\big\|_{L^{p}(I;L^{r'}(\varOmega;\R^3))}^{}\le C\DELTA\,.
\intertext{Moreover, the latter estimates in \eq{Euler-small-est-rho} and
\eq{Euler-small-est-E} imply in particular that}
&\big\|{\rm div}(\DELTA|\nabla\varrho_\EEps|^{r-2}\nabla\varrho_\EEps)\big\|_{L^2(I\times\varOmega))}^{}\le
C\sqrt{\DELTA}\ \ \ \text{ and }\ \ \
\|\DELTA\Delta\Ee_\EEps\|_{L^2(I\times\varOmega;\R^{3\times3}))}^{}\le C\sqrt{\DELTA}\,.
\end{align}\end{subequations}
From these estimates, by comparison from \eq{Euler-small-viscoelastodyn+0eps} and
from \eq{Euler-small-viscoelastodyn+2eps},
we can read also the estimates
\begin{subequations}\label{Euler-small-est+++}
\begin{align}
&\Big\|\pdt{\varrho_\EEps}\Big\|_{L^2(I\times\varOmega)}^{}\le C\ \ \ \text{ and }\ \ \ 
\Big\|\pdt{\Ee_\EEps}\Big\|_{L^2(I\times\varOmega;\R^{3\times3}))}^{}\le C\,;
\label{Euler-small-est-time-derivatives}
\intertext{here we have used also \eq{Euler-small-est-rho-E}.
Eventually, by comparison from \eq{Euler-small-viscoelastodyn+1eps}, we have also}
&\bigg\|\pdt{}\big(\varrho_\EEps\vv_\EEps\big)
\bigg\|_{L^{p'}(I;W^{2,p}(\varOmega;\R^3)^*)}^{}\le C\,.
\label{Euler-small-est-time-derivatives+}
\end{align}\end{subequations}

\medskip{\it Step 6: Limit passage for $\DELTA\to0$}.
By the Banach selection principle, we obtain a subsequence
converging weakly* with respect to the topologies indicated in
\eq{Euler-small-est}, \eq{Euler-small-est-rho-E-}, \eq{Euler-small-est-rho-E},
and \eq{Euler-small-est+++} to some limit
$(\varrho_\EPS,\vv_\EPS,\pp_\EPS,\Ee_\EPS)$. Specifically, 
\begin{subequations}\label{Euler-small-converge+}
\begin{align}
&\!\!\varrho_\EEps\to\varrho_\EPS&&\text{weakly* in $\ L^\infty(I;W^{1,r}(\varOmega))
\,\cap\, W^{1,p}(L^r(\varOmega))$\,,}
\\&\label{Euler-small-converge-bar-v+}
\!\!\vv_\EEps\to\vv_\EPS&&\text{weakly\ \;in $\ L^p(I;W^{2,p}(\varOmega;\R^3))$\,,}
\\&\label{Euler-small-converge-p+}
\!\!\pp_\EEps\to\pp_\EPS=\varrho_\EPS\vv_\EPS\!\!\!\!\!\!
&&\text{weakly\ \;in
$\ L^p(I;W^{1,r}(\varOmega;\R^3))\,\cap\,W^{1,p'}(I;W^{2,p}(\varOmega;\R^3)^*)$\,,}
\\
&\!\!\Ee_\EEps\to\Ee_\EPS\!\!\!&&\text{weakly* in $\ L^\infty(I;H^1(\varOmega;\Rsym))\,\cap\,H^1(I;L^2(\varOmega;\Rsym))$\,.}\!\!
\intertext{By the Arzel\`a-Ascoli and the Aubin-Lions theorems using
\eq{Euler-small-est+++}, we can further prove}
&\!\!\varrho_\EEps\to\varrho_\EPS&&\text{strongly in $\ C(I{\times}\barOmega)\,$,}
\\&\!\!\pp_\EEps\to\pp_\EPS
&&\text{strongly in $\ L^p(I;C(\barOmega;\R^3))\,$, and}
\\&\label{Euler-small-converge+E}
\!\!\Ee_\EEps\to\Ee_\EPS\!\!\!&&\text{strongly in 
$\ L^{a}(I;L^{2}(\varOmega;\Rsym))$ for $1\le a<\infty$ arbitrary}\,.
\intertext{Exploiting again the uniform monotonicity of the quasilinear operator
$\vv\mapsto{\rm div}({\rm div}(\HYPER|\Nabla^2\vv|^{p-2}\Nabla^2\vv)
-\bbD\strain(\vv))+\EPS|\vv|^{p-2}\vv$ and the convexity of the kinetic
energy as we did in \eq{Euler-small-strong-hyper-calc}--\eq{Euler-small-strong-hyper+++}
with a slight difference that the calculus in the inequality 
\eq{Euler-small-strong-hyper-calc} now turns simply to an equality, we can also
prove}
&\!\!\strain(\vv_\EEps)\to\strain(\vv_\EPS)\hspace*{-.5em}&&\text{strongly in
$\ L^2(I{\times}\varOmega;\Rsym)\ $ and}&&
\\&
\!\!\Nabla^2\vv_\EEps\to\Nabla^2\vv_\EPS\hspace*{-.5em}&&\text{strongly in
$\ L^p(I{\times}\varOmega;\R^{3\times3\times3})\,$}.&&
\end{align}\end{subequations}
Now, we can easily pass to the limit in the semilinear equations
\eq{Euler-small-viscoelastodyn+0eps} and \eq{Euler-small-viscoelastodyn+2eps}.

Furthermore, we can pass to the limit in the regularized
momentum equation \eq{Euler-small-viscoelastodyn+1eps}.
It is important that all three $\DELTA$-regularizing terms vanish for $\DELTA\to0$
due to the estimates \eq{Euler-small-est-rho-E}.

Altogether, we showed that $(\varrho_\EPS,\vv_\EPS,\Ee_\EPS)$ is a weak solution
to the original system \eq{Euler-small-viscoelastodyn+} with the boundary
conditions \eq{Euler-small-BC-hyper} but still with the momentum equation
\eq{Euler-small-viscoelastodyn+1} regularized by the ``dissipative force''
$\EPS|\vv|^{p-2}\vv$.

\medskip{\it Step 7: Further a-priori estimates and limit passage for $\EPS\to0$}.
Having now the mass-continuity equation nonregularized, i.e.\
\eq{Euler-small-viscoelastodyn+0}, we have a similar equation also for the
so-called {\it sparsity} $1/\rho$ at disposal, namely
\begin{align}\label{cont-eq-inverse}
\DT{\overline{\!\!\bigg(\frac1\varrho\bigg)\!\!}}\
=({\rm div}\,\vv)\frac{1}{\varrho}\,.
\end{align}
Having the initial density $\varrho_0\in W^{1,r}(\varOmega)$ qualified as
$\min_{\barOmega}^{}\varrho_0>0$, we have the initial condition
$1/\varrho_0\le1/\!\min_{\varOmega}\varrho_0$ with
$\nabla(1/\varrho_0)=-(\nabla\varrho_0)/\varrho_0^2\in W^{1,r}(\varOmega)$.
Using again \cite[Lemma\,5.1]{Roub24TVSE}, we obtain the a-priori bound
$\|1/\varrho_\EPS\|_{L^\infty(I;W^{1,r}(\varOmega))}\le C$. Using the embedding
$W^{1,r}(\varOmega)\subset L^\infty(\varOmega)$ relying on $r>3$, this shows
that $\min_{I\times\barOmega}^{}\varrho_\EPS>0$. We thus see that vacuum is well
ruled out due to the initial conditions and the hyperviscosity, taking into
account the fixed final time horizon $T$.

From the kinetic energy $\frac12\varrho_\EPS|\vv_\EPS|^2$, we then obtain also
the bound for $\vv_\EPS$ in $L^\infty(I;L^2(\varOmega;\R^3))$. Specifically,
from the first estimate in \eq{Euler-small-est1}
inherited for $\pp_\EPS/\!\sqrt{\varrho_\EPS}=\sqrt{\varrho_\EPS}\vv_\EPS$, 
we then have also
\begin{align}
&\|\vv_\EPS\|_{L^\infty(I;L^2(\varOmega;\R^3))}^{}\le
\|\sqrt{\varrho_\EPS}\vv_\EPS\|_{L^\infty(I;L^2(\varOmega;\R^3))}^{}
\bigg\|\frac1{\!\sqrt{\varrho_\EPS}}\bigg\|_{L^\infty(I\times\varOmega)}^{}\le C\,.
\label{basic-est-of-v}
\end{align}
As this estimate is uniform in $\EPS>0$, we can now re-make the previous
estimates which have been supported by only by the nonuniform 
$L^p(I{\times}\varOmega)$-estimate inherited for $\vv_\EEps$ from the latter estimate in
\eq{Euler-small-est1}. Actually, we have now even stronger bound
for $\vv_\EPS$ since
$L^\infty(I;L^2(\varOmega))\cap L^p(I;W^{2,p}(\varOmega))\subset
L^\infty(I;L^2(\varOmega))\cap L_{\rm w*}^p(I;L^\infty(\varOmega))\subset
L^{p+2}(I{\times}\varOmega)$. 

Therefore, we have at our disposal these other estimates from Steps~4-5 which
are now uniform for both $\EPS$ and $\DELTA$, i.e.\ the former estimates in
\eq{Euler-small-est-rho} and \eq{Euler-small-est-E} and also
\eq{Euler-small-est+++} at our disposal now for $(\varrho_\EPS,\vv_\EPS,\Ee_\EPS)$.
This concerns also \eq{Euler-small-est} except the later
estimate in \eq{Euler-small-est3} which is now ensured
by the assumed regularity of the initial conditions independently
of $\DELTA$ and, of course, on $\EPS$ too. Similarly, also
\eq{Euler-small-est10+} and \eq{Euler-small-est6+} are
now independent of $\DELTA$ and on $\EPS$ and thus inherited 
for $\pp_\EPS$ independently of $\EPS$. Having now the regularizing
$\DELTA$-Laplacians dismissed, we can improve \eq{Euler-small-est-time-derivatives} for
\begin{align}
\Big\|\pdt{\varrho_\EPS}\Big\|_{L^p(I;L^r(\varOmega))}^{}\le C\ \ \ \ \text{ and }\ \ \ \ 
\Big\|\pdt{\Ee_\EPS}\Big\|_{L^p(I;L^2(\varOmega;\R^{3\times3}))}^{}\le C\,.
\label{Euler-small-est-time-derivatives+++}
\end{align}

We now replicate the arguments in Step~6. Instead of
\eq{Euler-small-converge-bar-v+}, we have now a bit
finer convergence, namely $\vv_\EPS\to\vv$ weakly*
in $L^\infty(I;L^2(\varOmega;\R^3))\,\cap\,L^p(I;W^{2,p}(\varOmega;\R^3))$.
In addition, we need to pass to the limit of the regularizing dissipative force
$\EPS|\vv_\EPS|^{p-2}\vv_\EPS$, which is however simple because
we have the mentioned bound of $\vv_\EPS$ in $L^{p+2}(I{\times}\varOmega;\R^3)$
so that the norm of $\EPS|\vv_\EPS|^{p-2}\vv_\EPS$ in
$L^{(2+p)/(p-1)}(I{\times}\varOmega;\R^3)$ is $\mathscr{O}(\EPS)$ and
it surely converges to zero for $\EPS\to0$.

Altogether, we obtain a weak solution $(\varrho,\vv,\Ee)$ to the original system
\eq{Euler-small-viscoelastodyn+} with the boundary conditions \eq{Euler-small-BC-hyper}.

\medskip{\it Step 8: the energy-dissipation balance}.
It is now important that the tests and then all the subsequent calculations
leading to the energy balance \eq{Euler-small-energy-balance-stress}
integrated over a current time interval $[0,t]$ are analytically legitimate.

In particular, for the calculus \eq{Euler-small-divT.v++}, we rely on that 
$\TT=\varphi'(\Ee)+\varphi(\Ee)\bbI\in L^\infty(I;L^3(\varOmega;\Rsym))$
is surely in duality with $\strain(\vv)\in L^p(L^\infty(\varOmega;\Rsym))$
and $\varphi'(\Ee)\in L^\infty(I;L^6(\varOmega;\Rsym))$
is in duality both with $\ZJ\Ee\in L^p(I;L^2(\varOmega;\Rsym))$
and with ${\rm dev}\,\TT\in L^\infty(I;L^6(\varOmega;\Rsym))$, too.
Thus the calculus \eq{Euler-small-divT.v++} is indeed legitimate when
integrated over time.
Moreover, $\pdt{}(\varrho\vv)\in L^{p'}(I;W^{2,p}(\varOmega;\R^3)^*)$ due to
\eq{Euler-small-converge-p+} and
${\rm div}(\varrho\vv{\otimes}\vv)\in L^{1+p/2}(I{\times}\varOmega;\Rsym)$ as well
as ${\rm div}\DD\in L^{p'}(I;W^{2,p}(\varOmega;\R^3)^*)$ are in duality with
$\vv\in L^p(I;W^{2,p}(\varOmega;\R^3))$, as used when testing the momentum
equation by $\vv$, in particular in
\eq{rate-of-kinetic}--\eq{calculus-convective} and in \eq{Euler-test-momentum++}.
The calculus \eq{rate-of-kinetic} also relies on that both $\pdt{}\varrho$ and
${\rm div}(\varrho\vv)=\vv{\cdot}\nabla\varrho+\varrho{\rm div}\vv$ live in
$L^p(I;L^r(\varOmega))$ and thus are surely in duality with $|\vv|^2\in
L^{p/2}(I;L^\infty(\varOmega))$ for $p\ge3$.

\begin{remark}[{\sl More general stored energies $\varphi$'s}]\label{rem-general-phi}\upshape
The growth assumption \eq{Euler-small-ass-phi} can be weakened in order to
allow $\varphi$ with a super-quadratic growth when using the
Gagliardo-Nirenberg interpolation and the estimates (\ref{Euler-small-est}a--d).
Even bigger growth of $\varphi'$ and $\varphi$ can be allowed if one uses some
approximation of $\varphi$ in Steps~3-4 (e.g.\ by the Yosida approximation
like in the following Section~\ref{sec-polycon} parameterized here by $\DELTA$)
and then involving it in the limit passage in Step~7 while employed the estimate
of $\Ee_\EPS$ in $L^\infty(I;H^1(\varOmega;\Rsym))\subset L^\infty(I;L^6(\varOmega;\Rsym))$
obtained in Step~5, cf.\ \eq{Euler-small-est-E}.
\end{remark}

\begin{remark}[{\sl Boundary forcing}]\label{rem-f-ne-0}\upshape
One could think about a tangential forcing on the boundary, i.e.\ the
second boundary condition in \eq{Euler-small-BC-hyper} nonhomogeneous with
some ${\bm f}\ne\bm0$ in $L^2(I{\times}\varGamma;\R^3)$ with ${\bm f}\Cdot\nn=0$.
The energy balance \eq{Euler-small-energy-balance-stress} would then be augmented by 
a boundary power $\int_\varGamma{\bm f}\Cdot\vv\,\d S$. The analysis 
could be done by considering the Robin-type (Newton) boundary condition with
${\bm f}_\tau^k-\EPS\vv_\etau^k$ in the right-hand side of the middle condition
in \eq{BC-disc-a}, giving rise to a regularizing boundary dissipation term
$\EPS|\vv_\etau^k|^2$. After successive passage to the limit with $\tau\to0$ and
$\DELTA\to0$, one could estimate the boundary term
$\int_\varGamma{\bm f}\Cdot\vv_\EPS\,\d S$ 
by exploiting the bulk terms on the left-hand side, not relying on
that boundary dissipation term like we did in Step~6 for the
$\EPS$-regularizing force.
\end{remark}

\section{Finitely-strained visco-elastodynamics}\label{sec-finite-strain}

Having scrutinized most of technical difficulties of the simpler linearized
model in the previous Section~\ref{sec-linearized}, we will now present
the usage of the time-discretization method on the fully nonlinear Eulerian
visco-elastodynamics where additional technicalities naturally arise. Instead
of the symmetric small-strain
tensor $\Ee$, it uses the Eulerian deformation gradient $\FF_{\rm tot}$
(or the elastic distortion $\FF$), although the deformation $\yy$ itself
does not need to be explicitly involved in the model when formulated in rates;
actually $\yy$ even does not need to exist but, if it exists, then
$\FF_{\rm tot}=\nabla\yy$; cf.\ \cite[Remark~7]{Roub22VELS}. Without going
into details of this standard approach, we refer to many renowned monographs, e.g.\
to \cite{GuFrAn10MTC,Mart19PCM,Rubi21CMEF}. Specifically, the kinematic equation
$\ZJ\Ee=\strain(\vv)$ now turns into $\DT\FF_{\rm tot}=(\nabla\vv)\FF_{\rm tot}$
and $\ZJ\Ee=\strain(\vv)-R\,{\rm dev}\,\TT$ in \eq{Euler-small-flow-rule} turns into
$\DT\FF=(\nabla\vv)\FF-R\FF\MM$ with the so-called Mandel stress
$\MM={\rm dev}(\FF^\top\!\varphi'(\FF))$, i.e.\
$\DT\FF=(\nabla\vv)\FF-\FF\Lp$ with the trace-free inelastic distortion rate
$\Lp$ governed by the inelastic flow-rule $\Gm\Lp=-\MM$.

Here, however, the application of this method is complicated by the fact that, in
contrast to the linearized convective model in Sect.~\ref{sec-linearized}, the stored
energy is necessarily ``heavily'' nonconvex. Anyhow, in some situations this
nonconvexity can be overcome. In particular, this is the case when the actual
stored energy is polyconvex and independent of cofactors; later we will generalize
this to general polyconvex energies in Remark~\ref{rem-polyconvex-general} below
or even more general energies in Remark~\ref{rem-behind-polyconvexity}.

Another particular complication arises for $R>0$, i.e.\ when actually the Jeffreys
but not only the mere Kelvin-Voigt rheology is considered in the deviatoric part
of the model. Rigorous analysis requires some gradient in the inelastic distortion,
either in the conservative part \cite[Sect.9.4.1]{KruRou19MMCM}
or in the dissipative part, the latter option being more suitable because it avoids
spurious hardening under large slip, as pointed out in \cite{DaRoSt21NHFV},
and is easily implemented in the Eulerian frame \cite{Roub22QHLS}. It leads
to the inelastic flow rule for the inelastic distortion rate $\Lp$ enhanced
by the quasi-linear operator ${\rm div}(\nu|\nabla\Lp|^{q-2}\nabla\Lp)$ on
$\varOmega$ with the boundary condition $(\nn\Cdot\nabla)\Lp=\bm0$ on $\varGamma$,
i.e.\ the quasi-static elliptic problem $\Gm\Lp=
{\rm div}(\nu|\nabla\Lp|^{q-2}\nabla\Lp)-\MM$. Unlike Sect.~\ref{sec-linearized-system},
we will not eliminate the variable $\Lp$ from the system. For the influence of this
additional gradient on the dispersion (analysed on a linear 1D case)
cf.\ \cite[Sect\,3.4 and 3.6]{Roub24SGTL}.

\subsection{The system and its energetics}

Formulated in terms of the linear momentum $\pp$ as we did in
\eq{Euler-small-viscoelastodyn+trans}, the considered system reads as
\begin{subequations}\label{ED-Euler-large}\begin{align}
\label{ED-Euler-large0}&\pdt\varrho=-{\rm div}\,\pp\,,\\
\nonumber
     &\pdt\pp={\rm div}\big(\TT+\DD-\pp{\otimes}\vv\big)
    +\varrho\GRAVITY\ \ \ \ \ \ \text{ with }\ \pp=\varrho\vv\,,\ \ \
     \TT=\varphi'(\FF)\FF^\top\!+\varphi(\FF)\bbI\,,
 \\[-.2em]
    &\hspace*{3.5em}
 \text{ and }\ \ \DD=\bbD\strain(\vv)-{\rm div}\,\mathfrak{H}
 \ \ \text{ with }\ \ 
      \mathfrak{H}=\HYPER|\nabla^2\vv|^{p-2}\nabla^2\vv\,,
\label{ED-Euler-large1}
\\[-.3em]\label{ED-Euler-large2}
&\pdt\FF=(\nabla\vv)\FF-(\vv\Cdot\nabla)\FF-\FF\Lp\,,
\\[-.1em]\label{ED-Euler-large3}
&\Gm\Lp={\rm div}\big(\nu|\nabla\Lp|^{q-2}\nabla\Lp\big)-\MM\ \ \text{ with }\
\MM={\rm dev}(\FF^\top\!\varphi'(\FF))
\end{align}\end{subequations}
with the boundary conditions \eq{Euler-small-BC-hyper} enhanced by the boundary
condition $(\nn\Cdot\nabla)\Lp=\bm0$ on $\varGamma$. It should be emphasized that
$\varphi$ is meant in the actual frame in the physical units J/m$^3$. In contrast to
Sect.~\ref{sec-linearized}, here there are two other alternatives  that are even more
often considered in literature working with the stored energy $\upvarphi$ in a
referential frame either in the physical units J/m$^3$
(then $\varphi=\upvarphi/\!\det\FF$ and $\TT=\upvarphi'(\FF)\FF^\top\!/\!\det\FF$) or
in the physical units J/kg (then $\upvarphi=\varrho\varphi$ and
$\TT=\varrho\upvarphi'(\FF)\FF^\top$), cf.\
Remarks~\ref{rem-barotropic-fluids} and \ref{rem-referential} below. For $R=0$, the
model works with the Kelvin-Voigt rheology while $R>0$ describes the Jeffreys
rheology in the deviatoric part. In the former case, simply $\FF=\FF_{\rm tot}$.
In the latter case, this model is based on the
{\it multiplicative Kr\"oner-Lee-Liu decomposition} of
$\FF_{\rm tot}$ into the product of the elastic distortion $\FF$ and the
inelastic distrortion, and subsequent elimination of both the total deformation
gradient and the mentioned inelastic distortion.

The energetics behind the system \eq{ED-Euler-large} uses again
\eq{calculus-convective} and \eq{Euler-test-momentum++} while \eq{Euler-small-divT.v++}
is modified as
\begin{align}\nonumber
\!\!\!\!\int_\varGamma\!&\vv\Cdot\TT\nn\,\d S-\!\int_\varOmega\!{\rm div}\TT\Cdot\vv\,\d\xx
=\!\int_\varOmega\!\big(\varphi'(\FF)\FF^\top\big)\Colon\Nabla\vv
+\varphi(\FF){\rm div}\,\vv\,\d\xx
\\[-.4em]&\nonumber
=\int_\varOmega\!\varphi'(\FF)\Colon\big((\Nabla\vv)\FF\big)+\varphi(\FF){\rm div}\,\vv
=\!\int_\varOmega\!\varphi'(\FF)\Colon\big(\DT\FF\!{+}
\FF\Lp\big)+\varphi(\FF){\rm div}\,\vv\,\d\xx
\\&\nonumber
\int_\varOmega\!\varphi'(\FF)\Colon\DT\FF+\MM\Colon\Lp+\varphi(\FF){\rm div}\,\vv\,\d\xx
\\&\nonumber
=\int_\varOmega\!
\varphi'(\FF)\Colon\pdt{\FF}+\MM\Colon\Lp
+\!\!\!\!\lineunder{\varphi'(\FF)\Colon(\vv\Cdot\Nabla)\FF\!
+\varphi(\FF){\rm div}\,\vv}
{$\hspace{5em}=\int_\varGamma\varphi(\FF)(\vv\Cdot\nn)\,\d S=0\hspace{-5em}$}\!\!\!\d\xx
\\[-.9em]&
=\frac{\d}{\d t}\!\int_\varOmega\!\varphi(\FF)\,\d\xx+\!\int_\varOmega\!
\Gm|\Lp|^2\!+\nu|\nabla\Lp|^q\,\d\xx\,.
\label{Euler-hypoplast-test-momentum}\end{align}
Including still the calculations concerning the Stokes (hyper)viscosity
\eq{Euler-test-momentum++}, the energy-dissipation balance
\eq{Euler-small-energy-balance-stress} now turns into
\begin{align}
\frac{\d}{\d t}\int_\varOmega\!\!\!\!
\linesunder{\frac\varrho2|\vv|^2}{kinetic}{energy}
\!\!\!\!\!+\!\!\!\!\!\linesunder{\varphi(\FF)}{stored}{energy}\!\!\!\!
\d\xx
+\!\int_\varOmega\hspace{-.7em}
\linesunder{\bbD\strain(\vv)\Colon\strain(\vv){+}\HYPER|\nabla^2\vv|^p}{disipation rate due to}{the Stokes-type viscosity}\hspace{-.7em}
+\hspace{-.7em}\linesunder{\Gm\big|\Lp\big|^2\!+\nu\big|\nabla\Lp\big|^q
}{disipation rate due to}{the Maxwellian viscosity}\hspace{-.7em}
\d\xx
=\int_\varOmega\hspace{-1.4em}\linesunder{\varrho\GRAVITY\Cdot\vv}{power of}{gravity field}\hspace{-1.4em}
  \,\d\xx\,.
\label{Euler-large-energy-balance-stress}\end{align}

Let us still mention the equation for the Jacobian $J:=\det\FF$ arising
by application ``det'' to the equation \eq{ED-Euler-large2}. Exploiting
that the Mandel stress $\MM$ and thus also the inelastic distortion rate $\Lp$
are trace free and using the calculus
$\det'A={\rm Cof}\,A$ and the Cramer rule $A^{-1}={\rm Cof}\,A^\top\!/\!\det A$,
we obtain
\begin{align}\nonumber
\!\DT J&={\rm Cof}\FF\Colon\DT\FF=J\FF^{-\top}\!\!\Colon\DT\FF
=J\bbI\Colon\DT\FF\FF^{-1}\!=J\bbI\Colon\big(\Nabla\vv{-}\FF\Lp\FF^{-1}\big)
\\&=J{\rm div}\,\vv-J\FF^{-\top}\!\Colon\FF\Lp
=J{\rm div}\,\vv
-J\hspace{-1.7em}\lineunder{\bbI\Colon\Lp}{$\ \ \ \
={\rm tr}\Lp=0$}\hspace{-2em}=({\rm div}\,\vv)J.
\label{DT-det-extended}\end{align}
A similar equation holds also for the reciprocal Jacobian $A:=1/J$, i.e.\ the
determinant of the so-called distortion $\AA:=\FF^{-1}$. Indeed,
$\DT J/J={\rm div}\,v$, which is just \eq{DT-det-extended}, can
be written as $\DT{\overline{{\rm ln}\,J}}={\rm div}\,v$ and then, from the standard
calculus ${\rm ln}\,A=-{\rm ln}\,J$, we obtain
$\DT{\overline{{\rm ln}\,A}}=-{\rm div}\,v$, i.e.
\begin{align}\nonumber\\[-2em]
\DT{\overline{\!\!\Big(\frac1J\Big)\!\!}}=-({\rm div}\,\vv)\frac1J\,.
\label{DT-det-inverse}\end{align}

\subsection{Special polyconvex actual energies and time discretization}\label{sec-polycon}

We want to illustrate the time-discretization method first on the simplest
case of the actual stored energy $\varphi:{\rm GL}_3^+\to\R$ in the form
$\varphi(\FF)=\breve\varphi(\FF,\det\FF)$ with
$\breve\varphi=\breve\varphi(\FF,J)$ as a convex function
$\R^{3\times3}\times\R^+\to\R\cup\{+\infty\}$; here
${\rm GL}_3^+\subset\R^{3\times3}$ denotes the so-called identity component of the
general linear group, i.e.\ the set of matrices with positive determinant.
Later, in Section~\ref{sec-rem-exa}, we will generalize this ansatz
to allow for broader applications. Again, we consider the initial-value problem
for the system \eq{ED-Euler-large} with the initial conditions
\begin{align}\label{IC-large}
\varrho|_{t=0}^{}=\varrho_0^{}\,,\ \ \ \ \ \
\pp|_{t=0}^{}=\varrho_0^{}\vv_{0}^{}\,,\ \ \ \text{ and }\ \ {\FF}|_{t=0}^{}=\FF_0^{}
\,.\end{align}

\begin{definition}[Weak formulation of \eq{ED-Euler-large}]\label{def-ED-Ch5}
We call $(\varrho,\vv,\FF,\Lp)$ with
$\varrho\in L^\infty(I{\times}\varOmega)\cap W^{1,1}(I{\times}\varOmega)$,
$\vv\in L^p(I;W^{2,p}(\varOmega;\R^3))$,
$\FF\in W^{1,1}(I{\times}\varOmega;\R^{3\times3})$, and
$\Lp\in L^2(I;W^{1,q}(\varOmega;\R_{\rm dev}^{3\times3}))$ a weak solution to the system
\eq{ED-Euler-large} with the boundary conditions \eq{Euler-small-BC-hyper}
enhanced by the boundary condition $(\nn\Cdot\nabla)\Lp=\bm0$ on $\varGamma$ and
the initial conditions \eq{IC-large} if 
\begin{align}
&\nonumber
\int_0^T\!\!\!\int_\varOmega\bigg(\Big(\varphi'(\FF)\FF^\top\!+
\bbD\strain(\vv)-\varrho\vv{\otimes}\vv\Big){:}\strain(\widetilde\vv)
+\varphi(\FF){\rm div}\,\widetilde\vv
-\varrho\vv{\cdot}\pdt{\widetilde\vv}
\\[-.6em]&\hspace{7em}
+\HYPER|\nabla^2\vv|^{p-2}\nabla^2\vv\Vdots\nabla^2\widetilde\vv\bigg)\,\d\xx\d t
=\!\int_0^T\!\!\!\int_\varOmega\varrho\GRAVITY{\cdot}\widetilde\vv\,\d\xx\d t
+\!\int_\varOmega\!\varrho_0^{}\vv_0^{}{\cdot}\widetilde\vv(0)\,\d\xx
\label{def-ED-Ch5-momentum}\end{align}
holds for any $\widetilde\vv$ smooth with $\widetilde\vv{\cdot}\nn={\bm0}$ on
$I{\times}\varGamma$ and $\widetilde\vv(T)=0$, and if \eq{ED-Euler-large0},
\eq{ED-Euler-large2}, and \eq{ED-Euler-large3} hold a.e.\ on $I{\times}\varOmega$
together with the respective initial conditions for $\varrho$ and $\FF$ in
\eq{IC-large}.
\end{definition}

We impose the assumptions:
\begin{subequations}\label{Euler-large-ass}
\begin{align}\label{Euler-large-ass-phi}
&\breve\varphi\in C^1(\R^{3\times3}{\times}\R^+)\ \text{ convex},\ \ \ \ \ \ 
\lim_{\det F\to0}\breve\varphi(F,\det F)=+\infty
\,,
\ \ \ \ 
\\[-.1em]&\label{Euler-large-ass-IC}
\varrho_0\in W^{1,r}(\varOmega)\,,\ {\rm min}_{\barOmega}^{}\varrho_0>0\,,\ \
\vv_0\in L^2(\varOmega;\R^3)\,,\ \ \FF_0\in W^{1,r}(\varOmega;\R^{3\times3})\,,\ \
\ {\rm min}_{\barOmega}^{}\det\FF_0>0\,.\end{align}\end{subequations}
The blow-up condition in \eq{Euler-large-ass-phi} is the physically expected
property that ensures $\det\FF>0$ a.e., which is the {\it local non-interpenetration}.
The following assertion will be proved in Sect.\,\ref{sec-large-proof} below:

\begin{proposition}[Existence of weak solutions to \eq{ED-Euler-large}]
\label{prop-ED-Ch5-existence}
Let \eq{Euler-large-ass} with (\ref{Euler-small-ass}b,d) hold with $r>3$ and let
$\min(p,q)>3$. Then the initial-boundary-value
problem for the system \eq{ED-Euler-large} has at least one weak
solution $(\varrho,\vv,\FF,\Lp)$ in the sense of Definition~\ref{def-ED-Ch5} such that
also $\varrho\in L^\infty(I;W^{1,r}(\varOmega))\,\cap\,C(I{\times}\barOmega)$ with
$\min_{I{\times}\barOmega}\varrho>0$, $\vv\in L^\infty (I;L^2(\varOmega;\R^3))$,
$\FF\in L^\infty(I;W^{1,r}(\varOmega;\R^{3\times3}))
\,\cap\,C(I{\times}\barOmega;\R^{3\times3})$ with $\min_{I\times\barOmega}\det\FF>0$,
and ${\rm div}(|\nabla\Lp|^{q-2}\nabla\Lp)\in L^2(I{\times}\varOmega;\R_{\rm dev}^{3\times3})$
and the energy-dissipation balance \eq{Euler-large-energy-balance-stress} integrated
over the time interval $I$.
\end{proposition}

For the actual stored energy in the form $\varphi(\FF)=\breve\varphi(\FF,\det\FF)$
with $\breve\varphi=\breve\varphi(\FF,J)$ as a function
$\R^{3\times3}\times\R^+\to\R\cup\{+\infty\}$, we have the Cauchy and the Mandel
stresses in terms of $\breve\varphi$ as
\begin{subequations}\label{Cauchy-Mandel-stress-polyconvex}\begin{align}\nonumber
\TT&=\varphi'(\FF)\FF^\top\!+\varphi(\FF)\bbI
=\breve\varphi_\FF'(\FF,J)\FF^\top\!+\breve\varphi_J'(\FF,J)({\rm Cof}\FF)\FF^\top\!
+\breve\varphi(\FF,J)\bbI
\\&\hspace{7.7em}=\breve\varphi_\FF'(\FF,J)\FF^\top\!
+\Big(J\breve\varphi_J'(\FF,J){+}\breve\varphi(\FF,J)\Big)\bbI\ \ \ \ \text{ and}
\label{Cauchy-stress-polyconvex}\\\nonumber
\MM&={\rm dev}\big(\FF^\top\!\varphi'(\FF)\big)=
{\rm dev}\Big(\FF^\top\!\breve\varphi_\FF'(\FF,J)+\FF^\top\!\breve\varphi_J'(\FF,J)({\rm Cof}\FF)\Big)
\\&\hspace{7.7em}={\rm dev}\Big(\FF^\top\!\breve\varphi_\FF'(\FF,J)+\breve\varphi_J'(\FF,J)\bbI\Big)
={\rm dev}\big(\FF^\top\!\breve\varphi_\FF'(\FF,J)\big)\,.
\label{Mandel-stress-polyconvex}\end{align}\end{subequations}
The system \eq{ED-Euler-large} with \eq{Cauchy-Mandel-stress-polyconvex},
involving  \eq{DT-det-extended}, then reads as a system for the quintuple
$(\varrho,\pp,\FF,J,\Lp)$:
\begin{subequations}\label{ED-Euler-large-polyconvex}\begin{align}
\label{ED-Euler-large0-polyconvex}&\pdt\varrho=-{\rm div}\,\pp\,,\\
\nonumber
     &\pdt\pp={\rm div}\big(\TT+\DD-\pp{\otimes}\vv\big)
     +\varrho\GRAVITY\ \ \ \ \ \text{ with }\ \pp=\varrho\vv\,,
 \\[-.2em]\nonumber
    &\hspace*{4em}\text{where }\ \ 
    \TT=\breve\varphi_\FF'(\FF,J)\FF^\top\!+
    \Big(J\breve\varphi_J'(\FF,J){+}\breve\varphi(\FF,J)\Big)\bbI
\\&\hspace*{4em}
\text{and }\ \ \ \ \DD=\bbD\strain(\vv) -{\rm div}\,\mathfrak{H}
 \ \ \text{ with }\ \ 
      \mathfrak{H}=\HYPER|\nabla^2\vv|^{p-2}\nabla^2\vv\,,
\label{ED-Euler-large1-polyconvex}
\\[-.3em]\label{ED-Euler-large2-polyconvex}
&\pdt\FF=(\nabla\vv)\FF-(\vv\Cdot\nabla)\FF-\FF\Lp\,,
\\[-.0em]\label{ED-Euler-large3-polyconvex}
&\pdt J=({\rm div}\,\vv)J-\vv\Cdot\nabla J\,,
\\[-.1em]\label{ED-Euler-large4-polyconvex}
&\Gm\Lp={\rm div}\big(\nu|\nabla\Lp|^{q-2}\nabla\Lp\big)-\MM
\ \ \text{ with }\
\MM={\rm dev}\big(\FF^\top\!\breve\varphi_\FF'(\FF,J)\big)\,
\end{align}\end{subequations}
with the boundary conditions \eq{Euler-small-BC-hyper} enhanced by the
boundary condition $(\nn\Cdot\nabla)\Lp=\bm0$ on $\varGamma$ and the initial
conditions
\begin{align}\label{IC-large+}
\varrho|_{t=0}^{}=\varrho_0\,,\ \ \ \ \ \ \pp|_{t=0}^{}=\varrho_0\vv_{0}\,,\  
\ \ \ \ \ {\FF}|_{t=0}^{}=\FF_0\,,  \ \ \text{ and }\ \ \ J|_{t=0}^{}=\det\FF_0\,.
\end{align}

To reveal the energetics behind \eq{ED-Euler-large-polyconvex}, we
use \eq{calculus-convective} and \eq{Euler-test-momentum++} as before while
the calculations \eq{Euler-hypoplast-test-momentum} are then to be modified as 
\begin{align}\nonumber
\int_\varGamma&(\TT\nn)\Cdot\vv\,\d S-\int_\varOmega{\rm div}\,\TT\Cdot\vv\,\d\xx
=\int_\varOmega\TT\Colon\nabla\vv\,\d\xx
\\[-.3em]&\nonumber
\stackrel{\eq{Cauchy-stress-polyconvex}}{=}\,\int_\varOmega\Big(\breve\varphi_\FF'(\FF,J)\FF^\top
+\big(J\breve\varphi_J'(\FF,J){+}\breve\varphi(\FF,J)\big)\bbI\Big)\Colon\nabla\vv\,\d\xx
\\[-1.2em]&\nonumber
\!\!\UUU{=}{\eq{DT-det-extended}}{\eq{ED-Euler-large2-polyconvex}}\!\!\!\!
\int_\varOmega\breve\varphi_\FF'(\FF,J)\Colon\big(\DT\FF{+}
\FF\Lp\big)+\breve\varphi_J'(\FF,J)\DT J+\breve\varphi(\FF,J){\rm div}\,\vv\,\d\xx
\\&\nonumber
\stackrel{\eq{Mandel-stress-polyconvex}}{=}
\int_\varOmega\breve\varphi_\FF'(\FF,J)\Colon\DT\FF+
{\rm dev}\big(\FF^\top\!\breve\varphi_\FF'(\FF,J)\big)\Colon\Lp
+\breve\varphi_J'(\FF,J)\DT J-\nabla\breve\varphi(\FF,J)\Cdot\vv\,\d\xx
\\[-1.2em]&\nonumber
\!\!\UUU{=}{\eq{ED-Euler-large4-polyconvex+}}{\eq{ED-Euler-large4-polyconvex}}\!\!\!\!
\int_\varOmega\breve\varphi_\FF'(\FF,J)\Colon\DT\FF+\MM\Colon\Lp
+\breve\varphi_J'(\FF,J)\DT J-\breve\varphi_\FF'(\FF,J)\Colon(\vv\Cdot\nabla)\FF
-\breve\varphi_J'(\FF,J)\vv\Cdot\nabla J\,\d\xx
\\&\nonumber
\!\stackrel{\eq{ED-Euler-large4-polyconvex}}{=}\!\!
\int_\varOmega\breve\varphi_\FF'(\FF,J)\Colon\pdt{\FF}+
\Gm|\Lp|^2\!+\nu|\nabla\Lp|^q\!+\breve\varphi_J'(\FF,J)\pdt{J}\,\d\xx
\\&
\ \ =\frac{\d}{\d t}\int_\varOmega\!\!\!\!\lineunder{\breve\varphi(\FF,J)}{$=\varphi(\FF)$}\!\!\!\!\d\xx
+\int_\varOmega\!\!
\Gm|\Lp|^2\!+\nu|\nabla\Lp|^q\,\d\xx\,,
\label{ED-Euler-large-stress-calculus}\end{align}
where we used the calculus like \eq{Euler-small-calc}, namely
\begin{align}\nonumber
0=\int_\varOmega\!{\rm div}\big(\breve\varphi(\FF,J)\vv\big)\,\d\xx
&=\int_\varOmega\!\nabla\breve\varphi(\FF,J)\Cdot\vv+\breve\varphi(\FF,J){\rm div}\,\vv\,\d\xx
\\&\label{ED-Euler-large4-polyconvex+}=
\int_\varOmega\breve\varphi_\FF'(\FF,J)\Colon(\vv\Cdot\nabla\FF)+
\breve\varphi_J'(\FF,J)\vv\Cdot\nabla J+\breve\varphi(\FF,J){\rm div}\,\vv\,\d\xx\,.
\end{align}

We can then devise the recursive regularized time-discrete scheme in the spirit
of \eq{Euler-small-viscoelastodyn+disc} when discretizing also \eq{DT-det-extended}.
For the convex $\breve\varphi$, the original actual stored energy $\varphi$
is {\it polyconvex}. Yet, to replicate the analysis which exploited strong convexity
of $\varphi$ in \eq{Euler-ED-basic-energy-balance-disc}, we would need a
``strong polyconvexity'' of $\varphi$ in the sense that $\breve\varphi$ is
jointly strongly convex, which however would not be much realistic 
in view of applications, cf.\ the examples in Section~\ref{sec-rem-exa}
below. Moreover, it is desirable to regularize the singularity of
$F\mapsto\breve\varphi(F,\det F)$ for $\det F\to0+$. Thus, we also regularize
$\breve\varphi$ by defining
\begin{align}\nonumber
\breve\varphi_\EEps(F,J):=&[\mathscr{Y}_\EPS\breve\varphi](F,J)
+\frac{\sqrt\DELTA}2|F|^2+\frac{\sqrt\DELTA}2J^2
\\[-.4em]&\ \ \text{ with }\ \ [\mathscr{Y}_\EPS\breve\varphi](F,J)
=\inf_{\wt F\in\R^{3\times3},\,\wt J>0}
\Big(\breve\varphi(\wt F,\wt J)+\frac1{2\EPS}|\wt F{-}F|^2
+\frac1{2\EPS}|\wt J{-}J|^2\Big)\,.
\label{varphi-regularized}\end{align}
The convex function $\mathscr{Y}_\EPS\breve\varphi$ is the so-called {\it Yosida
approximation} of $\breve\varphi$ and has an at-most quadratic growth.
In contrast to $\breve\varphi$, this approximation is defined (and finite)
on the whole $\R^{3\times3}\times\R$, not only on $\R^{3\times3}\times\R^+$.
For any fixed $\DELTA>0$ and $\EPS>0$, the function $\breve\varphi_\EEps$ is
continuously differentiable with a bounded derivative, strongly convex, and
coercive with a quadratic growth. The convexity of $\breve\varphi$ guarantees
the convergence $[\mathscr{Y}_\EPS\breve\varphi]\to\breve\varphi$ in the so-called
epi-graphical sense. Moreover, for replacing $\breve\varphi$ by $\breve\varphi_\EEps$
in \eq{Cauchy-Mandel-stress-polyconvex}, let us abbreviate the ``stress-mappings''
$\mathscr{T}_\EEps,\mathscr{M}_\EEps:\R^{3\times3}\times\R\to\R^{3\times3}$ as continuous
mappings with at most quadratic growth defined by
\begin{subequations}\label{Euler-larger-polyconvex-stress}
\begin{align}\label{Euler-larger-polyconvex-Cauchy}
&\mathscr{T}_\EEps(F,J):=
\big([\mathscr{Y}_\EPS\breve\varphi]_F'(F,J){+}\sqrt\DELTA F\big)F^\top\!
+\Big(J[\mathscr{Y}_\EPS\breve\varphi]_J'(F,J){+}\sqrt\DELTA J^2{+}
\breve\varphi_\EEps(F,J)\Big)\bbI\ \ \ \text {and}
\\&\label{Euler-larger-polyconvex-Mandel}
\mathscr{M}_\EEps(F,J):=
{\rm dev}\Big(F^\top\!\big([\mathscr{Y}_\EPS\breve\varphi]_F'(F,J){+}\sqrt\DELTA F\big)\Big)\,.
\end{align}\end{subequations}

Then, using again the cut-off $\CUT$ from \eq{cut-off}, we devise the recursive
regularized time-discrete scheme as
\begin{subequations}\label{Euler-large-viscoelastodyn+disc}
\begin{align}\label{Euler-large-viscoelastodyn+0disc}
&\!\!\frac{\varrho_\etau^k{-}\varrho_\etau^{k-1}\!\!}\tau\,=
{\rm div}\big(\DELTA|\nabla\varrho_\etau^k|^{r-2}
\nabla\varrho_\etau^k-\,\CUT(\varrho_\etau^k)\pp_\etau^k\big)\,,
\\[-.2em]&\nonumber
\!\!\frac{\pp_\etau^k{-}\pp_\etau^{k-1}\!\!}\tau\,=
{\rm div}\Big(\TTtauk{+}\DD_\etau^k{-}\CUT(\varrho_\etau^k)\pp_\etau^k{\otimes}\vvk\Big)
+\varrho_\etau^k\GRAVITY_{\DELTA\tau}^k
\\[-.3em]&\nonumber\hspace*{16em}
-\EPS|\vvk|^{p-2}\vvk-\DELTA
|\nabla\varrho_\etau^k|^{r-2}(\nabla\vvk)\nabla\varrho_\etau^k\,,
\\[-.1em]&\nonumber\hspace*{5em}
\text{ where }\ \pp_\etau^k=\varrho_\etau^k\vvk\,,
\quad\text{ with }\ \ \TTtauk=\mathscr{T}_\EEps(\FF_{\!\etau}^k,J_{\etau}^k)\,,
\\[-.2em]&
\hspace*{5em}\text{ and }\ \ \ \DD_\etau^k=\bbD\strain(\vvk)
-{\rm div}\mathfrak{H}_\etau^k\ \ \text{ with }\ \ \mathfrak{H}_\etau^k=\HYPER\big|\nabla^2\vvk\big|^{p-2}\nabla^2\vvk\,,\!
\label{Euler-large-viscoelastodyn+1disc}
\\[-.1em]
&\!\!\frac{\FF_{\!\etau}^k{-}\FF_{\!\etau}^{k-1}\!\!}\tau\,
=(\nabla\vvk)\FF_{\!\etau}^k-(\vvk\Cdot\nabla)\FF_{\!\etau}^k-
\FF_{\!\etau}^k\Lp_{\etau}^k+\DELTA\Delta\FF_{\!\etau}^k\,,
\label{Euler-large-viscoelastodyn+2disc}
\\[-.5em]&\!\!\frac{J_\etau^k{-}J_\etau^{k-1}\!\!}\tau\,
=({\rm div}\,\vvk)J_\etau^k-\vvk\Cdot\nabla J_\etau^k+\DELTA\Delta J_\etau^k\,,\ \ \text{ and}
\label{Euler-large-viscoelastodyn+3disc}
\\[-.1em]\label{Euler-large-viscoelastodyn+4disc}
&\Gm\Lp_{\etau}^k=
{\rm div}\big(\nu|\nabla\Lp_{\etau}^k|^{q-2}\nabla\Lp_{\etau}^k\big)-\MM_{\etau}^k
\ \ \ \text{ with }\ \MM_{\etau}^k=\mathscr{M}_\EEps(\FF_{\!\etau}^k,J_{\etau}^k)
\end{align}\end{subequations}
with the boundary conditions \eq{BC-disc} with $\FF_{\!\etau}^k$ instead of
$\Ee_\etau^k$ and with $\nn\Cdot\nabla J_{\etau}^k=0$ and
$(\nn\Cdot\nabla)\Lp_{\etau}^k=\bm0$ on $I{\times}\varGamma$. This
boundary-value problem is recursive with the initial conditions
\begin{align}\label{IC-large-disc}
\varrho_\etau^0=\varrho_0\,,\ \ \ \ \ \ \ \ \pp_\etau^0=\varrho_0\vv_{0}\,,\  
\ \ \ \ \ {\FF}_{\!\etau}^0=\FF_0\,,  \ \ \text{ and }\ \ \ J_\etau^0=\det\FF_0\,.
\end{align}

It should be emphasized that the calculus \eq{DT-det-extended} does not hold
for time differences even if the regularizing $\Delta$-terms were
omitted in (\ref{Euler-large-viscoelastodyn+disc}c,d), so that
$J_\etau^k=\det\FF_{\!\etau}^k$ is {\it not} granted even if it is valid
for $k=0$, i.e.\ for the ``compatible'' initial conditions as in
\eq{IC-large-disc}. Therefore, to execute a discrete variant of
\eq{ED-Euler-large-stress-calculus}, we must indeed consider two separate
equations, i.e.\ \eq{Euler-large-viscoelastodyn+1disc} for $\FF_{\!\etau}^k$ and
\eq{Euler-large-viscoelastodyn+2disc} for $J_\etau^k$.

\subsection{Stability and convergence of the discrete scheme}\label{sec-large-proof}

We will follow the strategy from Sect.~\ref{sec-linearized-proof} and quite
directly exploit the arguments concerning the convexity of the kinetic energy
$\frac12|\pp|^2/\varrho$. Therefore, the relevant parts will be presented only briefly.
Here, in addition, we have to deal with nonconvexity, blow-up of and the ``non-strong
polyconvexity'' of $\varphi$. We divide the arguments into seven steps.

\medskip{\it Step 1: Choice of $\rhoMAX$}. Exploiting the energy-dissipation balance
\eq{Euler-large-energy-balance-stress} in the same way as in the Step~1 in
Section~\ref{sec-linearized-proof}, we can set $\rhoMAX$.

\medskip{\it Step 2: Basic stability of the scheme \eq{Euler-large-viscoelastodyn+disc}
and first a-priori estimates}.
The convexity of $\breve\varphi$ and thus the strong convexity of $\breve\varphi_\EEps$
is then used in the discrete variant of \eq{ED-Euler-large-stress-calculus}
as an inequality
\begin{align}\nonumber
&\int_\varOmega
\frac{\breve\varphi_\EEps(\FF_{\!\etau}^k,J_\etau^k){-}\breve\varphi_\EEps(\FF_{\!\etau}^{k-1},J_\etau^{k-1})}
\tau\,\d\xx
\le
\int_\varOmega\bigg(\Big([\mathscr{Y}_\EPS\breve\varphi]_\FF'(\FF_{\!\etau}^k,J_\etau^k)
+\sqrt\DELTA\FF_{\!\etau}^k\Big)\Colon\frac{\FF_{\!\etau}^k{-}\FF_{\!\etau}^{k-1}\!\!}\tau
\\[-.5em]\nonumber
&\hspace{19em}+
\Big([\mathscr{Y}_\EPS\breve\varphi]_J'(\FF_{\!\etau}^k,J_\etau^k)
+\sqrt\DELTA J_{\!\etau}^k\Big)\frac{J_\etau^k{-}J_\etau^{k-1}\!}\tau\,\bigg)\,\d\xx
\\[-.3em]\nonumber&
=\int_\varOmega\bigg(\Big([\mathscr{Y}_\EPS\breve\varphi]_\FF'(\FF_{\!\etau}^k,J_\etau^k)
+\sqrt\DELTA\FF_{\!\etau}^k\Big)\Colon
\Big((\nabla\vvk)\FF_{\!\etau}^k-(\vvk\Cdot\nabla)\FF_{\!\etau}^k-
\FF_{\!\etau}^k\Lp_{\etau}^k+\DELTA\Delta\FF_{\!\etau}^k\Big)
\\[-.5em]\nonumber&\hspace{5em}
+\Big([\mathscr{Y}_\EPS\breve\varphi]_J'(\FF_{\!\etau}^k,J_\etau^k)
+\sqrt\DELTA J_{\!\etau}^k\Big)
\Big(({\rm div}\,\vvk)J_\etau^k-\vvk\Cdot\nabla J_\etau^k+\DELTA\Delta J_\etau^k\Big)
\bigg)\,\d\xx
\\[-.1em]\nonumber&
=-\int_\varOmega\bigg(\hspace{-.7em}\lineunder{\big[\DELTA[\mathscr{Y}_\EPS\breve\varphi]''(\FF_{\!\etau}^k,J_\etau^k)\big]
\Big(\!\!\begin{array}{c}\nabla\FF_{\!\etau}^k\\\nabla J_{\etau}^k\end{array}\!\!\Big)
\Vdots\Big(\!\!\begin{array}{c}\nabla\FF_{\!\etau}^k\\\nabla J_{\etau}^k\end{array}\!\!\Big)}{$\ge0$}
\\[-1.2em]\nonumber&\hspace{15em}
+\DELTA^{3/2}|\nabla\FF_{\!\etau}^k|^2\!+\DELTA^{3/2}|\nabla J_{\etau}^k|^2\!
+\Gm\big|\Lp_{\etau}^k\big|^2\!+\nu\big|\nabla\Lp_{\etau}^k\big|^q\bigg)
\,\d\xx\,.
\end{align}
Here, like in \eq{Euler-small-calc-Delta}, we have formally used the 2nd-order
derivative $[\mathscr{Y}_\EPS\breve\varphi]''$ but the convex $C^1$-function
$\mathscr{Y}_\EPS\breve\varphi$ yields such estimate by a smoothening argument, too.

Thus, like \eq{Euler-ED-basic-energy-balance-disc}, we obtain the discrete
(and regularized) analog of \eq{Euler-large-energy-balance-stress}, specifically
\begin{align}\nonumber
\!\!\!&\int_\varOmega\!\frac{|\pp_\etau^k|^2\!}{2\varrho_\etau^k\!}
+\breve\varphi_\EEps(\FF_{\!\etau}^k,J_\etau^k)\,\d\xx
+\tau\sum_{l=1}^k\int_\varOmega\!\bigg(\bbD\strain(\vv_\etau^l)\Colon\strain(\vv_\etau^l)
+\Gm\big|\Lp_{\etau}^k\big|^2+\HYPER|\nabla^2\vv_\etau^l|^p\!
\\[-.7em]&\hspace{15em}\nonumber
+\nu\big|\nabla\Lp_{\etau}^k\big|^q+\DELTA|\vv_\etau^l|^p
+\DELTA^{3/2}|\nabla\FF_{\!\etau}^l|^2+\DELTA^{3/2}|\nabla J_{\etau}^l|^2\bigg)\d\xx
\\[-.7em]&
\hspace{7em}
\le\int_\varOmega\!
\frac{|\pp_0|^2\!}{2\varrho_0\!}+\breve\varphi(\FF_0,J_0)
+\frac{\sqrt\DELTA}2|\FF_0|^2+\frac{\sqrt\DELTA}2J_0^2\,\d\xx
+\tau\sum_{l=1}^k
\int_\varOmega\!\varrho_\etau^l\GRAVITY_{\DELTA\tau}^l\Cdot\vv_\etau^l\,\d\xx\,;
\label{Euler-large-energy-balance-disc}\end{align}
note that we have taken into account the convexity of
$[\mathscr{Y}_\EPS\breve\varphi]$ which grants the positive semi-definiteness of
$[\mathscr{Y}_\EPS\breve\varphi]''(\FF_{\!\etau}^k,J_\etau^k)$ so that this term was
simply omitted and also that $\mathscr{Y}_\EPS\breve\varphi\le\breve\varphi$ so that
surely $[\mathscr{Y}_\EPS\breve\varphi](\FF_0,J_0)\le\breve\varphi(\FF_0,J_0)$.
Again, we have \eq{Euler-small-est-1+} and, like \eq{Euler-small-est-Gronwall},
we can see the a-priori estimates (\ref{Euler-small-est}a,b,e) and now also
\begin{subequations}\label{Euler-large-est}\begin{align}
&\label{Euler-large-est3}
\|\overlineFetau\|_{L^\infty(I;L^2(\varOmega;\R^{3\times3}))}^{}\le C\,, \ \ 
\ \ \ \,
\|\nabla\overlineFetau\|_{L^2(I{\times}\varOmega;\R^{3\times3\times3})}^{}\le\frac{C}{\DELTA^{3/4}}\,,
\\&\label{Euler-large-est3+}
\|\overlineJetau\|_{L^\infty(I;L^2(\varOmega))}^{}\le C\,,
\hspace{4em}\|\nabla\overlineJetau\|_{L^2(I{\times}\varOmega;\R^{3})}^{}\le \frac{C}{\DELTA^{3/4}}\,,
\\&\label{Euler-large-est3++}
\|\overlineLpetau\|_{L^2(I{\times}\varOmega;\R^{3\times3}))}^{}\le C\,,\ \text{ and }\ 
\|\nabla\overlineLpetau\|_{L^q(I{\times}\varOmega;;\R^{3\times3\times3}))}^{}\le C\,,
\intertext{and, by the Gagliardo-Nirenberg interpolation between
these estimates in \eq{Euler-large-est3} and similarly in
\eq{Euler-large-est3+}, we have also}\label{Euler-large-est3-10/3}
&\|\overlineFetau\|_{L^{10/3}(I\times\varOmega;\R^{3\times3}))}^{}\le C_\EEps
\ \ \ \text{ and }\ \ \
\|\overlineJetau\|_{L^{10/3}(I\times\varOmega))}^{}\le C_\EEps\,.
\intertext{From this and using the at-most linear growth
of $[\mathscr{Y}_\EPS\breve\varphi]'$, we have also }
&\|\overlineTetau\|_{L^\infty(I;L^1(\varOmega;\R^{3\times3}))\,\cap\,L^{5/3}(I\times\varOmega;\R^{3\times3})}^{}\!
\le C_\EEps\ \:\text{ and }\ \:
\|\overlineMetau\|_{L^\infty(I;L^1(\varOmega;\R^{3\times3}))\,\cap\,L^{5/3}(I\times\varOmega;\R^{3\times3})}^{}\!
\le C_\EEps\,.\!
\end{align}\end{subequations}

Like in Section~\ref{sec-large-proof}, except the latter estimates in
\eq{Euler-small-est2} and in \eq{Euler-large-est3++}, all the above estimates show
the stability of the discrete regularized scheme even for the simple-material variant
with $\HYPER=0$ and $\nu=0$. However, for the rest, we will need the non-simple
material model with $\HYPER>0$ and $\nu>0$ and $\min(p,q)>3$.

\medskip{\it Step 3: Further a-priori estimates relying on $\HYPER>0$ and
on $\nu>0$}.
By the same arguments as in Step~3 in Sect.\,\ref{sec-linearized-proof}, we can show
also 
\eq{Euler-small-est+} as well as the existence
of $(\varrho_\etau^k,\vv_\etau^k,\FF_\etau^k,J_\etau^k,\Lp_\etau^k)$.
Also \eq{Euler-small-est6} is directly valid here, too.
Using the above estimate of $\overlineTetau$ in $L^\infty(I;L^1(\varOmega;\R^{3\times3}))$,
we also prove \eq{Euler-small-est6+}. Like \eq{Euler-small-est7}, we have now
\begin{subequations}\label{Euler-large-est+}\begin{align}
\nonumber
&\Big\|\pdt{\FF_{\!\etau}}\Big\|_{L^{6/5}(I;H^1(\varOmega;\R^{3\times3})^*)}^{}
=\!\!\sup_{\|\wt\FF\|_{L^{5}(I;H^1(\varOmega;\R^{3\times3}))}\le1}\int_0^T\!\!\!\int_\varOmega\bigg(
\DELTA\nabla\overlineFetau\Vdots\nabla\wt\FF
\\[-.2em]&\hspace{11em}
+\Big((\overlinevvtau\Cdot\nabla)\overlineFetau\!
-(\nabla\overlinevvtau)\overlineFetau\!
+
\overlineFetau\overlineLpetau
\Big)\Colon\wt\FF\bigg)\,\d\xx\d t\le C_\EEps\,,
\intertext{where we used the above estimates
(\ref{Euler-small-est}a,b) and (\ref{Euler-large-est}a,c), and thus in particular
$\overlineFetau\overlineLpetau$ in $L^2(I{\times}\varOmega;\R^{3\times3})$.
Similarly,}
&\Big\|\pdt{J_{\etau}}\Big\|_{L^{6/5}(I;H^1(\varOmega)^*)}^{}\le C_\EEps\,.
\end{align}\end{subequations}

The last estimate allows us to improve \eq{Euler-large-est3++} when realizing
that $\overlineLptau$ solves the elliptic problem
${\rm div}\big(\nu|\nabla\overlineLptau|^{q-2}\nabla\overlineLptau\big)
-\Gm\overlineLptau=\overlineMetau$ with the boundary condition
$(\nn\Cdot\nabla)\overlineLptau=\bm0$ on $I{\times}\varGamma$, cf.\
\eq{Euler-large-viscoelastodyn+4disc}. Namely,
\begin{align}\label{Euler-large-est-Lp}
\|\overlineLptau\|_{L^\infty(I;W^{1,q}(\varOmega;\R^{3\times3}))}^{}\le C_\EEps\,.
\end{align}

\medskip{\it Step 4: Limit passage for $\tau\to0$}.
By the Banach selection principle, like \eq{Euler-small-converge}, we obtain a
subsequence converging weakly* with respect to the topologies indicated in
(\ref{Euler-small-est}a,b,e), \eq{Euler-small-est+}, \eq{Euler-large-est},
\eq{Euler-large-est+}, and \eq{Euler-large-est-Lp} to some limit
$(\varrho_\EEps,\pp_\EEps,\vv_\EEps,\FF_\EEps,J_\EEps,\Lp_\EEps)$. Specifically,
beside (\ref{Euler-small-converge}a--e) we have now also
\begin{subequations}\label{Euler-large-converge}
\begin{align}
&&&\!\!\overlineFetau\to\FF_\EEps\!\!\!&&\text{weakly* in $\ 
L^2(I;H^1(\varOmega;\R^{3\times3}))\,\cap\,L^\infty(I;L^2(\varOmega;\R^{3\times3}))$\,,}\!\!
\\
&&&\!\!\FF_\etau\to\FF_\EEps\!\!\!&&\text{weakly\ \;in $\ 
L^2(I;H^1(\varOmega;\R^{3\times3}))\,\cap\,W^{1,6/5}(I;H^1(\varOmega;\R^{3\times3})^*)$\,.}&&
\\
&&&\!\!\overlineJetau\to J_\EEps\!\!\!&&\text{weakly* in $\ 
L^2(I;H^1(\varOmega))\,\cap\,L^\infty(I;L^2(\varOmega))$\,,}\!\!
\\
&&&\!\!J_\etau\to J_\EEps\!\!\!&&\text{weakly\ \;in $\ L^2(I;H^1(\varOmega))
\,\cap\,W^{1,6/5}(I;H^1(\varOmega)^*)$\,, and}&&
\\\label{Euler-large-converge-Lp}
&&&\!\!\overlineLpetau\to\Lp_\EEps\!\!\!&&\text{weakly* in $\
L^\infty(I;W^{1,q}(\varOmega;\R^{3\times3}))$\,.}&&
\end{align}\end{subequations}

By the compact embedding $W^{1,r}(\varOmega)\subset C(\barOmega)$ and the
(generalized) Aubin-Lions theorem, relying on the compact embedding
$W^{1,r}(\varOmega)\subset C(\barOmega)$ for $r>3$,
we have again (\ref{Euler-small-converge-strong}a--c,e),
(\ref{Euler-small-converge-strong+}d,e), \eq{Euler-small-nabla-rho-strong}
and now also
\begin{subequations}\label{Euler-large-converge-strong}
\begin{align}\label{Euler-large-converge-strong-F}
&&&\!\!\overlineFetau\to\FF_\EEps\!\!\!\!\!&&\text{strongly in }\
L^a(I;L^2(\varOmega;\R^{3\times3}))\ \text{ and}
\\&&&\!\!\overlineJetau\to J_\EEps\!\!\!\!\!&&\text{strongly in }\
L^a(I;L^2(\varOmega))\ \text{ for any }\ 1\le a<\infty\,.
\end{align}\end{subequations}
Like (\ref{Euler-small-converge-strong}a,b), we have now
\begin{subequations}\label{Euler-large-converge-strong+}
\begin{align}\label{Euler-large-converge-strong-T}
&\!\!\overlineTetau\to\mathscr{T}_\EEps(\FF_\EEps,J_\EEps)\!\!\!\!\!\!
&&\text{strongly in }\ L^a(I;L^1(\varOmega;\R^{3\times3}))\ \text{ and}
\\&\!\!\overlineMetau\to\mathscr{M}_\EEps(\FF_\EEps,J_\EEps)\!\!\!\!\!\!
&&\text{strongly in }\ L^a(I;L^1(\varOmega;\R^{3\times3}))\ \text{ for any }\ 1\le a<\infty\,.
\label{Euler-large-converge-strong-M}
\intertext{From the uniform monotonicity
of the quasi-linear operator $\Lp\mapsto\Gm\Lp-{\rm div}(\nu|\nabla\Lp|^{q-2}\nabla\Lp)$
and from \eq{Euler-large-converge-strong-M} when realizing that
$L^a(I;L^1(\varOmega;\R^{3\times3}))\subset L^a(I;W^{1,q}(\varOmega;\R^{3\times3})^*)$
for $q>3$, we can improve the convergence \eq{Euler-large-converge-Lp} to}
&\!\!\overlineLpetau\to\Lp_\EEps&&\text{strongly in }\
L^a(I;W^{1,q}(\varOmega;\R^{3\times3}))\ \ \text{ for any }\ 1\le a<\infty\,.
\label{Euler-large-converge-strong-Lp}
\end{align}\end{subequations}

We can then repeat the argumentation
\eq{Euler-small-strong-hyper+}--\eq{Euler-small-strong-hyper+++} to make a limit
passage in the momentum equation. Based on \eq{Euler-large-converge-strong-T}, 
we can pass to the limit in the equation for $\FF$-tensor. The limit passage in
the kinematic equation for $J$ is even easier. The convergence in the quasi-static
equation for the inelastic distortion rate follow easily from
(\ref{Euler-large-converge-strong+}b,c).

The limit $(\varrho_\EEps,\pp_\EEps,\vv_\EEps,\FF_\EEps,J_\EEps,\Lp_\EEps)$ solves 
the $(\EPS,\DELTA)$-regularized system:
\begin{subequations}\label{ED-Euler-large-polyconvex-eps}\begin{align}
\label{ED-Euler-large0-polyconvex-eps}
&\pdt{\varrho_\EEps}={\rm div}\big(\DELTA|\nabla\varrho_\EEps|^{p-2}\nabla\varrho_\EEps
-\CUT(\varrho_\EEps)\pp_\EEps\big)\,,\\
\nonumber
     &\pdt{\pp_\EEps}={\rm div}\big(\TT_\EEps+\DD_\EEps-\CUT(\varrho_\EEps)\pp_\EEps{\otimes}\vv_\EEps\big)
     +\varrho_\EEps\GRAVITY_\DELTA
    -\EPS|\vv_\EEps|^{p-2}\vv_\EEps
     -\DELTA|\nabla\varrho_\EEps|^{p-2}(\nabla\vv_\EEps)\nabla\varrho_\EEps
 \\[-.1em]\nonumber
    &\hspace*{4em}\text{ with }\ \pp_\EEps=\varrho_\EEps\vv_\EEps\,,\ \ \text{where }\,
    \TT_\EEps=\mathscr{T}_\EEps(\FF_\EEps,J_\EEps)
\\&\hspace*{4em}
\text{ and }\ \ \ \,\DD_\EEps=\bbD\strain(\vv_\EEps)
 -{\rm div}\,\mathfrak{H}_\EEps
 \ \ \text{ with }\ \ 
      \mathfrak{H}_\EEps=\HYPER|\nabla^2\vv_\EEps|^{p-2}\nabla^2\vv_\EEps\,,
\label{ED-Euler-large1-polyconvex-eps}
\\[-.3em]
&\pdt{\FF_\EEps}=(\nabla\vv_\EEps)\FF_\EEps-(\vv_\EEps\Cdot\nabla)\FF_\EEps-
\FF_\EEps\Lp_\EEps+\DELTA\Delta\FF_\EEps\,,
 \label{ED-Euler-large2-polyconvex-eps}
\\[-.1em]\label{ED-Euler-large3-polyconvex-eps}
&\pdt{J_\EEps}=({\rm div}\,\vv_\EEps)J_\EEps-\vv_\EEps\Cdot\nabla J_\EEps+\DELTA\Delta J_\EEps\,,
\\[-.1em]\label{ED-Euler-large4-polyconvex-eps}
&
\Gm\Lp_\EEps={\rm div}\big(\nu|\nabla\Lp_\EEps|^{q-2}\nabla\Lp_\EEps\big)-\MM_\EEps
   \ \text{ with }\
\MM_\EEps=\mathscr{M}_\EEps(\FF_\EEps,J_\EEps)
\end{align}\end{subequations}
again with the boundary conditions \eq{Euler-small-BC-hyper} written
for $(\vv_\EEps,\TT_\EEps,\DD_\EEps,\mathfrak{H}_\EEps)$ together with
 $\nn\Cdot\nabla\varrho_\EEps=0$, $(\nn\Cdot\nabla)\FF_\EEps=\bm0$,
$\nn\Cdot\nabla J_\EEps=0$, and $(\nn\Cdot\nabla)\Lp_\EEps=\bm0$ on $I{\times}\varGamma$.

\medskip{\it Step 5: Elimination of the cut-off $\CUT$ and further a-priori estimates}.
By the same arguments as in Step~5 in Section~\ref{sec-linearized-proof}, we can
forget $\CUT(\varrho_\EEps)\equiv1$ in (\ref{ED-Euler-large-polyconvex-eps}a,b).

Having the time-continuous problem, we can now obtain a-priori estimates independent
of $\EPS$ and $\DELTA$ when exploiting the Gronwall inequality in a standard way.
Then we can take the same argumentation leading to \eq{Euler-small-est-rho}.
Further, we can improve also the estimates \eq{Euler-large-est3} for
$\FF$ when testing \eq{ED-Euler-large2-polyconvex-eps} by $\Delta\FF_\EEps$.
In contrast to \eq{Euler-small-est-E}, we could not rely on a mere
$L^2(I{\times}\varOmega)$-bound for $\Lp_\EEps$, which is why we need to involve the
gradient modification of \eq{ED-Euler-large3} with $\nu>0$. Relying on the
$L^2(I;W^{1,q}(\varOmega))$-bound, cf.\ \eq{Euler-large-est3++}, we can execute
the test of \eq{ED-Euler-large2-polyconvex-eps} by $\Delta\FF_\EEps$
and, exploiting the embedding $W^{1,q}(\varOmega)\subset L^\infty(\varOmega)$ with the
norm $N_q$ for $q>3$, modify the arguments in \eq{Euler-small-est-E-} as 
\begin{align}\nonumber
&\!\int_\varOmega\nabla(\FF_\EEps\Lp_\EEps)\Vdots\nabla\FF_\EEps\,\d\xx
=\int_\varOmega\big(\nabla\FF_\EEps\Cdot\Lp_\EEps+\FF_\EEps\Cdot\nabla\Lp_\EEps\big)
\Vdots\nabla\FF_\EEps\,\d\xx
\\[-.1em]&\nonumber\quad\le\! N_q\|\Lp_\EEps\|_{L^\infty(\varOmega;\R^{3\times3})}^{}
\|\nabla\FF_\EEps\|_{L^2(\varOmega;\R^{3\times3\times3})}^2\!
\\[-.1em]&\nonumber\hspace{13.3em}
+\|\FF_\EEps\|_{L^6(\varOmega;\R^{3\times3})}^{}\|\nabla\Lp_\EEps\|_{L^3(\varOmega;\R^{3\times3\times3})}^{}
\|\nabla\FF_\EEps\|_{L^2(\varOmega;\R^{3\times3\times3})}^{}\ \ \ \ 
\\[-.1em]&\nonumber\quad\le\! N_q\|\Lp_\EEps\|_{L^\infty(\varOmega;\R^{3\times3})}^{}
\|\nabla\FF_\EEps\|_{L^2(\varOmega;\R^{3\times3\times3})}^2\!
\\[-.2em]&\hspace{2em}+C\Big(\|\FF_\EEps\|_{L^2(\varOmega;\R^{3\times3})}^{}\!+
\|\nabla\FF_\EEps\|_{L^2(\varOmega;\R^{3\times3\times3})}^{}\Big)
\|\nabla\Lp_\EEps\|_{L^3(\varOmega;\R^{3\times3\times3})}^{}
\|\nabla\FF_\EEps\|_{L^2(\varOmega;\R^{3\times3\times3})}^{}
\label{Euler-large-est-F-}\end{align}
at each time instant $t\in I$. This can then be treated by the Gronwall inequality,
based on the already obtained former estimates in \eq{Euler-large-est3}
and \eq{Euler-large-est3++} inherited for $\FF_\EEps$ and $\Lp_\EEps$. The terms
$(\nabla\vv_\EEps)\FF_\EEps$ and $(\vv_\EEps\Cdot\nabla)\FF_\EEps$ tested by
$\Delta\FF_\EEps$ should be treated by the Green formula as in \cite{Roub22QHLS}.
This yields the estimate analogous to \eq{Euler-small-est-E}, i.e.
\begin{subequations}\label{Euler-large-est++}\begin{align}
\label{Euler-large-est++F}
&\|\FF_\EEps\|_{L^\infty(I;H^1(\varOmega;\R^{3\times3}))}^{}\le C\ \ \ \text{ and }\ \ \
\|\Delta\FF_\EEps\|_{L^2(I\times\varOmega;\R^{3\times3}))}^{}\le\frac C{\sqrt{\DELTA}}\,.
\intertext{Testing \eq{ED-Euler-large3-polyconvex-eps} by $\Delta J_\EEps$ is even
simpler and yields}
&\label{Euler-large-est++J}
\|J_\EEps\|_{L^\infty(I;H^1(\varOmega;\R^{3\times3}))}^{}\le C\ \ \ \text{ and }\ \ \
\|\Delta J_\EEps\|_{L^2(I\times\varOmega;\R^{3\times3}))}^{}\le\frac C{\sqrt{\DELTA}}\,.
\end{align}\end{subequations}
Further estimates \eq{Euler-small-est-rho-E} now hold, too, for $\FF_\EEps$ and
$J_\EEps$ in place of $\Ee_\EEps$.

\medskip{\it Step 6: Limit passage for $\DELTA\to0$}.
By the Banach selection principle, we obtain a subsequence converging weakly* with
respect to the topologies specified in Step~5
to some limit $(\varrho_\EPS,\pp_\EPS,\vv_\EPS,\FF_{\!\EPS},J_\EPS)$. Specifically, we
have the convergences (\ref{Euler-small-converge+}a--c,e--f,h,i) and also
(\ref{Euler-small-converge+}d,j) hold for $\FF_\EEps$ and
$J_\EEps$ in place of $\Ee_\EEps$.

In contrast to \eq{Euler-small-est3} showing the blow-up rate
$\mathscr{O}(\DELTA^{-1/2})$, we have now a faster blow-up rate
$\mathscr{O}(\DELTA^{-3/4})$ in (\ref{Euler-large-est}a,b) inherited for
$(\FF_{\!\EEps},J_\EEps)$. This will ensure convergence in the pressure term
$\breve\varphi_\EEps(\FF_{\!\DELTA},J_\DELTA)\bbI$ in the Cauchy stress with the
regularized $\breve\varphi_\EEps$ from \eq{varphi-regularized}. More in detail,
we use the interpolation between the $L^\infty(I;L^2(\varOmega))$- and
$L^2(I;H^1(\varOmega))$-estimates \eq{Euler-large-est3} with the weight $1/2$. This
gives the $L^4(I;L^3(\varOmega))$-estimates, i.e.\ here 
\begin{align}
\|\FF_\EEps\|_{L^4(I;L^3(\varOmega;\R^{3\times3}))}^{}\le\frac{C}{\DELTA^{3/8}}
\ \ \:\text{ and }\ \ \:\|J_\EEps\|_{L^4(I;L^3(\varOmega))}^{}\le \frac{C}{\DELTA^{3/8}}\,.
\label{Euler-large-est3*}\end{align}
Due to these estimates, we have
$\|\sqrt\DELTA|\FF_{\!\EEps}|^2{+}\sqrt\DELTA J_\EEps^2\|_{L^2(I;L^{3/2}(\varOmega))}
=\mathscr{O}(\DELTA^{1/8})\to0$ for $\DELTA\to0$. Recalling the at-most quadratic
growth of $\mathscr{Y}_\EPS\breve\varphi$, 
we have the convergence in the pressure in the Cauchy stress
\begin{align*}
\!\breve\varphi_\EEps(\FF_{\!\EEps},J_\EEps)
=[\mathscr{Y}_\EPS\breve\varphi](\FF_{\!\EEps},J_\EEps)
{+}\frac{\sqrt\DELTA}2|\FF_{\!\EEps}|^2{+}\frac{\sqrt\DELTA\!}2J_\EEps^2
\to
[\mathscr{Y}_\EPS\breve\varphi](\FF_{\!\DELTA},J_\DELTA)
\ \text{ strongly in }\,L^2(I;L^{3/2}(\varOmega))\,.
\end{align*}
The other regularizing terms in \eq{Euler-larger-polyconvex-Cauchy}, i.e.\
$\sqrt\DELTA\FF_{\!\DELTA}\FF_{\!\DELTA}^\top$ and $\sqrt\DELTA J_\DELTA^2\bbI$, converge to
$\bm0$ by the same arguments, while the remaining term in
\eq{Euler-larger-polyconvex-Cauchy} are continuous. Similar arguments apply to the
Mandel stress \eq{Euler-larger-polyconvex-Mandel}. Altogether, when defining
\begin{subequations}\label{Euler-larger-polyconvex-stress-}
\begin{align}\label{Euler-larger-polyconvex-Cauchy+}
&\mathscr{T}_\EPS(F,J):=[\mathscr{Y}_\EPS\breve\varphi]_F'(F,J)F^\top\!
+\Big(J[\mathscr{Y}_\EPS\breve\varphi]_J'(F,J){+}[\mathscr{Y}_\EPS\breve\varphi](F,J)\Big)\bbI
\ \ \ \text {and}
\\&\label{Euler-larger-polyconvex-Mandel+}
\mathscr{M}_\EPS(F,J):={\rm dev}\big(F^\top\![\mathscr{Y}_\EPS\breve\varphi]_F'(F,J)\big)\,,
\end{align}\end{subequations}
and when realizing the at-most linear growth of
$[\mathscr{Y}_\EPS\breve\varphi]'$, we summarize that
\begin{subequations}\label{Euler-larger-polyconvex-stress+}
\begin{align}
&&&\TT_\EEps=\mathscr{T}_\EEps(\FF_{\!\EEps},J_\EEps)\to\mathscr{T}_\EPS(\FF_{\!\EPS},J_\EPS)&&
\text{strongly in }\ L^2(I;L^{3/2}(\varOmega;\R^{3\times3}))\ \text{ and}&&
\\
&&&\MM_\EEps=\mathscr{M}_\EEps(\FF_{\!\EEps},J_\EEps)
\to\mathscr{M}_\EPS(\FF_{\!\EPS},J_\EPS)\!\!\!&&
\text{strongly in }\ L^2(I;L^{3/2}(\varOmega;\R^{3\times3}))\,.
\end{align}\end{subequations}

Thus, the limit $(\varrho_\EPS,\pp_\EPS,\vv_\EPS,\FF_{\!\EPS},J_\EPS,\Lp_\EPS)$ solves the
$\EPS$-regularized system
\begin{subequations}\label{ED-Euler-large-polyconvex-eta}\begin{align}
\label{ED-Euler-large0-polyconvex-eta}
&\pdt{\varrho_\EPS}=-{\rm div}\,\pp_\EPS\,,\\
\nonumber
&\pdt{\pp_\EPS}={\rm div}\big(\TT_\EPS+\DD_\EPS-\pp_\EPS{\otimes}\vv_\EPS\big)
     +\varrho_\EPS\GRAVITY -\EPS|\vv_\EPS|^{p-2}\vv_\EPS
   \ \text{ with }\ \pp_\EPS=\varrho_\EPS\vv_\EPS\,,\ \ \text{where }\,
\\&\hspace*{4em}
\TT_\EPS=\mathscr{T}_\EPS(\FF_\EPS,J_\EPS)\ \text{ and }\ \DD_\EPS=\bbD\strain(\vv_\EPS)
 -{\rm div}\,\mathfrak{H}_\EPS
 \ \text{ with }\ 
      \mathfrak{H}_\EPS=\HYPER|\nabla^2\vv_\EPS|^{p-2}\nabla^2\vv_\EPS\,,\!\!
\label{ED-Euler-large1-polyconvex-eta}
\\[-.3em]
&\pdt{\FF_\EPS}=(\nabla\vv_\EPS)\FF_\EPS-(\vv_\EPS\Cdot\nabla)\FF_\EPS-\FF_\EPS\Lp_\EPS\,,
 \label{ED-Euler-large2-polyconvex-eta}
\\[-.1em]\label{ED-Euler-large3-polyconvex-eta}
&\pdt{J_\EPS}=({\rm div}\,\vv_\EPS)J_\EPS-\vv_\EPS\Cdot\nabla J_\EPS\,,
\\[-.1em]\label{ED-Euler-large4-polyconvex-eta}
&
\Gm\Lp_\EPS={\rm div}\big(\nu|\nabla\Lp_\EPS|^{q-2}\nabla\Lp_\EPS\big)-\MM_\EPS
   \ \ \ \text{ with }\ \
\MM_\EPS=\mathscr{M}_\EPS(\FF_\EPS,J_\EPS)
\end{align}\end{subequations}
with the boundary conditions \eq{Euler-small-BC-hyper} written for $\vv_\EPS$,
$\TT_\EPS$, $\DD_\EPS$, and $\mathfrak{H}_\EPS$ together with
$(\nn\Cdot\nabla)\Lp_\EEps=\bm0$ on $I{\times}\varGamma$. Noteworthy, the equation
\eq{ED-Euler-large3-polyconvex-eta} can now be eliminated when realizing that
$J_\EPS=\det\FF_\EPS$ is now inherited during the whole evolution from the
initial conditions $J|_{t=0}=\det\FF_0$, cf.\ \eq{IC-large+}.

Like \eq{Euler-large-energy-balance-stress}, we have now the energy-dissipation
balance
\begin{align}\nonumber
\frac{\d}{\d t}&\int_\varOmega\frac{\varrho_\EPS}2|\vv_\EPS|^2
+[\mathscr{Y}_\EPS\breve\varphi](\FF_\EPS,J_\EPS)\,\d\xx
\\[-.3em]&\quad
+\!\int_\varOmega\!\bbD\strain(\vv_\EPS)\Colon\strain(\vv_\EPS)+\HYPER|\nabla^2\vv_\EPS|^p\!
+\Gm\big|\Lp_\EPS\big|^2\!+\nu\big|\nabla\Lp_\EPS\big|^q\!+\EPS|\vv_\EPS|^q\,\d\xx
=\int_\varOmega\!\varrho_\EPS\GRAVITY\Cdot\vv_\EPS\,\d\xx\,.
\label{Euler-large-energy-balance-stress+}\end{align}

\medskip{\it Step 7: Limit passage for $\EPS\to0$ and the energy-dissipation balance}.
From \eq{Euler-large-energy-balance-stress+}, we get the estimates for $\nabla\vv_\EPS$
uniform with respect to $\EPS>0$ and thus we have also the ``strong positivity''
of $\varrho_\EPS$ leading to \eq{basic-est-of-v}.

In view of the qualification of $\FF_0$ in \eq{Euler-small-ass-IC} with $r>3$,
applying \cite[Lemma\,5.1]{Roub24TVSE} in combination with \cite[Sect.4]{Roub22QHLS} to
\eq{ED-Euler-large2-polyconvex-eta}, we get the bound for $\FF_\EPS$ in
$L^\infty(I;W^{1,r}(\varOmega;\R^{3\times3}))$; here we need $\nu>0$ and $q>3$ in the
gradient term in \eq{ED-Euler-large4-polyconvex-eta} to execute the test of
$\FF_\EPS\Lp_\EPS$ by ${\rm div}(|\nabla\FF_\EPS|^{r-2}\nabla\FF_\EPS)$ as in
\cite{Roub22QHLS}. By this way, we obtain also the $L^\infty(I;W^{1,r}(\varOmega))$-bound
for $J_\EPS$ solving \eq{ED-Euler-large3-polyconvex-eta}.

We now replicate the arguments in Step~6, obtaining a limit
$(\varrho,\pp,\vv,\FF,J,\Lp)$.
Counting also an information about time-derivative $\pdt{}\FF_\EPS$, we have 
\begin{align}
\FF_\EPS\to\FF\quad\text{strongly in }\ C(I{\times}\barOmega;\R^{3\times3})
\end{align}
and thus also $J_\EPS=\det\FF_\EPS\to J=\det\FF$ strongly in
$C(I{\times}\barOmega)$.

Having the regularizing $\DELTA$-Laplacean in \eq{ED-Euler-large3-polyconvex-eps}
dismissed in \eq{ED-Euler-large3-polyconvex-eta}, we can read also the
equation for $1/J_\EPS$, namely: 
\begin{align}\nonumber\\[-2em]\label{1/J-continuity-eq}
&\pdt{}\frac1{J_\EPS}=-({\rm div}\,\vv_\EPS)\frac1{J_\EPS}-\vv_\EPS\Cdot\nabla\frac1{J_\EPS}\,,
\end{align}
cf.\ \eq{DT-det-inverse}. Applying again \cite[Lemma\,5.1]{Roub24TVSE} now to
\eq{1/J-continuity-eq}, we obtain a bound for $1/J_\EPS$ in $L^\infty(I;W^{1,r}(\varOmega))$
uniform with respect to $\EPS$. Realizing the embedding $W^{1,r}(\varOmega)\subset
L^\infty(\varOmega))$ for $r>3$, this bound implies that
\begin{align}\nonumber\\[-2em]\label{J-min}
\inf_{\EPS>0}\min_{I\times\barOmega}J_\EPS=:J_{\min{}}>0\,.
\end{align}
This effectively eliminates the singularity of $\breve\varphi$ at $J=0$ as
assumed in \eq{Euler-large-ass-phi}. As $r>3$, we can see that $(\FF_\EPS,J_\EPS)$ is
valued in a compact set in ${\rm GL}_3^+\times\R^+$. On the convex hull of this set,
$\breve\varphi$ is convex and continuously differentiable 
and we have $\mathscr{Y}_\EPS\breve\varphi\to\breve\varphi$ as well as
$[\mathscr{Y}_\EPS\breve\varphi]'\to\breve\varphi'$ uniformly for $\EPS\to0$,
cf.\ \cite[Chap.3, Prop.2.29 and 4.33]{HuPap97HMAT}.

We have improved the mode of the convergence \eq{Euler-larger-polyconvex-stress+}
of stresses, i.e.\ now as
\begin{subequations}\label{Euler-larger-polyconvex-stress++}
\begin{align}
&\mathscr{T}_\EPS(\FF_{\!\EPS},J_\EPS)\to\breve\varphi_\FF'(\FF,J)\FF^\top\!
+\Big(J\breve\varphi_J'(\FF,J){+}\breve\varphi(\FF,J)\Big)\bbI
\ \ 
\text{strongly in }\,C(I{\times}\barOmega;\R^{3\times3})\,\text{ and}&&
\\
&\mathscr{M}_\EPS(\FF_{\!\EPS},J_\EPS)\to{\rm dev}\big(\FF^\top\!\breve\varphi_\FF'(\FF,J)\big)\ \ 
\text{ strongly in }\ C(I{\times}\barOmega;\R^{3\times3})\
\text{ with }\ J=\det\FF\,.
\end{align}\end{subequations}

In view of \eq{Cauchy-Mandel-stress-polyconvex}, we showed that
$(\varrho,\vv,\FF,\Lp)$ is a weak solution due to Definition~\ref{def-ED-Ch5}.
Satisfaction of \eq{Euler-large-energy-balance-stress}
integrated over the time interval $I$ then copies the arguments from Step~8 in
Section~\ref{sec-linearized-proof}, noting that, in addition, the inelastic
flow rule \eq{ED-Euler-large3} is satisfied in the sense of
$L^2(I{\times}\varOmega;\R_{\rm dev}^{3\times3})$ and its test by
$\Lp\in L^2(I;W^{1,q}(\varOmega;\R_{\rm dev}^{3\times3}))$ is also legitimate.

\def\rhoR{\rho_\text{\sc r}}

\begin{remark}[{\sl Barotropic viscoelastic fluids}]\label{rem-barotropic-fluids}\upshape
Fluids are mechanically characterized by $\varphi(\FF)=\breve\varphi(\det\FF)$ for some
$\breve\varphi=\breve\varphi(J)$ and the Cauchy stress \eq{Cauchy-stress-polyconvex}
reduces to the hydro-static part of $\TT=-p\bbI$ with the pressure
$p=-J\breve\varphi'(J)-\breve\varphi(J)
=-(J\breve\varphi(J))'$. Taking into account $\varrho=\rhoR/J$ with
$\rhoR$ denoting the referential mass density reveals the so-called
{\it isentropic state equation}\index{state equation!isentropic} in the form
$p=-[(J\breve\varphi(J))'](\rhoR/\varrho)$. For the referential
$\breve\upvarphi(J):=J\breve\varphi(J)$, we have simply
\begin{align}\nonumber
p=-\breve\upvarphi'\Big(\frac{\rhoR}\varrho\Big)\,.
\end{align}
Such fluids are called
{\it barotropic},\index{fluid!barotropic} i.e.\ their density depends only on pressure
or, vice versa, the pressure depends on density. For example, choosing
$\breve\varphi(J)=J^{-\gamma}$ yields $(J\varphi(J))'=
(J^{1-\gamma})'=(1{-}\gamma)J^{-\gamma}$, which further gives
$p=-(1{-}\gamma)(\rhoR/\varrho)^{-\gamma}=a\varrho^\gamma$
with $a=(\gamma{-}1)/\rhoR^{\gamma}$. For this case,
without the hyperviscosity, the fully implicit time-discretization method
was used in \cite{Kar13CFEM,FeKaPo16MTCV} with
$\gamma>3$ and in \cite{Zato12ASCN} with $\gamma>2$. Actually, the critial
exponent $\gamma$ for other less constructive methods is $3/2$. Here, involving
hyperviscosity into our model needs  merely the convexity of
$\breve\varphi(J)=J^{-\gamma}$, which also holds for $\gamma>0$ 
although $0<\gamma\le1$ leads to a nonphysical coefficient $a\le0$.
\end{remark}

\subsection{Some generalizations of polyconvex stored energies and examples}\label{sec-rem-exa}

The polyconvex ansatz for the actual stored energy $\varphi$ in Section~\ref{sec-polycon}
served well for the lucid presentation of the method, but in fact has rather limited
applications. For this reason, we end this paper by few remarks widening the
applicability and accompany them by representative examples. 

\begin{remark}[{\sl General polyconvex\index{stored energy!polyconvex} stored
energy}]\label{rem-polyconvex-general}\upshape
Denoting $\HH:={\rm Cof}\FF$, a general polyconvex stored energy $\varphi$
uses some convex 
$\breve\varphi:\R^{3\times3}\times\R^{3\times3}\times\R\to\R\cup\{+\infty\}$
so that $\varphi(\FF)=\breve\varphi(\FF,\HH,J)$. The calculations
\eq{Cauchy-Mandel-stress-polyconvex} leading to the Cauchy and the Mandel
stresses must be extended as
\begin{subequations}\label{Cauchy-Mandel-stress-polyconvex+}\begin{align}\nonumber
\TT&=\varphi_\FF'(\FF)\FF^\top\!\!+\varphi(\FF)\bbI
\\&\nonumber
=\breve\varphi_\FF'(\FF,\HH,J)\FF^\top\!\!+\breve\varphi_\HH'(\FF,\HH,J)({\rm Cof}'\FF)\FF^\top\!\!
+\breve\varphi_J'(\FF,\HH,J)({\rm Cof}\FF)\FF^\top\!+\breve\varphi(\FF,\HH,J)\bbI
\\[-.2em]&=\breve\varphi_\FF'(\FF,\HH,J)\FF^\top\!\!
+\breve\varphi_\HH'(\FF,\HH,J)({\rm Cof}'\FF)\FF^\top\!\!
+\Big(J\breve\varphi_J'(\FF,\HH,J){+}\breve\varphi(\FF,\HH,J)\Big)\bbI\ \text{ and}
\label{Cauchy-stress-polyconvex+}
\\[-.2em]\MM&={\rm dev}\big(\FF^\top\!\varphi_\FF'(\FF)\big)={\rm dev}\Big(
\FF^\top\!\breve\varphi_\FF'(\FF,\HH,J)
+\FF^\top\!\breve\varphi_\HH'(\FF,\HH,J){\rm Cof}'\FF\Big)\,,
\label{Mandel-stress-polyconvex+}\end{align}\end{subequations}
where ${\rm Cof}'=\det{}''$ is a linear mapping $\R^{3\times3}\to\R^{3\times3}$
being the Hessian of the cubic form $\det:\R^{3\times3}\to\R$; actually
the rather short-hand formulas \eq{Cauchy-Mandel-stress-polyconvex+} should
be written more exactly componentwise. The flow rule for
$\HH={\rm Cof}\FF=J\FF^{-\top}$ can be revealed from the calculation
\begin{align}\nonumber
\DT\HH=\DT{\overline{J\FF^{-\top}}}\!
&=\DT J\FF^{-\top}\!\!+J\DT{\overline{\FF^{-\top}}}
=\DT J\FF^{-\top}\!\!+J\DT{\overline{(\HH/J)}}
\\[-.2em]&\nonumber
=\DT J\FF^{-\top}\!\!+({\rm Cof}'\FF)\DT\FF-(\DT J/J)\HH
\\&=\DT J\FF^{-\top}\!\!+({\rm Cof}'\FF)\DT\FF-\DT J\FF^{-\top}\!
=({\rm Cof}'\FF)\Big((\nabla\vv)\FF-\FF\Lp\Big)\,.
\label{H-flow-rule}\end{align}
This is to be included, in the corresponding time-discrete $\DELTA$-regularized form,
into \eq{Euler-large-viscoelastodyn+disc} as
\begin{align}\label{Euler-large-viscoelastodyn-H}
\frac{\HH_\etau^k{-}\HH_\etau^{k-1}\!\!}\tau\,=({\rm Cof}'\FF_{\!\etau}^k)
\Big((\nabla\vv_\etau^k)\FF_{\!\etau}^k-\FF_{\!\etau}^k\Lp_\etau^k\Big)
-(\vv_\etau^k\Cdot\nabla)\HH_\etau^k+\DELTA\Delta\HH_\etau^k
\end{align}
while the initial condition $\HH_\etau^0={\rm Cof}\,\FF_0$ is to be included into
\eq{IC-large-disc}. Analytically, making the regularization like
\eq{varphi-regularized} also in terms of $\HH$, we have the the coercivity of
$\breve\varphi_\etau$ also in terms of $\HH$ and can rely on the estimate
like \eq{Euler-large-est3} also for $\overlineHetau$, namely
\begin{align}
&\label{Euler-large-est3-H}
\|\overlineHetau\|_{L^\infty(I;L^2(\varOmega;\R^{3\times3}))}^{}\le C \ \ \,\text{ and }
\ \ \ \,
\|\nabla\overlineHetau\|_{L^2(I{\times}\varOmega;\R^{3\times3\times3})}^{}\le\frac{C}{\DELTA^{3/4}}\,.
\end{align}
From \eq{Euler-large-est3}, we have the bound of ${\rm Cof}'\overlineFetau$ in
$L^\infty(I;L^2(\varOmega;\R^{3\times3}))\,\cap\,L^2(I;L^6(\varOmega;\R^{3\times3}))
\subset L^{10/3}(I{\times}\varOmega;\R^{3\times3})$. Thus, when the bound for
$\nabla\overlinevvtau\in L^p(I;W^{1,\infty}(\varOmega;\R^3))$ is
proved as in \eq{Euler-small-est2}, we have
$({\rm Cof}'\overlineFetau)((\nabla\overlinevvtau)\overlineFetau{-}\overlineFetau\overlineLpetau)$ is bounded in
$L^p(I;L^1(\varOmega;\R^{3\times3}))\cap L^{15/14}(I;L^{5/3};\R^{3\times3}))$, which
eventually yields also a bound for $\pdt{}\HH_\etau$.
\end{remark}

\begin{example}[{\sl A polyconvex\index{stored energy!polyconvex} stored energy}]\label{exa-polyconvex-general}\upshape
A popular model is 
the {\it neo-Hookean} 
stored energy
\begin{align}
\varphi(\FF)&
=\frac12K_\text{\sc e}^{}\Big(\det\FF{-}1{+}\frac{\epsilon}{K_\text{\sc e}^{}}\Big)^2
+G_\text{\sc e}^{}
\frac{{\rm tr}(\FF\FF^\top)}{(\det\FF)^{2/3}}+\begin{cases}\epsilon{\rm ln}(1/\!\det\FF)&\text{for }\ \det\FF>0,\\\ \ \infty&\text{otherwise }\,. \end{cases}
\label{neo-Hookean-potential-actual}\end{align}
For any $\epsilon>0$ typically ``small'' (in Pa\,=\,J/m$^3$), this ansatz
exhibits the physically relevant blow-up
$\lim_{{\rm det}\FF\to0+}^{}\varphi(\FF)=+\infty$ and simultaneously attains its
minimum on the orbit 
${\rm SO}_3=\{Q\in\R^{3\times3};\ Q^\top\!Q=QQ^\top\!=\bbI\}$, i.e.\ the
so-called special orthogonal group. Also, this function obviously has the form
$\varphi(\FF)=\breve\varphi(\FF,\det\FF)$ with $\breve\varphi(\FF,J)
=\frac12K_\text{\sc e}^{}(J{-}1{+}\epsilon/K_\text{\sc e}^{})^2+
G_\text{\sc e}^{}J^{-2/3}|\FF|^2-\epsilon{\rm ln}J$ for $J>0$, otherwise
$\breve\varphi(\FF,J)=+\infty$. Such $\breve\varphi$ is
convex. Indeed, realizing that that both $J\mapsto(J{-}1)^2$
and $J\mapsto-{\rm ln}J$ are certainly convex, the proof of the convexity of
$\breve\varphi$ is based on the analysis of the 2nd-order derivative (Hessian) of
the function $\R^{3\times3}{\times}\R^+\to\R:(\FF,J)\mapsto{\rm tr}(\FF\FF^\top)/J^p=
|\FF|^2/J^p$ which is the $(3^2{+}1)\times(3^2{+}1)$-matrix
$$
\bigg(\frac{|\FF|^2}{J^p}\bigg)''=
\bigg(\begin{array}{c}2\FF/J^p\\-p|\FF|^2/J^{p+1}\\\end{array}\bigg)'
=\bigg(\begin{array}{cc}2J^{-p}\bbI{\times}\bbI\ \ \ , &\ -2p\FF/J^{p+1}
\\-2p\FF/J^{p+1}\,, &\ p(p{+}1)|\FF|^2/J^{p+2}\end{array}\bigg)\,.
$$
The desired convexity needs the positive semi-definiteness of this matrix.
Since surely the matrix $J^{-p}\bbI{\times}\bbI$ is positive definite,
by Sylvester's criterion, it suffices to verify that the
determinant of this $(3^2{+}1)\times(3^2{+}1)$-matrix is
non-negative. By the formula for the ``block determinant'', this overall
determinant equals
\begin{align}\nonumber
&\!\!\det\Bigg(\bigg(\frac{|\FF|^2}{J^p}\bigg)''\Bigg)=
\det\Big(2\frac{\bbI{\times}\bbI}{J^{p}}\Big)\:
\det\Big(\frac{p(p{+}1)|\FF|^2}{J^{p+2}}
-\frac{2p\FF^\top}{J^{p+1}}[2J^{-p}\bbI{\times}\bbI)]^{-1}
\frac{2p\FF}{J^{p+1}}\Big)
\\[-.3em]&\hspace{3.em}
=\big(2^{9}J^{-9p}\big)\Big(p(p{+}1)J^{-p-2}|\FF|^2\!-2p^2J^{-p-2}\FF^\top\!\Colon\FF\Big)
=2^9p(1{-}p)J^{-10p-2}|\FF|^2\ge0;\!
\label{convexity-neo-Hookean}\end{align}
note that the desired non-negativity holds for $0\le p\le 1$. For 
\eq{neo-Hookean-potential-actual}, the relevant choice is $p=2/3$.
For a similar proof 
see also \cite[Prop.\,6]{CDHL88ETSC} or \cite[Lemma\,2.1]{HarNef03PGPH}. 
For many other polyconvex functions relevant in large-strain elasticity
(as well as some counterexamples) we refer to \cite{HarNef03PGPH} and the references
therein. In particular,
the last term in \eq{neo-Hookean-potential-actual} admits various
modifications, e.g.\ by replacing $-\epsilon{\rm ln} J$ with 
$+\epsilon/J^\varkappa$ with $\varkappa\ge1$.
\end{example}

\begin{remark}[{\sl Referential stored energies}]\label{rem-referential}\upshape
In literature, particular models are formulated rather in terms of the referential
stored energy $\upvarphi(\FF):=\varphi(\FF)\det\FF$ than the actual stored energy
$\varphi$ itself. Then the Cauchy and the Mandel stresses take the form
\begin{subequations}\label{Cauchy-Mandel-stress++}\begin{align}\nonumber
\TT&=\varphi'(\FF)\FF^\top\!\!+\varphi(\FF)\bbI
=\Big(\frac{\upvarphi(\FF)}{\det\FF}\Big)'\FF^\top\!\!+\frac{\upvarphi(\FF)}{\det\FF}\bbI
\\&
=\frac{\upvarphi'(\FF)\FF^\top\!\!}{\det\FF}
-\frac{\upvarphi(\FF)({\rm Cof}\FF)\FF^\top\!\!}{\det\FF^2}
+\frac{\upvarphi(\FF)}{\det\FF}\bbI=\frac{\upvarphi'(\FF)\FF^\top\!\!}{\det\FF}
\ \ \text{ and}
\label{Cauchy-stress++}
\\[-.0em]\MM&={\rm dev}\big(\FF^\top\!\varphi'(\FF)\big)={\rm dev}\Big(
\FF^\top\!\Big(\frac{\upvarphi(\FF)}{\det\FF}\Big)'\Big)
={\rm dev}\Big(\frac{\FF^\top\!\upvarphi'(\FF)}{\det\FF}\Big)\,.
\label{Mandel-stress++}\end{align}\end{subequations}
\end{remark}

\begin{remark}[{\sl Behind polyconvexity}]\label{rem-behind-polyconvexity}\upshape
The polyconvexity is a concept for static large-strain elasticity
in the reference configuration, i.e.\ in the Lagrangian frame, and is
not directly relevant in the Eulerian elastodynamics. In fact, the convexity of
$\breve\varphi=\breve\varphi(\FF,J)$ can be considered in a different $\FF$ than in
Section~\ref{sec-polycon}-\ref{sec-large-proof} and in the previous
Remark~\ref{rem-polyconvex-general}. E.g., one can think about e.g.\
$\HH:=J^\alpha\FF$ to be understood as a ``re-calibrated'' deformation gradient.
An analogy of \eq{H-flow-rule} reads as
\begin{align}\nonumber
\DT\HH=\DT{\overline{J^\alpha\FF}}\!
&=\alpha J^{\alpha-1}\DT J\FF+J^\alpha\DT\FF
\\&=\alpha J^{\alpha-1}({\rm div}\,\vv)J\FF+J^\alpha\big((\nabla\vv)\FF-\FF\Lp\big)
=\alpha({\rm div}\,\vv)\HH+(\nabla\vv)\HH-\HH\Lp\,.
\end{align}
Instead of \eq{Euler-large-viscoelastodyn-H}, one should expand the system
\eq{Euler-large-viscoelastodyn+disc} by
\begin{align}\label{Euler-large-viscoelastodyn-H+}
\frac{\HH_\etau^k{-}\HH_\etau^{k-1}\!\!}\tau\,=\alpha({\rm div}\,\vvk)\HH_\etau^k+
(\nabla\vvk)\HH_\etau^k-\HH_\etau^k\Lp_\etau^k-(\vv_\etau^k\Cdot\nabla)\HH_\etau^k
+\DELTA\Delta\HH_\etau^k
\end{align}
while the initial condition $\HH_\etau^0=(\det\,\FF_0)^\alpha\FF_0$ is to be
included into \eq{IC-large-disc}, relying on $\det\,\FF_0>0$ so that
$(\det\,\FF_0)^\alpha$ is well defined. The actual stored energy
$\varphi(\FF)=\breve\varphi(\HH,J)$ or even $\breve\varphi(\FF,\HH,J)$ with
$\breve\varphi$ convex can then be treated as before. Such $\varphi$ does not need
to be polyconvex, cf.\ Example~\ref{exa-neo-Hookean-ref} below with the calculations
like \eq{Cauchy-Mandel-stress++} modified accordingly. In general, one can devise such
scheme for ``polynomially polyconvex'' stored energies, i.e.\ such $\varphi$'s that
admit a representation $\breve\varphi$ which is convex in terms of various polynomial
functions of minors of $\FF$. 
\end{remark}

\begin{example}[{\sl A referential neo-Hookean ansatz}]\label{exa-neo-Hookean-ref}\upshape
Following the previous remark, also the neo-Hookean ansatz
\eq{neo-Hookean-potential-actual} is standardly understood rather in the referential
frame than in the actual frame. The above Example~\ref{exa-polyconvex-general} can
therefore be considered relevant rather only for slightly compressible media only,
if at all. Interestingly, the polyconvex referential stored energy $\upvarphi$ having
the form \eq{neo-Hookean-potential-actual},
in fact, does not lead to a polyconvex actual energy $\varphi=\upvarphi/J$. This
is seen from the calculations \eq{convexity-neo-Hookean}, which shows non-convexity
for $p=1+2/3$ and, furthermore, the last term in \eq{neo-Hookean-potential-actual}
leads to the contribution to the actual stored energy as $J\mapsto-({\rm ln}J)/J$
which is also nonconvex.
{Anyhow, we can cast the evolution system directly for the new variable
$K:=({\rm ln}\,J)/J=A\,{\rm ln}\,J$ with $A=1/J$, specifically
$\DT K=({\rm div}\,\vv)(A{-}K)$ with $\DT A=-({\rm div}\,\vv)A$ as in
\eq{DT-det-inverse} and then the contribution to the actual stored energy
\eq{neo-Hookean-potential-actual} is just linear, i.e. $-\epsilon K$.}
Alternatively, the variant $1/J^\varkappa$ instead of $-{\rm ln} J$ mentioned in the
previous Example~\ref{exa-polyconvex-general}, now considered in the referential frame,
would give a convex contribution $1/J^{\varkappa+1}$ in the actual
$\varphi$ if $\varkappa\ge0$. The mentioned resulting non-polyconvex term
$J^{-1-2/3}|\FF|^2$ in the actual stored energy can be handled by
Remark~\ref{rem-behind-polyconvexity} when considering $\alpha=-5/6$, i.e.\
$\HH:=J^{-5/6}\FF$, so that
\begin{align}\nonumber
\varphi(\FF)=\frac{\upvarphi(\FF)}{\det\FF}=
\breve\varphi\Big(\frac{\FF}{\det\FF^{5/6}},\det\FF\Big)
\ \ \text{ with }\ \ 
\breve\varphi(\HH,J)=\frac12K_\text{\sc e}^{}(J{-}2{+}1/J)+G_\text{\sc e}^{}|\HH|^2\,;
\end{align}
such $\breve\varphi$ is obviously convex, even strongly in terms of $\HH$,
which allows for a simpler regularization in \eq{varphi-regularized},
cf.\ Remark~\ref{rem-neo-Hookean-ref}. Alternatively, we could put
$\HH:=J^{-1/3}\FF$, i.e.\ the isochoric part of $\FF$, and then consider
$\breve\varphi(\HH,J)=\frac12K_\text{\sc e}^{}(J{-}2{+}1/J)
+G_\text{\sc e}^{}|\HH|^2/J$; such $\breve\varphi$ is convex, although
not strongly.
\end{example}

\begin{example}[{\sl A referential Mooney–Rivlin ansatz}]\label{exa-Mooney–Rivlin-ref}
\upshape
The compressible Mooney–Rivlin model, which works more realistically for large elastic
strains than the neo-Hookean model, expands \eq{neo-Hookean-potential-actual} by the
term $G_\text{\sc mr}|{\rm Cof}\FF|^2/J^{4/3}$ in the referential frame. However, this
term is not convex, see \cite[Sect.\,2.2]{CDHL88ETSC}. This means
that \eq{neo-Hookean-potential-actual}, beside $\Ge|\FF_{\rm iso}|^2$ with the
isochoric part $\FF_{\rm iso}=\FF/J^{1/3}$ of the elastic distortion $\FF$, is expanded by
the term $G_\text{\sc mr}|\FF_{\rm iso}^{-1}|^2=G_\text{\sc mr}J^{2/3}\FF^{-1}\FF^{-\top}$. In
the actual Eulerian frame, this means
$G_\text{\sc mr}|\FF_{\rm iso}^{-1}|^2/J=G_\text{\sc mr}J^{-1/3}\FF^{-1}\FF^{-\top}
=G_\text{\sc mr}J^{-1/3}|{\rm Cof}\FF/J|^2=G_\text{\sc mr}|\HH|^2$ with
$\HH=J^{-7/6}{\rm Cof}\FF$. 
For such $\HH$, the flow rule \eq{H-flow-rule} is then modified as
\begin{align}\nonumber
\DT\HH&=\DT{\overline{J^{-7/6}{\rm Cof}\FF}}\!
=J^{-7/6}\DT{\overline{{\rm Cof}\FF}}-\frac76J^{-13/6}\DT J{\rm Cof}\FF
=J^{-7/6}({\rm Cof}'\FF)\DT\FF-\frac76J^{-13/6}\DT J{\rm Cof}\FF
\\[-.1em]&\nonumber
=\frac{{\rm Cof}'\FF\!}{J^{7/6}}\Big((\nabla\vv)\FF{-}\FF\Lp\Big)
-\frac{7{\rm div}\vv\!}{6J^{7/6}\!}\,{\rm Cof}\FF
=({\rm Cof}'\FF)\Big((\nabla\vv)\HH{-}\HH\Lp\Big)-\frac{7{\rm div}\vv}6\HH.
\end{align}
This discrete kinematics \eq{Euler-large-viscoelastodyn-H+} is then modified as
\begin{align}\nonumber
\frac{\HH_\etau^k{-}\HH_\etau^{k-1}\!\!}\tau\,=
({\rm Cof}'\FF_\etau^k)\Big((\nabla\vvk)\HH_\etau^k{-}\HH_\etau^k\Lp_\etau^k\Big)
-\frac{\!7{\rm div}\,\vvk\!}6\HH_\etau^k\!-(\vv_\etau^k\Cdot\nabla)\HH_\etau^k\!
+\DELTA\Delta\HH_\etau^k.
\end{align}
The discrete kinematics for $\HH$, when accompanied by the initial condition
$\HH_\etau^0=J_0^{-7/6}{\rm Cof}\FF_0$, should then expand the system
\eq{Euler-large-viscoelastodyn+disc} with $\mathscr{T}_\EEps$ and 
$\mathscr{M}_\EEps$ modified accordingly.
\end{example}

\begin{remark}[{\sl A simplified approximation $\breve\varphi_\EEps$}]\label{rem-neo-Hookean-ref}\upshape
As also seen above, $\varphi$ typically has a specific form of the type
$\varphi(\FF)=f(\FF,\HH,J)+s(J)$ for some $\HH$ (such as
$\HH=J^{-5/6}\FF$ as in Example~\ref{exa-neo-Hookean-ref} or
$\HH=J^{-7/6}{\rm Cof}\FF$ as in Example~\ref{exa-Mooney–Rivlin-ref}) and for
$J=\det\FF$ with some $f:\R^{3\times3}\times\R^{3\times3}\times\R$ convex with
a polynomial growth and with some $s:\R^+\to\R\cup\{+\infty\}$ convex
with a singularity for $J\to0+$ like we have in Example~\ref{exa-neo-Hookean-ref}.
This allows for a simpler regularization $\breve\varphi_\EEps$ than that one which
was devised in \eq{varphi-regularized}. Actually, instead of the Yosida
approximation $\mathscr{Y}_\EPS\breve\varphi$ used in \eq{varphi-regularized}, we
can simply keep the original $f$ and only modify the singular part $s$ for arguments
below $J_{\min{}}$ from \eq{J-min} e.g.\ by defining a $C^1$ function
$s_\EPS(J)=s_\EPS(J_{\min{}})+s_\EPS'(J_{\min{}})(J{-}J_{\min{}})$ for $J<J_{\min{}}$
while keeping $s_\EPS(J)=s(J)$ for $J\ge J_{\min{}}$. Then the convergence in Step~6 in
Section~\ref{sec-large-proof} can completely ignore such (actually inactive)
modification and would thus be as simple as in Section~\ref{sec-linearized-proof}.
In addition, the $\DELTA$-regularization in \eq{varphi-regularized} can be omitted
in those variables at which $\breve\varphi$ is strongly convex. It occurs specifically
in $\HH$ in Examples~\ref{exa-neo-Hookean-ref} and \ref{exa-Mooney–Rivlin-ref}.
\end{remark}

{
\section*{{\large Acknowledgments}}

\vspace*{-1em}

This research has been partially supported also from the CSF (Czech Science
Foundation) project GA22-00863K
and by the institutional support RVO: 61388998 (\v CR).

\medskip
}

\baselineskip12pt


\end{document}